\documentclass[12pt,eqno]{article}
\usepackage{amsmath}
\usepackage{amssymb}

\addtocounter{section}{-0}
\newtheorem{guess}{Theorem}[section]
\newtheorem{proposition}[guess]{Proposition}

\newtheorem{remark}[guess]{Remark}
\newtheorem{lemma}[guess]{Lemma}
\newtheorem{corollary}[guess]{Corollary}
\newtheorem{theorem}[guess]{Theorem}

\numberwithin{equation}{section}

\newcommand{\be}{\begin{equation}}
\newcommand{\ee}{\end{equation}}

\newcommand{\proof}{{\it Proof.\ }}
\newcommand{\qed}{\hspace*{\fill }$\square $}

\newcommand{\Ninf}{$N\to\infty$}

\def\I{{\cal I}}

\begin{document}

\title{Dynamics of periodic Toda chains with a large number of particles}

\author{D. Bambusi\footnote{Dipartimento di Matematica, Universit\`a degli Studi di
Milano, Via Saldini 50, I-20133 Milano}, T. Kappeler\footnote{Institut f\"ur Mathematik,
Universit\"at Z\"urich, Winterthurerstrasse 190, CH-8057 Z\"urich. Supported in part by
the Swiss National Science Foundation.}, T. Paul\footnote{CNRS
and CMLS, \'Ecole Polytechnique, F-91128 Palaiseau  }}

\maketitle

\date{}

\begin{abstract}
\noindent For periodic Toda chains with a large number $N$ of particles
we consider states which are $N^{-2}-$close to the equilibrium and constructed 
by discretizing any given $C^2-$functions with mesh size $N^{-1}$.
For such states we derive asymptotic expansions of the Toda frequencies
$(\omega^N_n)_{0 < n < N}$ and the actions $(I^N_n)_{0 < n < N},$ both listed in the standard way, in powers of $N^{-1}$ as $N \to \infty$. 
%listed in accordance with the ordering of the frequencies at the equilibrium, 
%$(2 \sin \frac{n\pi } {N})_{0 < n < N}$. 
At the two edges $n \sim 1$ and $N -n \sim 1$, the expansions of the frequencies are computed up to order $N^{-3}$ with an error term of higher order. 
Specifically, the coefficients of the expansions of $\omega^N_n$ and $\omega^N_{N-n}$  at order $N^{-3}$ 
are given by a constant multiple of the n'th KdV frequencies $\omega^-_n$ and $\omega^+_n$
of two periodic potentials, $q_{-}$ respectively $q_+$, constructed in terms of the states considered.
The frequencies $\omega^N_n$ for $n$ away from the edges are shown to be asymptotically close to the frequencies of the equilibrium.
For the actions $(I^N_n)_{0 < n < N},$ asymptotics of a similar nature are derived.

\end{abstract}

%\tableofcontents

%%%%%%%%%%%%%%%%%%%%%%%%%%%%%%%%%%%%%%%%%%%%%%%%%%%%%%%%%%%%%%%%%%%%%%%%%%%%%%%%%

\section{Introduction}\label{intro}

%{\tt to be inserted somewhere the definition of O and that all the
%  estimates are uniform in bounded sets of $\alpha$ and $\beta$. To be
%inserted also that $L:=F(M)\equiv M'$.}

In this paper we study the asymptotics of the dynamics of periodic Toda chains with
a large number of particles of equal mass for initial data close to the equilibrium. 
Toda chains were introduced
by Toda \cite{T} as a class of special Fermi Pasta Ulam (FPU) chains with the main
feature that they are {\em integrable} Hamiltonian systems.
The Hamiltonian of the Toda chain with $N$ particles is given by
\[
\mathcal H^{\mbox{}}=\frac 1 2\sum_{n=1}^Np_n^2+\sum_{n=1}^Ne^{q_n-q_{n+1}}
\]
where $q_n$ denotes the displacement of the $n$'th particle from its equilibrium position
and $p_n$ its momentum. It is convenient to define $(q_n,p_n)$ for any $n\in\mathbb Z$
by requiring that $(q_{n+N},p_{n+N}) = (q_n,p_n) \,\, \forall \, n\in\mathbb Z$.
Note that the total momentum $\sum_{i=1}^N p_i$ is a conserved quantity of the Toda flow
and hence the Hamiltonian equations of motion
imply that the center of mass $\frac{1}{N}\sum_{n=1}^N q_n$ moves at constant speed. Hence when considered relative to the motion of its center of mass, the Toda chain has $N-1$ degrees of freedom.\\
Flaschka \cite{F} introduced the variables 
$b = (b_n)_{1 \leq n \leq N}$, $a = (a_n)_{1 \leq n \leq N}$ defined by
$ b_n=-p_n \in \mathbb R$, $a_n=e^{\frac 12(q_n-q_{n+1})} >0$
and showed that when evolved along the Toda flow the corresponding equations of motion 
of $b$, $a$ can be described by
a Lax pair $(L, B)$,
\begin{equation}\label{laxx1}
\dot L = [B,L]
\end{equation}
where the $N \times N$ matrices $L = L(b,a)$ and $B = B(a)$ are given by
\[ \begin{pmatrix}
b_1 &a_1 &0&\ldots&0&a_N \\
a_1 &b_2 &a_2&\ldots &0&0 \\
0 &a_2 &b_3&\ldots&0&0\\
\vdots&&&&&\vdots\\
0&\ldots &&\ldots &b_{N-1}&a_{N-1}\\
a_N &0&\vdots &&a_{N-1}&b_N
\end{pmatrix}
\mbox{ and }
\begin{pmatrix}
0 &a_1 &0&\ldots&0&-a_N \\
-a_1 &0 &a_2&\ldots &0&0 \\
0 &-a_2 &0&\ldots&0&0\\
\vdots&&&&&\vdots\\
0&\ldots &&\ldots &0&a_{N-1}\\
a_N &0&\vdots &&-a_{N-1}&0
\end{pmatrix}
\]
respectively. Note that by the periodicity of the $q_n$'s, $\prod_{n=1}^N a_n= 1$
and hence the position of the center of mass cannot be recovered from the $a_n$'s.
However the latter is easily reconstructed from its initial condition and the total momentum
$-\sum_{n=1}^N b_n$. 
In the sequel, we exclusively concentrate on the system described by \eqref{laxx1}, viewing it
as equivalent to the Toda chain and referring to it by the same name. Actually, it is convenient
to choose a slightly larger space for the variable $a$ by requiring only that $a= (a_n)_{1 \le n \le N} \in \mathbb R^N_{> 0}.$

First we note that $b,$ $a$ are not canonical variables. In fact, when expressing the Hamiltonian
equations of motion of the Toda chain in terms of the variables $b,$ $a,$
the corresponding Poisson bracket is degenerate.  The level sets 
\[
\mathcal L_{ c_{1}, c_{2}} := \{(b,a) \in \mathbb R^N \times \mathbb R^N_{>0} :  \,\,
\sum_{n=1}^N b_n = c_1; \,\,  \prod_{n=1}^N a_n= c_2 \}
\] 
are the corresponding symplectic leaves. As already mentioned, when restricted to such a leaf, the system has $N-1$ degrees of freedom.
Using \eqref{laxx1}, it can be shown to be an integrable system of $N-1$ coupled oscillators, 
meaning that it admits globally defined Birkhoff coordinates  -- see \cite{HK} for details. In particular, its invariant manifolds are smooth tori. As a consequence, the dynamics
of such chains is quasi-periodic in time.
Furthermore, in terms of the spectrum of $L(b,a)$, which by \eqref{laxx1} is preserved by the Toda flow,  there is a canonical choice of $N-1$ globally defined actions
$I = (I_n)_{ 0 < n < N} \in \mathbb R_{\ge 0}^{N-1}$
-- see \cite{HK} for details as well as Section \ref{actions}. They parametrize the invariant tori 
$\mathcal T_{I}$ and allow to read off their dimension,  $\mbox{dim } \mathcal T_{I} = |\{0 < n < N | I_n > 0\}|$. Furthermore, the Toda Hamiltonian $\mathcal H$ can be expressed as a function of $I$ which we again denote by  $\mathcal H$. In fact, $\mathcal H$ is a real analytic function of $I$ and the total momentum. 
The frequencies corresponding to this canonical choice of actions are denoted by
$\omega_n(I) = \partial_{I_{n}} \mathcal H $. We refer to 
$\omega (I) = (\omega_n(I))_{0< n < N}$ as the frequency vector corresponding to $I$.
As an aside we mention that actions are uniquely determined only up to unimodular transformations. Nevertheless, the dynamics on such a torus can be described in a coordinate free way in terms of the frequency map.
More precisely, let $x= (b,a)$ be an arbitrary point on a leaf $\mathcal L _{ c_{1}, c_{2}}$
and denote by $\mathcal T(x)$ the invariant torus containing $x$. Let  $\mathcal M \equiv \mathcal M (x)$ be the frequency module of $\mathcal T(x)$,
consisting of all integer combinations of $\omega_1 , \dots , \omega_{N-1},$
\[
\mathcal M = \{ k \cdot \omega : k \in \mathbb Z^{N-1} \} \subset \mathbb R.
\]
Note that $\mathcal M$ doesn't depend on the choice of the actions and thus is invariantly defined.
 It plays an important role when analyzing data of the evolution of systems such as Toda chains 
-- see e.g. \cite{L}.
 For any given complex valued function $f$ defined on $\mathcal L_{ c_{1}, c_{2}}$ and any
$x \in \mathcal L_{ c_{1}, c_{2}}$
let $f_x(t) := f (\Phi^t(x))$ where $x \mapsto \Phi^t(x)$ denotes the Toda flow. 
The support of the Fourier transform of $f_x(t)$ with respect to the real variable $t$ is then
contained in the frequency module $\mathcal M$ and if $f$ is sufficiently regular,  $f_x(t)$ takes the form
\[
f_x(t) = \sum_{k \in \mathbb Z^{N-1}} a_x(k) \exp(i \varphi_x(k) t)
\]
where $\varphi_x: \mathbb Z^{N-1} \to \mathcal M$ is a linear map and 
$a_x: \mathbb Z^{N-1} \to \mathbb C$ has appropriate decay conditions.  In the case at hand, $\varphi_x$ can be chosen
to be $\varphi_x(k) = k \cdot \omega(I)$ where $I$ is the action of $\mathcal T(x)$. For more details see e.g. \cite{KP}, Section 3.\\
The aim of this paper is to compute the asymptotics as $N \to \infty$ of the
frequencies and actions of states which are $N^{-2}$-close
to the equilibrium $0_N= (0, ... , 0), \,\, 1_N = (1, ...,1).$
As KdV frequencies will play an important role for describing these asymptotics we first need to review some features of the KdV equation in the periodic set-up. It turns out that for our purposes, the period of the space variable is equal to $\frac{1}{2}$.
Recall that the (generalized) KdV equation with parameters $d_1, d_2 \in \mathbb R$
\[
\partial_t q = d_1 ( - \partial^3_x q  + 6q\partial_x q) + d_2 \partial_x q \  ,  
\qquad t \in \mathbb R, \,\, x \in \mathbb R/ (\mathbb Z /2) 
\]
is an integrable Hamiltonian PDE with Hamiltonian 
\[
{\mathcal H}_{d_1, d_2} = d_1  \int ^{\frac{1}{2}}_0 \big( \frac{1}{2}(\partial _x q)^2 + q^3 \big) dx
+ d_2 \frac{1}{2} \int^\frac{1}{2}_0 q^2dx
\]
Note that given any solution,  the average is a conserved quantity. For our purposes it suffices to consider solutions in Sobolev spaces $H_0^s, s \ge 2,$ of elements in $H^s$ with mean zero.
Such solutions exist globally in time and are almost periodic,
evolving on invariant sets which are tori, generically of infinite dimension. These invariant tori are parametrized by globally defined actions which have the property that on $H_0^1,$ the KdV Hamiltonian
can be expressed as a real analytic function of them alone. Furthermore, the KdV frequencies are given by the partial derivatives of the KdV Hamiltonian with respect to these actions. 
See \cite{KP} for details. 
Note that ${\mathcal H}_{d_1, d_2}$ is a linear combination of the first two Hamiltonians
${\mathcal H}_1$ and ${\mathcal H}_2$ of the so called KdV hierarchy, where 
\[
{\mathcal H}_1(q) = \frac{1}{2} \int ^{\frac{1}{2}}_0 q^2dx \quad \mbox{and} \quad 
{\mathcal H}_2 (q) = \int ^{\frac{1}{2}}_0 \big( \frac{1}{2}(\partial _x q)^2 + q^3 \big) dx.
\]

To state our first result, assume that $\beta,$ $\alpha$ are in the space 
$C^2_0(\mathbb T) \equiv C_0^2(\mathbb T, \mathbb R)$ of one periodic, real valued functions of class $ C^2$ and mean $0$. The space $C_0^2(\mathbb T)$ is endowed with the standard supremum norm $\| f \|_{C^2}$.\\
For any $N \ge 3,$ we then introduce the periodic Toda chain with $N$ particles, defined in terms of Flaschka coordinates by
\be\label{1Da.300} b^N_n = \frac1{4N^2} \beta \big( \frac{n}{N} \big) \quad \mbox{and} \quad a^N_n =
      1 + \frac1{4N^2} \alpha \big( \frac{n}{N} \big).
 \ee
Alternatively, one can consider
\[
p^N_n = - \frac1{4N^2} \beta \big( \frac{n}{N} \big) \quad \mbox{and} \quad 
q^N_n = - \frac{2}{4N} \xi \big( \frac{n}{N} \big)
\]
where $\xi$ is the element in $C_0^3(\mathbb T)$, satisfying $\xi ' = \alpha$.
Using that
\[
\exp( \frac{q^N_{n} - q^N_{n+1}}{2} ) = 1 + \frac{q^N_{n} - q^N_{n+1}}{2} +O(N^{-4}) 
= a_n^N + O(N^{-3})
\]
one can show that our results stated below hold for either of the two discretizations.
See the remark at the end of Section \ref{spec} for details.
In this paper, we concentrate on the case where the data is given in the form \eqref{1Da.300}\\
Denote by $\omega^N_n$ and  $I^N_n,\ 0 < n < N,$ the frequencies and actions, corresponding
to $(b^N, a^N)$.
To describe their asymptotics as $N \to \infty$ 
introduce the two potentials 
\[q_\pm(x) = -2\alpha(2x)\mp\beta(2x). 
\]
Note that they have period $\frac{1}{2}$
and are of class $C^2$. 
Let $I^{\pm}_j$ and $\omega^\pm_j = \partial_{I^{\pm}_j}\mathcal H_2,$  $j \ge 1,$ be the corresponding KdV actions and frequencies and define the following sequence of KdV Hamiltonians $(N \ge 2)$
\[
{\mathcal H}^N_{KdV} := \frac{1}{2N} {\mathcal H}_1 - \frac{1}{24}
\frac{1}{(2N)^3} {\mathcal H}_2
\]
whose frequencies can be computed to be
 \[ \partial _{I^{\pm}_n} {\mathcal H}^N_{KdV} (q_{\pm})= 
\frac{2 \pi n}{N} - \frac{1}{24} \frac{1}{(2N)^3} \omega ^{\pm}_n. 
  \]
  
To describe the asymptotics of the frequencies $\omega^N_n$ for $n \sim 1$ and $n \sim N-1$ consider functions 
$F : \mathbb N \rightarrow {\mathbb R}_{\geq 1}$
satisfying
\[
(F) \qquad \lim _{N \rightarrow \infty} F(N) = \infty;\qquad   F \mbox{ increasing; }  \qquad
 F(N)  \leq  N^\eta  \,\, \mbox{ with } \,\, \eta > 0.
\]
%%%%%%%%%%%%%%%%%%%%%%%%%%%%%%%%%%%%%%%%%%%%%%%%%%%%%%%%%%%%%%%%%%%%%%%%
\begin{guess}
\label{thm1.4} 
%\begin{guess}
\label{thm1.5inSec6} 
 Let $F : {\mathbb N} \rightarrow {\mathbb R}_{\geq 1}$ 
satisfy (F) with $\eta \le 1/3$ and set  $M =~ [ F(N) ]$.
% and $L= [ F(M) ]$ .
Then the asymptotics of the frequencies $(\omega^N_n)_{0 < n < N}$
of $(b^N, a^N)$ defined by \eqref{1Da.300}
are as follows: 

at the left and right edges: for $1 \leq n \leq F(M)$
\begin{equation}
\label{freq1inSec6}
%\label{freq1} 
\omega^N_n= 
\partial _{I^{-}_n}{\mathcal H}^N_{KdV}
%\frac{1}{N^3}
%\frac{2 \pi n}{N} - \frac{1}{24} \frac{1}{(2N)^3} %\omega ^{-}_n
+ 
O\left(\frac{1}{N^3}\big( \frac{n^2F(M)}{M^{1/2}} + \frac{1}{F(M)^{5/2}}\big) \right)
\end{equation}
\begin{equation}
\label{freq2inSec6}
%\label{freq2}
\omega^N_{N-n} = 
\partial _{I^{+}_n}{\mathcal H}^N_{KdV}
% \frac{2 \pi n}{N} - \frac{1}{24} \frac{1}{(2N)^3} %\omega ^{+}_n
 + 
O\left(\frac{1}{N^3}\big( \frac{n^2F(M)}{M^{1/2}} + \frac{1}{F(M)^{5/2}}\big) \right)
\end{equation}
in the bulk: $M < n < N -M $
\begin{equation}
\label{freq3inSec6}
%\label{freq3}
\omega^N_n=2\sin{\frac{\pi n}N} \Big( 1+  O\big( \frac{\log M}{M^2} \big) \Big).
\end{equation}
%in the bulk: $[N^\fffg] < n < N-[N^\fffg]$
%\be\label{freq3}
%\omega^N_n=2\sin{\frac{\pi n}N} \left( 1+ O\left( N^{-\fffg}+N^{2\fff-1} \right) \right).
%\ee
Finally, for $F(M) < n \leq M$
\be
\label{freq4inSec6}
%\label{freq4}
0<\omega^N_n=\frac{2\pi n}N +O \left( \frac{n^3}{N^{3}} \right), \qquad 0 <\omega^N_{N-n}=
\frac{2\pi n}N +O\left(\frac{n^3}{N^{3}} \right).
\ee
These estimates hold uniformly in $0 < n < N$ and uniformly on bounded
subsets of functions $\alpha , \beta $ in $C_0^2(\mathbb T)$.
\end{guess}

%Let $\nu_0 \leq\fff<\tres$ with $\nu_0 = \frac{2}{17} ???$ and $\tres = ???.$
%% and $\mu=\min\{-3\fff^2+\frac\fff 2,-\fff^2-2\fff+1\}$.
%Then the asymptotics of the frequencies of the periodic Toda chains $(b^N, a^N)$ defined by \eqref{1Da.300} are as follows:
%
%\medskip
%
%at the left and right edges: for $ 1\leq n\leq N^\fffg$
%\be\label{freq1}
%\omega^N_n= \partial _{I^{-}_n}{\mathcal H}^N_{KdV}+ O(N^{-3 - \freq})
%% \partial _{I^{-}_n}{\mathcal H}^N_{KdV} = \frac{2 \pi n}{N} - 
%%\frac{1}{24} \frac{1}{(2N)^3} \omega ^{-}_n
%\ee
%\be\label{freq2}
%\omega^N_{N-n} = \partial _{I^{+}_n}{\mathcal H}^N_{KdV}+ O(N^{-3 - \freq})
%%\partial _{I^{+}_n}{\mathcal H}^N_{KdV} = \frac{2 \pi n}{N} - 
%%\frac{1}{24} \frac{1}{(2N)^3} \omega ^{+}_n 
%\ee
%in the bulk: for $N^\fff < n < N-N^\fff$
%\be\label{freq3}
%\omega^N_n=2\sin{\frac{\pi n}N} \left( 1+  O\left( \frac{\log N} {N^{2\fff}} \right) \right).
%\ee
%%in the bulk: $[N^\fffg] < n < N-[N^\fffg]$
%%\be\label{freq3}
%%\omega^N_n=2\sin{\frac{\pi n}N} \left( 1+ O\left( N^{-\fffg}+N^{2\fff-1} \right) \right).
%%\ee
%Finally, for $N^\fffg < n \leq N^\fff$
%\be\label{freq4}
%0<\omega^N_n=\frac{2\pi n}N +O \big( \frac{n^3}{N^{3}} \big), \qquad 0 <\omega^N_{N-n}=
%\frac{2\pi n}N +O\big(\frac{n^3}{N^{3}} \big).
%\ee
%
%These estimates hold uniformly in $1 \leq n \leq N - 1$ and uniformly on bounded
%subsets of functions $\alpha , \beta $ in $C_0^2(\mathbb T)$.
%Note that the best remainder for \eqref{freq1}-\eqref{freq2} is obtained for $\fff=\frac2{17} ???$ and is $O(N^{-3-\frac5{289}}) ???$.
%\end{guess}
We remark that uniformly for $1\leq n\leq F(M)$, the error terms in the asymptotics \eqref{freq1inSec6} and \eqref{freq2inSec6} have stronger decay as $N\to\infty$ than the principal terms 
$\partial _{I^{\pm}_n}{\mathcal H}^N_{KdV} = 
\frac{2 \pi n}{N} - \frac{1}{24} \frac{1}{(2N)^3} \omega ^{\pm}_n$  in view of the asymptotics
$\omega ^{\pm}_n = (4n\pi)^3 + O(1)$ (cf Section \ref{frequencies}).

\medskip

To state the asymptotics of the actions $I^N_n,$ $0 < n < N$, we first need to introduce some additional notation. Recall that the Hill operators $H_\pm := - \partial ^2_x + q _\pm$, associated to the potentials 
$q_\pm,$ come up in the formulation of the KdV equation in terms of a Lax pair. In particular, 
the spectrum of 
$H_\pm,$ when considered with periodic boundary conditions on the interval $[0,1]$,  is preserved 
by the KdV flow. It is pure point and consists of real eigenvalues which when
listed in increasing order and with their multiplicities satisfy
   \[ \lambda ^\pm _0 < \lambda ^\pm _1 \leq \lambda ^\pm _2 < \cdots.
   \]
For any $n\ge 1$, the difference 
$\gamma^{\pm}_n := \lambda^{\pm}_{2n} - \lambda^{\pm}_{2n-1}$ is referred to as n'th gap length. It is well know that the decay properties of $\gamma^\pm_n$
as $n \to \infty$ are related to the smoothness of $q_\pm.$ In particular, as $q_\pm$ are of class $C^2,$ one has $\sum_{n=1}^\infty n^{4} (\gamma_n^\pm)^2 < \infty$ , hence in particular $\gamma_n^\pm = O(n^{-2})$.
%%%%%%%%%%%%%%%%%%%%%%%%%%%%%%%%%%%%%%%%%%%%%%%%%%%%%%%%%%%%%%%%%%%%%%%%%
\begin{guess}
\label{thm1.5} 
 Let $F : {\mathbb N} \rightarrow {\mathbb R}_{\geq 1}$ 
satisfy (F) with $\eta < 1/2$ and set $M=~ [ F(N) ]$.
Then the asymptotics of the actions $I^N = (I^N_n)_{0 < n < N}$
of the states $(b^N, a^N)$ defined by \eqref{1Da.300}
are as follows: 

\medskip

at the left and right edges: for $1 \leq n \leq F(M)$
\begin{equation}
\label{1.6} 8N^2 I^N_n = I^-_n + 
O \Big( \frac{M^2}{N} \frac{F(M)} {M^{1/2}} +\frac{M^3}{N^{3/2}} +
\gamma ^-_n (\frac{F(M)}{M^{1/2}} + \frac{M}{N^{1/2}} ) \Big) 
\end{equation}
\begin{equation}
\label{1.7} 8N^2  I^N_{N-n} = I^+_n +
O \Big( \frac{M^2}{N} \frac{F(M)}{M^{1/2}} + \frac{M^3}{N^{3/2}}  + \gamma ^+_n (\frac{F(M)}{M^{1/2}} +\frac{M}{N^{1/2}} )\Big) 
\end{equation}
in the bulk: $M < n \leq N/2$
\[ I^N_{n}, I^N_{N-n} = O \Big(  \frac{1}{n M^2} \frac{1}{N^{2}}\Big)
\]
whereas for $F(M) < n \le M$,
\[ I^N_n = O \Big( \frac{1}{n} \big( (\gamma ^-_n)^2 + \frac{M^4}{N^2} \big)\frac{1}{N^2} \Big),
\qquad
 I^N_{N-n} = O \Big(\frac{1}{n} \big( (\gamma ^+_n)^2 +\frac{M^4}{N^2}\big) \frac{1}{N^2} \Big).
\]
These estimates hold uniformly in $0 < n < N$ and uniformly on bounded
subsets of functions $\alpha , \beta $ in $C_0^2({\mathbb T})$.
\end{guess}
%\begin{guess}
%\label{thm1.3} For any $0 < \nu < 1/2$, the asymptotics of the actions  of the periodic Toda chains $(b^N, a^N)$  defined by \eqref{1Da.300}
%are as follows: 
%
%\medskip
%
%at the left and right edges: 
%for $1 \leq n \leq N^{\nu^2}$
%\begin{equation}
%\label{1.6} I^N_n = \frac{1}{8N^2} 
%\Big(  I^-_n + 
%O \big( N^{-1+\frac32\nu+\nu^2} + \gamma ^-_n N^{-\frac\nu2+\nu^2} \big)
%\Big)
%\end{equation}
%\begin{equation}
%\label{1.7} I^N_{N-n} = \frac{1}{8N^2} \Big( I^+_n+ 
%O \big(N^{-1+\frac32\nu+\nu^2}  + \gamma ^+_n N^{-\frac\nu2+\nu^2} \big)\Big)
%\end{equation}
%
%in the bulk: $N^\nu < n \leq N/2$
%\[ I^N_{n}, I^N_{N-n} = O \big( \frac{1}{N^2} \frac{1}{n N^{2\nu}} \big).
%\]
%
%Finally, for $N^{\nu^2} < n \leq N^\nu$
%\[ I^N_n = O \Big( \frac{1}{N^2} \frac{1}{n} \big( (\gamma ^-_n)^2 + 
%{N^{-2+4\nu}} \big) \Big),
%\quad 
% I^N_{N-n} = O \Big( \frac{1}{N^2} \frac{1}{n} \big( ( \gamma ^+_n)^2 +
%N^{-2+4\nu}\big) \Big).
%\]
%These estimates hold uniformly in $0 \leq n \leq 2N-1$ and uniformly on bounded
%subsets of functions $\alpha , \beta $ in $C^2_0(\mathbb T).$
%\end{guess}

\medskip  

The following remark comments on the implications of Theorem \ref{thm1.4} and Theorem \ref{thm1.5} on the approximation of the Toda chain with initial data given by $(b^N, a^N)$ by KdV type solutions. For simplicity we assume that $F(N) = N^\eta$ with $0 < \eta \le 1/3$. Then $M = [N^\eta]$ and $F(M) \sim N^{\eta^2}.$

\begin{remark} 
In Birkhoff coordinates, the n'th component $(x^N_n(t) , y^N_n(t))$ of the solution, 
$ 0 < n < N,$ is of the form 
\[
(x^N_n(t) , y^N_n(t) )= \sqrt{2I^N_n} ( \cos(\theta^N_n +t \omega^N_n), \sin(\theta^N_n +t \omega^N_n))
\quad \mbox{if} \quad I^N_n \ne 0
\]
and zero otherwise where $\theta^N_n$ is the n'th angle coordinate, determined by the initial data.
For any given $N$, it evolves on the torus 
\[
\mathcal T^N = \{(x^N_n, y^N_n)_{0 < n < N} | 
\, (x^N_n)^2 + (y^N_n)^2 = 2 I^N_n  \,\,\,\forall 0 < n < N \}.
\]
First note that according to Theorem~\ref{thm1.5} the size of the components 
of the solution in the bulk is small in the sense that
\[
\sum_{ F(M) < n < N -1 -F(M)}  I^N_n = \, 
 O\Big( \frac{1}{ N^2} \big( \frac{1}{N^{5 \eta^{2}}} + \frac{\log N}{N^{2\eta}} \big) \Big)
\]
%whereas
%\[
%\sum_{1 \le n \le F(M)}  I^N_n + I^N_{N-n} 
%= \frac{1}{8N^2}  \sum_{1 \le n \le F(M)} I^-_n  + I^+_n \, + \,
%O \Big( \frac{1}{N^2} \big( N^{-\eta/2 + \eta^2} +  N^{-1 +\frac{3}{2}\eta +2 \eta^2} \big) \Big).
%\]
Here we used the fact mentioned above that $\sum_{n=1}^\infty n^{4} (\gamma_n^\pm)^2 < \infty$.
Hence for $N$ sufficiently large, the solution of the Toda chain with initial data $(b^N, a^N)$ can be viewed as a small perturbation of a long wave obtaind from the solution by setting the components 
$(x^N_n(t) , y^N_n(t) )$ with $F(M) < n < N-1 - F(M)$ to zero. \\
Secondly we note that it follows from Theorem \ref{thm1.4} that on a time interval of size larger than $N^{3}$, these long waves are approximated by two KdV type solutions. More precisely, in the case 
where $0 < \eta < \frac{1}{11},$ one has for any $1 \le n \le F(M)$
by \eqref{freq1inSec6} and \eqref{freq2inSec6}, 
\[
\omega^N_n - \partial _{I^{-}_n}{\mathcal H}^N_{KdV},\,\,\,
\omega^N_{N-n} - \partial _{I^{+}_n}{\mathcal H}^N_{KdV} = O(\frac{1}{N^3 N^{5\eta^2/2}}).
\]
Hence the approximation of solutions of Toda chains considered above
is valid on a longer time interval than the one obtained in \cite{BP} (cf also
\cite{SchW}).\\
The asymptotics of the Toda frequencies up to order 5 -- which in principle is
possible -- will incorporate in \eqref{freq1inSec6}-\eqref{freq2inSec6} further terms,
conjecturely involving the third Hamiltonian of the KdV hierarchy, and thus will provide
an approximation beyond the one by KdV type solutions.
\end{remark}
\medskip
%Typically, asymptotics as the ones stated in Theorem \ref{thm1.4} and Theorem
%\ref{thm1.3} are out of reach by current techniques.
%Two key ingredients allowed us in the case of periodic Toda chains to overcome the
%difficulties typically encountered in deriving such asymptotic expansions. 
{\em Outline of the proofs of Theorem \ref{thm1.4} and Theorem \ref{thm1.5}: } 
A first key ingredient are the asymptotic expansions of the eigenvalues of Jacobi
matrices $L(b^N, a^N)$ obtained in our paper \cite{bkp2}, having the novel feature that
they involve the periodic eigenvalues of two independent Hill operators. Secondly,
we rely on the rich structure of periodic Toda chains and the KdV equation related
to their integrability. In particular, we use that their spectral curves are of a
very similar nature.
The proof of Theorem \ref{thm1.5} involves formulas, representing the Toda and KdV
actions as periods of differentials on the corresponding spectral curves where the
cycles involved are closely related to the spectrum of the underlying Jacobi
matrices respectively Hill operators. The spectral curve as well as the
differentials are expressed in terms of the discriminants, associated to the Toda
chains respectively KdV equation and the asymptotics of the Toda actions 
are obtained from the asymptotic expansion of the
discriminant
$\Delta _N(\mu )$, associated to the sequence of Toda chains \eqref{1Da.300}, 
derived in \cite{bkp2} -- see Section \ref{spec} for a precise statement of this result.
The proof of Theorem \ref{thm1.4} is more involved. Again the Toda and the KdV
frequencies can both be represented as periods of certain differentials on the corresponding
spectral curves. Unlike in the case of the actions there are however significant structural
differences between the differentials involved in the two cases.
To continue we first need to introduce additional notation. 
 It turns out to be more convenient
to double the size of the Jacobi matrix $L(b^N, a^N)$ and to consider
\[
Q^{\alpha , \beta }_N = L \big( (b^N, b^N), (a^N, a^N) \big).
\]
The eigenvalues of the symmetric $2N \times 2N$ matrix $Q^{\alpha , \beta }_N$
when listed in increasing order and with their multiplicities then satisfy
\[
\lambda ^N_0 < \lambda ^N_1 \leq \lambda ^N_2 < \cdots < \lambda ^N_{2N-3} \leq
\lambda ^N_{2N-2} < \lambda ^N_{2N-1} .
\]
Precise asymptotics as $N\to\infty$ of the spectrum of $Q^{\alpha , \beta }_N$ were derived in \cite{bkp2} and are recalled in Section \ref{spec}.
Furthermore, let $\chi_N(\mu)$ be the characteristic polynomial of $Q_N^{\alpha, \beta}$
\[
\chi_N (\mu) = \prod _{0 \le k \le 2N-1} (\mu - \lambda ^{N}_k ).
\]
By Floquet theory,  it  equals $\Delta^2 _N(\mu ) -4$ up to a $\mu-$independent factor
where $\Delta _N(\mu )$ is the discriminant of the Toda chain $(b^N, a^N)$ (see Section \ref{spec}).
This factor has been computed in \cite{bkp2}, 
\begin{equation} \label{frakq}
\Delta ^2_N(\mu ) - 4 = \frak q^{-2}_N \chi_N(\mu ) \quad \mbox{ with } \quad  
\frak q_N = \prod _{n=1}^{N} a_n^N.
\end{equation}
Finally denote by $\sqrt[c]{\Delta^2_N(\mu ) - 4}$ the canonical root of 
$\Delta^2_N(\mu )- 4$, defined for $\mu \in \mathbb C \setminus ((- \infty, \lambda^N_0]\cup [\lambda^N_1, \lambda^N_2]\cup \dots \cup [\lambda^N_{2N-1}, \infty))$ by the sign condition
\begin{equation}
\label{cNroot}
\sqrt[c]{\Delta _N(\mu + i0)^2 - 4} > 0 \qquad \forall\,\, \mu > \lambda ^N_{2N - 1}
\end{equation}
and set $\sqrt[c]{\chi _N(\mu)} = \frak q_N \sqrt[c]{\Delta _N(\mu)^2 - 4} $.
The representation of the Toda frequencies used for deriving their asymptotics is the following one
\begin{equation}
\label{intro.4}
iN\omega ^N_n = \int ^{\lambda ^N_{2n-1}}_{\lambda ^N_0}
\frac{(\mu - \frak p_N/N) \dot \Delta _N(\mu )d \mu}{\sqrt[c]{\Delta _N^2 (\mu - i0) - 4}}
- \sum _{\underset{I_k^N \ne 0}{0 < k < N}} I^N_k  \omega ^N_k
\int ^{\lambda ^N_{2n-1}}_{\lambda ^N_0} \frac{\varphi ^N_k(\mu )d \mu}
{\sqrt[c]{\chi_N (\mu - i0)}}.
\end{equation}
Here $\frak p_N$ is the trace of $L(b^N , a^N)$, $\frak p_N=tr L(b^N , a^N)$,  
$\dot \Delta_N(\mu )=\frac d{d\mu}\Delta _N(\mu )$ and for any $ 0 < n < N$,
$\varphi ^N_n(\mu )$ is a polynomial of degree $N-2,$
\[ \varphi ^N_n(\mu ) =
\prod _{\underset{k\not=n}{0 < k < N}} (\mu - \sigma ^{N,n}_k ),
\]
whose $N-2$ zeroes  $(\sigma ^{N,n}_k)_{0 < k \ne n <N}$ are uniquely determined
determined by the $N-2$ normalization conditions
\be\label{boot} \frac{1}{2\pi } \int _{\Gamma ^N_k} 
\frac{\varphi ^N_n(\mu)}{\sqrt[c]{\chi_N (\mu)}} d\mu = 0 \qquad \forall \,\, 0 < k\neq n < N.
\ee
Here $(\Gamma ^N_k)_{0 < k < N}$, denote pairwise disjoint counterclockwise
oriented contours in $\mathbb C,$ chosen in such a way that $\lambda ^N_{2k-1},
\lambda ^N_{2k}$
are inside of $\Gamma ^N_k$ whereas all other eigenvalues of $Q^{\alpha , \beta }_N$
are outside. (Actually, $(\frac{\varphi ^N_k(\mu )} {\sqrt{\chi_N (\mu)}}d \mu)_{0 < k < N}$ 
are differential forms on the spectral curve defined by $\chi_N(\mu)$ -- 
see Section \ref{frequencies} for more details.)

Unlike for the corresponding formula for the KdV frequencies, \eqref{intro.4} is a
{\em system of equations} for the frequencies $(\omega ^N_n)_{0 < n < N}$ rather
than a formula
as the Toda frequencies also appear on the right hand side of the identity
\eqref{intro.4}.
Fortunately, the two terms on the right hand side of \eqref{intro.4} are not of the
same order as \Ninf. It turns out that the first one dominates the second one,
suggesting to use a two step approach to prove the claimed asymptotics. In a first
step (Section \ref{priori}), we compute the leading order of the asymptotics of
$\omega ^N_n$, using the following alternative
representation of the Toda frequencies.
%%%%%%%%%%%%%%%%%%%%%%%%%%%%%%%%%%%%%%%%%%%%%%%%%%%%%%%%%%%%%%%%%%%%%%%%
\begin{guess}\label{lemnouv}
The frequencies $\omega ^N_n,\
0 < n < N,$ of the states $(b^N, a^N)$ defined by
\eqref{1Da.300} are strictly positive, $\omega ^N_n> 0,$ and given by the expression
\begin{equation}
\label{B.1} \omega ^N_n = \Big( \frac{1}{2\pi } \int _{\Gamma ^N_n}
\frac{\varphi ^N_n(\mu )}{\sqrt[c]{\chi_N(\mu )}} d\mu \Big) ^{-1}.
\end{equation}
\end{guess}
We use Theorem \ref{lemnouv} to prove the asymptotics \eqref{freq3inSec6} of the frequencies in the bulk in Section \ref{priori}, 
%, The leading order asymptotics of the frequencies
%are etablished in Section \ref{priori}. In particular, \eqref{freq3inSec6} is proved in 
Proposition \ref{Proposition B.1}. 
In a second step we establish the asymptotics 
\eqref{freq1inSec6}, \eqref{freq2inSec6}, and \eqref{freq4inSec6} of Theorem \ref{thm1.4}.
The proof of these asymptotics is based on formula \eqref{intro.4}, uses the first order asymptotics mentioned above and relies on a comparison
with corresponding formulas for the KdV frequencies (Section \ref{frequencies}) as well as the asymptotics of the zeroes of the functions $\varphi_n^N$ (Section \ref{zeroes})  and  the Toda actions $I^N_n$ (Section \ref{actions}). The formulas \eqref{intro.4} of the Toda frequencies are reviewed in Section \ref{frequencies}. There we also prove Theorem \ref{lemnouv}.
%The proof of Theorem \ref{thm1.5inSec6} is quite involved and hence we split it up 
%in three subsequent sections. The leading order asymptotics of the frequencies
%are etablished in Section \ref{priori}. In particular, \eqref{freq3inSec6} is proved in 
%Proposition \ref{Proposition B.1}.
 The higher order asymptotics \eqref{freq4inSec6} 
are shown in Proposition \ref{Proposition B.3A}.
Using a symmetry of the Toda chain, discussed in Appendix \ref{TodaSymmetry}, the proof of 
\eqref{freq2inSec6} can be reduced to the one of \eqref{freq1inSec6}. The latter asymptotics are proved in Section \ref{ProofFrequ} using auxilary estimates derived in Appendix \ref{tkdv}.

{\em Related work:}
In the fifties, Fermi, Pasta, and Ulam introduced
and studied a model consisting of a chain of particles interacting with their nearest
neighbors through nonlinear strings, referred to as FPU chain.
Expecting that the thermalization of small energy solution of such chains is valid
they wanted to compute in numerical experiments the partition of energy
among the normal modes of the linearized chain to see
at what rate equipartition of energy is achieved. Much to their surprise,
instead of equipartition of energy, they observed recurrent features of these solutions.
Their report \cite{FPU} had a far reaching impact.
To our days, the questions raised with their experiments have been intensely studied.
However, the question if thermalization occurs has not been settled so far. It requires
to compute delicate asymptotics for large number of particles. For status reports
on the investigations of the FPU chains at the occasion of the 50th anniversary
of the FPU experiments see the articles/books \cite{BI}, \cite{CGG}, \cite{G}. 

Formally, the dynamics of FPU chains of particles corresponding to long wave initial
data can be approximated by solutions of the KdV equation at least for some interval of time --
see e.g. \cite{ZK}, \cite{T}, \cite{C}.
But only recently this has been rigorously proved for certain classes of initial data
(\cite{SchW}, \cite{BP}). More precisely, in \cite{SchW} it was proved that for long
waves on the real line with initial date of size $\epsilon^{2}$ and sufficiently
strong decay,
the corresponding solutions of the FPU chain
can be approximated up to a translation by the superposition of a right--going and a
left--going wave, both evolving according to KdV, over a time interval of order
$\epsilon^{-3}$.
Despite the fact that in the periodic set-up, waves interact much stronger, an
analogous result for special type of analytic initial data could be proved in this
set-up as well by \cite{BP}.

In the case of the Toda chain, Toda himself formally addressed the problem of
approximating solutions of the
Toda chain by solutions of KdV with the aim of explaining the recurrent features
of FPU chains at least in this integrable case -- see \cite{T}. Typically, in
investigations of this kind,
a limit of Toda chains of the type describe above is obtained by computing the
formal limit of the Jacobi matrices $L$ for a sequence of initial data, obtained by
discretizing two given smooth functions, one for the diagonal and the other  for the
upper (and hence also lower) diagonal of $L$. In such a way one gets formally that
$L$ approaches a Hill operator which via the Lax pair formalism then leads to a
solution of the KdV equation. The fundamental question then arose how a second Hill
operator can be found
providing the second solution of the KdV equation needed to describe the above
mentioned approximation of long waves of periodic Toda chains.

To answer this question, we used as a key ingredient the discovery made in
\cite{BGPU} that the Jacobi matrix describing the Toda chain with $N$ particles can
be viewed as the (geometric) quantization with Planck constant $\hbar$ of a
Hamiltonian system on the torus $\mathbb T^2$, and the limit of such periodic Toda
chains as the semiclassical limit of this quantization where the number $N$ of
particles of a chain coincides with the inverse of the mesh size of the
discretization as well as with $\frac {1}\hbar$. This observation was used in a
crucial way in paper \cite{bkp2} (cf also \cite{BTP}) where among other results we proved an
asymptotic expansion in $N^{-1}$ of the eigenvalues of the Jacobi matrices
in terms of the eigenvalues of two Hill operators with potentials $q_+$ respectively $q_-$, 
introduced above.\\

{\em Acknowledgments:} The authors would like to thank the CNRS, the \'Ecole
polytechnique,
the Swiss National Science Foundation, the University of Milan and the University of
Z\"urich for financial support during the elaboration of this work.

\medskip

%%%%%%%%%%%%%%%%%%%%%%%%%%%%%%%%%%%%%%%%%%%%%%%%%%%%%%%%%%%%%%%%%%%%%%%%%%%%%%%%%%%%%%%
%%%%%%%%%%%%%%%%%%%%%%%%%%%%%%%%%%%%%%%%%%%%%%%%%%%%%%%%%%%%%%%%%%%%%%%%%%%%%%%%%%%%%%%%%
%%%%%%%%%%%%%%%%%%%%%%%%%%%%%%%%%%%%%%%%%%%%%%%%%%%%%%%%%%%%%%%%%%%%%%%%%%%%%%%%%%%%%%%%%%

\section{Spectral asymptotics}\label{spec}

For the convenience of the reader, we recall in this section results established in \cite{bkp2} which will be used throughout in the sequel.
 The first result concerns the asymptotics of the eigenvalues 
\[ 
\lambda^N_0 < \lambda^N_1 \le \lambda^N_2 < \cdots < \lambda^N_{2N-3} \le 
\lambda^N_{2N-2} < \lambda^N_{2N-1}
\] 
 of the Jacobi matrices 
\be\label{jmq}
Q_N^{\alpha, \beta} \equiv Q(b^N, a^N) = L \big( (b^N, b^N), (a^N, a^N) \big)
\ee
 where $(b^N, a^N)$ are the states introduced in \eqref{1Da.300},
\[
b^N_n = \frac1{4N^2} \beta \big( \frac{n}{N} \big), \quad a^N_n =
      1 + \frac1{4N^2} \alpha \big( \frac{n}{N} \big) \quad \mbox{with}\quad  
\beta, \alpha \in C^2_0(\mathbb T).
\]
 Explicitly, 
the symmetric $2N \times 2N$ matrix $Q_N^{\alpha, \beta} $ is given by     
\[
Q_N^{\alpha, \beta}  =
 \begin{pmatrix} 
   b^N_1 &a^N_1 &0&\ldots&\ldots&\ldots&\ldots&\ldots&0&a^N_N \\ 
   a^N_1 &b^N_2 &a^N_2&0 &\ldots&\ldots&\ldots&\ldots&\ldots&0 \\ 
   0 &a^N_2 &b^N_3&a^N_3&0&\ldots&\ldots&\ldots&\ldots&0\\
   \vdots&&&&&&&&&\vdots\\
   0&\ldots&0&a^N_{N-1}&b^N_N&a^N_N&0&\ldots&\ldots&0\\
   0&\ldots&\ldots&0&a^N_N&b^N_1&a^N_1&0&\ldots&0\\
   \vdots&&&&&&&&&\vdots\\
   0&\ldots &\ldots&\ldots &\ldots&\ldots&0&a^N_{N-2}&b^N_{N-1}&a^N_{N-1}\\
      a^N_N &0&\ldots &\ldots&\ldots&\ldots&\ldots&0&a^N_{N-1}&b^N_N 
      \end{pmatrix}
\]

Using Floquet theory (cf. e.g. \cite{HK1}) one sees that 
the eigenvalues of $ L \big( b^N, a^N)$ can be identified among 
the ones of $Q(b^N, a^N)$ as follows: 
if $N$ is even they are 
\begin{equation}\label{rootsA}
\lambda_0^N < \lambda_3^N \le \lambda_4^N < \cdots < \lambda_{2N-5}^N \le 
\lambda_{2N-4}^N < \lambda_{2N-1}^N
\end{equation}
and hence 
\[
\Delta_N(\mu) -2 = \frak q_N^{-1} (\mu - \lambda ^{N}_0)(\mu - \lambda ^{N}_{2N-1} )
 \prod _{0 < k < N/2} (\mu - \lambda ^{N}_{4k} )(\mu - \lambda ^{N}_{4k-1} )
\] 
whereas if $N$ is odd
\begin{equation}\label{rootsB}
\lambda_1^N \le \lambda_2^N < \lambda_5^N \le \lambda_6^N < \cdots < 
\lambda_{2N-5}^N \le \lambda_{2N-4}^N  < \lambda_{2N-1}^N
\end{equation}
 leading to
\[
\Delta_N(\mu) -2 = \frak q_N^{-1} (\mu - \lambda ^{N}_{2N-1} )
 \prod _{0 < k < N/2} (\mu - \lambda ^{N}_{4k-2} )(\mu - \lambda ^{N}_{4k-3} ).
\] 
To describe the asymptotics of $\lambda ^N_n$ at the edges, $n \sim 1$ or $n \sim 2N-1$,
we introduced in \cite{bkp2} two Hill operators $H_\pm := - \partial ^2_x + q _\pm$ with potentials
   \begin{equation}\label{pot} q_\pm (x) = - 2\alpha (2x) \mp \beta (2x).
   \end{equation}
   
The periodic  eigenvalues $(\lambda ^\pm_n)_{n \geq 0}$ of $H_\pm $ on $[0,1]$, 
when listed in increasing order and with multiplicities are  known to satisfy
\[
\lambda ^\pm _0 < \lambda ^\pm _1 \leq \lambda ^\pm _2 < 
\lambda ^\pm _3 \leq \lambda ^\pm _4 < \cdots.
\]
As $q_{\pm}$ have average zero, one has in addition that $\lambda ^\pm _0 < 0.$
%To describe the edges of the spectrum of $Q(b^N, a^N)$ we consider functions 
%$F : \mathbb N \rightarrow {\mathbb R}_{\geq 1}$
%satisfying
%\[
%(F1) \,\,\, \lim _{N \rightarrow \infty} F(N) = \infty \quad \mbox{ and } \quad 
%(F2) \,\,\, F(N) \leq N^\eta \,\, \mbox{ with } \,\, \eta > 0.
%\]

\begin{guess}
\label{Theorem 1.1} Let $F : {\mathbb N} \rightarrow {\mathbb R}_{\geq 1}$  
satisfy (F) with $\eta < 1/2$ and set $M =~ [ F(N) ]$. Then for any $\alpha , \beta$  in 
$C_0^2(\mathbb T),$ the asymptotics of the eigenvalues 
$(\lambda ^N_n)_{0 \le n \le 2N-1}$ of $Q_N^{\alpha, \beta}$ are as follows:

\medskip

at the left and right edges: for $0 \leq n \leq 2M$
\[ \lambda ^N_n = - 2 + \frac{1}{4N^2} \lambda ^-_n + O(\frac{M^2}{ N^{3}}),
   \qquad \lambda ^N_{2N-1-n} = 2 - \frac{1}{4N^2} \lambda ^+_n + O(\frac{M^2}{N^{3}})
\]

\medskip

in the bulk: for $ n = 2\ell , 2 \ell - 1$ with $M < \ell < N - M$,
    \[ \lambda ^N_n = - 2 \cos \frac{\ell \pi }{N}
       + O(\frac{1}{MN^{2}}) .
   \]
These estimates hold uniformly in $0 < n < N$ and uniformly on bounded subsets of
functions $\alpha , \beta $ in $C_0^2({\mathbb T})$.
\end{guess}

\medskip

The second result which we recall from \cite{bkp2} concerns the asymptotics of
the discriminant associated to the difference equation 
\be\label{14bis}
a^N_{k-1} y(k-1) + b^N_k y(k) + a^N_ky(k+1) = \mu y(k) \quad \forall k \in \mathbb Z.
\ee
The discriminant is the trace of the Floquet operator of \eqref{14bis} and given by
\[ 
\Delta _N(\mu ) = y_1(N,\mu ) + y_2(N + 1, \mu )
\]
 where $y_1^N$ and $y_2^N$ are the fundamental solutions of \eqref{14bis} determined by
\[ 
y_1(0,\mu ) = 1,\,\, y_1(1,\mu ) = 0 \quad \mbox { and } \quad y_2(0,\mu ) = 0, \,\, y_2(1,\mu )= 1.
\]
Analogously,   the discriminant $\Delta_\pm(\lambda)\equiv\Delta(\lambda,q_\pm)$ of 
\be\label{14ter}
-y''(x,\lambda)+q_\pm(x)y(x,\lambda)=\lambda y(x,\lambda),
\ee
is defined as the trace of the Floquet operator associated to \eqref{14ter},
\[
\Delta_\pm(\lambda)=y_1^\pm(1/2,\lambda)+(y_2^\pm)'(1/2,\lambda)
\]
where $y_1^\pm(x,\lambda)$ and $(y_2^\pm)'(x,\lambda)$ are the fundamental solutions of \eqref{14ter} determined by
\[
y_1^\pm(0,\lambda)=1,\,\, (y_1^\pm(0,\lambda))'=0 \quad \mbox{ and } \quad
y_2^\pm(0,\lambda)=0,\,\, (y_2^\pm(0,\lambda))'=1.
\]

For $F: \mathbb N \to \mathbb  R_{\ge 1}$ as in Theorem \ref{Theorem 1.1} and  $M = [F(N)]$ let
   \[  \Lambda^{\pm , M}_2:= [\lambda ^\pm _0 - 2, \lambda
      ^\pm_{2[F(M)]} + 2] + i[-2,2]
   \]
and choose $N_0 \in {\mathbb Z}_{\geq 1}$ so that
   \be\label{81bis} \lambda ^\pm_{2k+1} - \lambda ^\pm _{2k} \geq 6 \quad \forall \
   k \geq F(F(N_0))
   \ee
By the Counting Lemma for eigenvalues of Hill operators, $N_0$ can be chosen uniformly for
subsets of function $\alpha , \beta $ in $C_0^2({\mathbb T})$.

\begin{guess}
\label{Theorem9.1} 
 Let $F : {\mathbb N} \rightarrow {\mathbb R}_{\geq 1}$ 
satisfy (F) with $\eta < 1/2$ and set $M =~ [ F(N) ]$.
Then for any $\alpha, \beta \in C^2_0(\mathbb T)$
and $N \ge N_0$, one has uniformly for $\lambda$ in $\Lambda^{-,M}_2$,
   \begin{equation}
    \label{9.1} \Delta_N \big(-2+\frac{1}{4N^2}\lambda\big) = (-1)^N\Delta_-(\lambda) + O
                \big(\frac{F(M)^2}{M} \big).
    \end{equation}
 Similarly, uniformly for $\lambda$ in $\Lambda^{+,M}_2$
  \begin{equation}
  \label{9.2} \Delta_N (2-\frac{1}{4N^2}\lambda) = \Delta_+(\lambda) + O\big(\frac{F(M)^2}
              {M}\big).
  \end{equation}
 These estimates hold uniformly on bounded subsets of functions $\alpha, \beta $ in 
$C_0^2({\mathbb T})$.
\end{guess}

\medskip

The latter theorem has been applied in \cite{bkp2} to obtain asymptotics for 
the derivatives of the descriminant. In this paper we will need such asymptotic estimates for
the first derivative,
$\dot \Delta_N(\mu):=\frac d {d\mu}\Delta_N(\mu)$. Let
\[
\Lambda^{\pm , M}_1:= [\lambda ^\pm _0 - 1, \lambda^\pm_{2[F(M)]} + 1] + i[-1,1].
\]
Then, under the same assumptions as in Theorem \ref{Theorem9.1} the following holds:
\begin{corollary}
\label{Corollary 8.7} Uniformly for $\lambda $ in
$\Lambda ^{-,M}_1$
   \[ \frac{1}{4N^2} \dot \Delta _N\big( -2 + \frac{1}{4N^2} \lambda \big) 
       = (-1)^N \dot \Delta _-(\lambda ) + O \Big( \frac{F(M)^2}{M} \Big)
   \]
and similarly, for $\lambda $ in $\Lambda ^{+ , M}_1$
   \[ -\frac{1}{4N^2} \dot \Delta _N\big( 2 - \frac{1}{4N^2} \lambda
      \big) = \dot \Delta _+(\lambda ) + O \Big( \frac
      {F(M)^2}{M} \Big) .
   \]
These estimates hold uniformly on bounded subsets of functions $\alpha , \beta $ in
$C_0^2({\mathbb T})$.
\end{corollary}

For later reference we record that in view of (\ref{frakq}) and (\ref{rootsA}) - (\ref{rootsB}), 
$ \dot \Delta _N(\mu )$ admits the product
representation
\begin{equation}\label{DeltaDot}
\dot \Delta _N(\mu ) = N \frak q_N^{-1} \cdot\prod _{0 < k < N} (\mu - \dot \lambda ^{N}_k)
\end{equation}
where the zeros $(\dot \lambda ^N_n)_{0 < n <N}$ of $ \dot \Delta _N(\mu )$
 are listed in increasing order and satisfy
\[
 \lambda ^N_{2n-1} \le \dot \lambda ^N_n \le \lambda ^N_{2n} \qquad \forall \,\, 0 < n < N.
\]

\begin{remark}
\label{Remark0}
As mentioned in the introduction, instead of considering the states $(b^N,a^N)$, alternatively, one could consider $(d^N,c^N)$, obtained by expressing $(p^N, q^N)$ in Flaschka variables. 
Recall that the components of $(p^N, q^N)$ were defined in Section 1 by
\[
p^N_n = - \frac1{4N^2} \beta \big( \frac{n}{N} \big) \quad \mbox{and} \quad 
q^N_n = - \frac{2}{4N} \xi \big( \frac{n}{N} \big)
\]
where $\xi$ is the element in $C_0^3(\mathbb T)$, satisfying $\xi ' = \alpha$.
Note that $d^N = b^N$ whereas $c^N_n = a^N_n + O(N^{-3})$ for any $1 \le n \le N $. Indeed, as
\[q^N_{n} - q^N_{n+1} =  \frac{1}{4N} \int_0^{N^{-1}} \alpha (\frac{n}{N} + t) dt \,
\mbox{ and } \, \alpha (\frac{n}{N} + t) = \alpha(\frac{n}{N}) + \int_0^t \alpha ' (\frac{n}{N} +s)ds
\]
one concludes that
\[
c^N_n = \exp( \frac{q^N_{n} - q^N_{n+1}}{2} ) = 1 + \frac{q^N_{n} - q^N_{n+1}}{2} +O(N^{-4}) 
= a_n^N + O(N^{-3})
\]
where the error term is uniform for $\alpha$ in a bounded subset of the Banach space 
$C^2_0(\mathbb T)$. 
Let $Q(b^N, a^N) = L((b^N, b^N), (a^N,a^N))$ and similarly, define $Q(d^N, c^N)$.
The operator norm of the difference $Q(b^N, a^N) - Q(d^N, c^N) $ can thus be estimated by
\[
\| Q(b^N, a^N) - Q(d^N, c^N) \| = \| Q(0_N, a^N - c^N) \| = O(N^{-3})
\]
Therefore the spectrum of $Q(c^N, d^N)$ is $N^{-3}$-close to the one of $Q(b^N, a^N)$.
Theorem \ref{Theorem 1.1} thus remains true if $(b^N, a^N)$ is replaced by $(d^N, c^N).$
In view of Proposition 8.1 in \cite{bkp2}, it also follows that Theorem \ref{Theorem9.1}
and Corollary \ref{Corollary 8.7} remain true for $(d^N, c^N).$
We leave it to the reader to go through the arguments used in the proofs of the results of the paper to verify that they remain valid if $(b^N, a^N)$ is replaced by $(d^N, c^N)$.
\end{remark}

%We finish this section by recalling well known bounds of the fundamental solutions 
%$y_i(x, \lambda)\equiv y_i(x,\lambda, q)$, $i=1,2,$ of the equation $-y'' + qy = \lambda q$
%with $q \in L^2[0,1]$(cf e.g. \cite{PT}) and apply them to get bounds for the discriminant 
%$\Delta(\lambda) \equiv \Delta(\lambda, q)$, defined in our set-up by
%$ \Delta(\lambda) = y_1(\frac{1}{2}, \lambda) + y_2'(\frac{1}{2}, \lambda)$.
%For any $q \in L^2[0,1]$, let $\|q \|$ denote the 
%$L^2-$norm of $q$, $\|q\| = (\int_0^1 q(x)^2 dx)^{1/2}.$

%\begin{lemma}\label{boundDelta}
%For any $q \in L^2[0, 1]$ and any $\lambda \in \mathbb R,\,\,  0 \le x \le 1$
%\[
%|y_1(x,\lambda)| \le \exp(| \Im \sqrt{\lambda} | + \|q\|) 		
%\]
%\[
%|y'_2(x,\lambda)| \le \exp(| \Im \sqrt{\lambda} |) +\|q\|  \exp(|2 \Im \sqrt{\lambda} | + \|q\|). 
%\]
%As a consequence, for any $\lambda \in \mathbb R,$ 
%\[
%|\Delta(\lambda, q)| \le  (2+ \|q\|) \exp(|2 \Im \sqrt{\lambda} | + \|q\|). 
%\]
%\end{lemma}

%{\em Proof:} See \cite{PT}, Theorem 1 in Chapter 1, 

%%%%%%%%%%%%%%%%%%%%%%%%%%%%%%%%%%%%%%%%%%%%%%%%%%%%%%%%%%%%%%%%%%%%%%%%
%%%%%%%%%%%%%%%%%%%%%%%%%%%%%%%%%%%%%%%%%%%%%%%%%%%%%%%%%%%%%%%%%%%%%%%%
%%%%%%%%%%%%%%%%%%%%%%%%%%%%%%%%%%%%%%%%%%%%%%%%%%%%%%%%%%%%%%%%%%%%%%%%%

\section{Asymptotics of the actions}\label{actions}
The aim of this section is to prove the asymptotics of the action variables for the states $(b^N, a^N)$ stated in Theorem~\ref{thm1.5}.
%They follow from the following slightly more general result by choosing 
%$F(N) = N^\eta$ for the function $F$ defining the edges of the spectrum of $Q^N_{\alpha, \beta}$.
Without further reference, we will use the notation introduced in the previous sections and
%\begin{guess}
%\label{thm1.5} 
% Let $F : {\mathbb N} \rightarrow {\mathbb R}_{\geq 1}$ be an increasing function
%satisfying (F) with $\eta < 1/2$ and set $M = [ F(N) ]$.
%Then the asymptotics of the actions $I^N = (I^N_n)_{0 < n < N}$
%of the states $(b^N, a^N)$ defined by \eqref{1Da.300}
%are as follows: 
%
%\medskip
%
%at the left and right edges: for $1 \leq n \leq F(M)$
%\begin{equation}
%\label{1.6} 8N^2 I^N_n = I^-_n + 
%O \Big( \frac{M^2}{N} \frac{F(M)} {M^{1/2}} +\frac{M^3}{N^{3/2}} +
%\gamma ^-_n (\frac{F(M)}{M^{1/2}} + \frac{M}{N^{1/2}} ) \Big) 
%\end{equation}
%\begin{equation}
%\label{1.7} 8N^2  I^N_{N-n} = I^+_n +
%O \Big( \frac{M^2}{N} \frac{F(M)}{M^{1/2}} + \frac{M^3}{N^{3/2}}  + \gamma ^+_n (\frac{F(M)}{M^{1/2}} +\frac{M}{N^{1/2}} )\Big) 
%\end{equation}
%in the bulk: $M < n \leq N/2$
%\[ I^N_{n}, I^N_{N-n} = O \Big(  \frac{1}{n M^2} \frac{1}{N^{2}}\Big)
%\]
%whereas for $F(M) < n \le M$,
%\[ I^N_n = O \Big( \frac{1}{n} \big( (\gamma ^-_n)^2 + \frac{M^4}{N^2} \big)\frac{1}{N^2} \Big),
%\qquad
% I^N_{N-n} = O \Big(\frac{1}{n} \big( (\gamma ^+_n)^2 +\frac{M^4}{N^2}\big) \frac{1}{N^2} \Big).
%\]
%These estimates hold uniformly in $0 < n < N$ and uniformly on bounded
%subsets of functions $\alpha , \beta $ in $C_0^2({\mathbb T})$.
%\end{guess}
assume that the assumptions of Theorem~\ref{thm1.5} hold.
To prove it we first derive asymptotic estimates for $I^N_n$ not involving KdV actions and valid for any $0 < n < N.$
In the case where $\lambda^N_{2n-1} < \lambda^N_{2n},$ the action variable $I^N_n$ is given by 
(cf \cite{FM}, \cite[Section 3]{HK1})
\begin{equation}
\label{Arnold}
  I^N_n = \frac{1}{\pi } \int ^{\lambda ^N_{2n}} _{\lambda ^N_{2n-1}} \frac{(\mu
      - \dot \lambda ^N_n) \dot \Delta _N(\mu )}{\sqrt[c]{\Delta ^2_N(\mu - i0 ) - 4}} d\mu
\end{equation}
and is zero otherwise. Recall that 
$\dot \lambda ^N_n$ denotes the n'th zero of 
$ \dot \Delta _N(\mu )$ (Section~\ref{spec})
and that the canoncial root ${\sqrt[c]{\Delta ^2_N(\mu ) - 4}}$ was introduced in Section~\ref{intro}.
The above formula for the actions leads to the following estimates.

\begin{proposition}
\label{Proposition C.1}
{\rm (i)} For any $1 \leq n \leq M$
\[
I^N_n =  O \Big( \frac{1}{n} \big( (\gamma_n^-)^2
+ \frac{M^4}{N^2} \big)  \frac{1}{N^2} \Big), \quad 
I^N_{N-n} = O \Big( \frac{1}{n}  \big( (\gamma_n^+)^2 + \frac{M^4}{N^2} \big) \frac{1}{N^2} \Big).
\]
{\rm (ii)} For any  $ M < n \leq \frac{N}{2}$
\[
I^N_{n}, \,  I^N_{N-n}  = O \Big( \frac{1}{nM^2} \frac{1}{N^2} \Big).
\]
The estimates hold uniformly in $n$ and uniformly on bounded
subsets of functions $\alpha , \beta $ in $C_0^2({\mathbb T})$.
\end{proposition}

{\it Proof of Proposition \ref{Proposition C.1}:}  Clearly, we only need to consider
$0 < n < N$ with $\lambda^N_{2n-1} < \lambda^N_{2n},$ as otherwise, $I^N_n = 0$.
As by (\ref{DeltaDot}) and (\ref{frakq}), 
\[
\dot \Delta _N(\mu ) = N \frak q_N^{-1} \cdot\prod _{0 < k < N} (\mu - \dot \lambda ^{N}_k), 
\qquad
\Delta ^2_N(\mu ) - 4 =  \frak q_N^{-2} \prod _{0 \le k \le 2N-1} (\mu - \lambda ^{N}_k)
\]
one then gets by the above formula for the actions
\[
 I^N_n = \frac{N}{\pi } \int ^{\lambda ^N_{2n}}_{\lambda ^N_{2n-1}} 
   \frac{(\mu - \dot \lambda ^N_n)^2}{\sqrt[+]{(\lambda ^N_{2N-1}-\mu)(\mu - \lambda ^N_0)}} 
     \prod _{k\not= n} \frac{\mu - \dot \lambda ^N_k}{w_k^N(\mu)} 
   \cdot \frac{1}{\sqrt[+]{(\lambda ^N _{2n}-\mu)(\mu - \lambda ^N_{2n-1})}}  d\mu
\]
where $w_k^N(\mu)$ denotes the standard root, 
$w_k^N(\mu) = \sqrt[s]{(\mu - \lambda ^N_{2k})(\mu - \lambda ^N_{2k-1})}$,
defined on $\mathbb C \setminus [\lambda^N_{2k-1}, \lambda^N_{2k}]$ -- see Section \ref{zeroes}
for the determination of the sign. Using the identity
$\int ^{\lambda ^N_{2n}}_{\lambda ^N_{2n-1}} 
\frac{1}{\sqrt[+]{(\lambda ^N _{2n}-\mu)(\mu - \lambda ^N_{2n-1})}} d\mu = \pi$ 
and the mean value theorem (cf  Lemma \ref{atif}) one sees that there exists $\lambda ^N_{2n-1}
\leq \mu _\ast \leq \lambda ^N_{2n}$ so that
\[
    I^N_n  \leq \frac{N(\gamma ^N_n)^2}{(\lambda ^N _{2N-1}-\lambda ^N_{2n})^{1/2}
      (\lambda ^N_{2n-1}- \lambda ^N  _0)^{1/2}} 
       \prod _{k\not= n} \frac{\mu _\ast - \dot \lambda ^N_k}{w^N_k(\mu _\ast)}.
  \]
 Using that $1 + x \leq e^x$ for $x \geq 0$ one has, with $\mu \equiv \mu_\ast$,
\[
\prod_{k < n} \frac{(\mu - \dot \lambda ^N_k)^2}{w^N_k(\mu)^2} 
      \leq \prod_{k < n} \frac{\mu -\lambda ^N_{2k-1}}{\mu - \lambda ^N_{2k}} = 
      \prod_{k < n} \left( 1 + \frac{\gamma ^N_k}{\mu  - \lambda ^N_{2k}} \right) \leq 
      \exp \left( \sum _{k < n} \frac{\gamma ^N_k}{\lambda ^N _{2n-1} - \lambda ^N_{2k}} \right),
\]
and similarly
$ \prod_{k > n} \frac{(\mu - \dot \lambda ^N_k)^2}{w^N_k(\mu)^2} 
  \leq \exp \left( \sum _{k > n} \frac{\gamma ^N_k}{\lambda ^N _{2k-1} - \lambda ^N_{2n}} \right).$
By Proposition \ref{basicEstimate} it then follows that 
 $ \prod_{k \not= n} \frac{\mu _\ast - \dot \lambda ^N_k}{w^N_k(\mu _\ast)} =O(1),$
yielding
   \begin{equation}
   \label{C.1}   I^N_n = O \Big( \frac{N(\gamma ^N_n)^2 }{(\lambda ^N_{2N-1} - \lambda ^N
                 _{2n})^{1/2} (\lambda ^N_{2n-1} - \lambda ^N_0)^{1/2} } \Big) .
   \end{equation}
We now discuss the asymptotics of $I^N_n, N \rightarrow \infty $, in the four cases
listed in the statement. To get the claimed estimate for $I^N_n$ with $1 \le n \le M$ note that 
by Theorem~\ref{Theorem 1.1},
   $(\gamma ^N_n)^2 = O
      \big(\frac{1}{N^4}  ((\gamma_n^-)^2+ \frac{M^4}{N^2} ) \big)
   $
and
   \begin{equation}
   \label{C.2} (\lambda ^N_{2N-1}-\lambda ^N_{2n})^{-1/2} = O(1), \qquad  
                    (\lambda ^N_{2n-1}-\lambda ^N_0)^{-1/2} = O\left( \frac{N}{n} \right)
   \end{equation}
yielding the claimed estimate
$    I^N_n = O \big( \frac{1}{N^4}  \frac{1}{n} ( (\gamma_n^-)^2 +
               \frac{M^4}{N^2} ) \big).$
%   \end{equation}
The one for  $I^N_{N-n}$ is shown similarly.
%   = O \left( \frac{1}{N^2} \frac{1}{n} \left( (\gamma_n^+)^2 +
%      \frac{M^4}{N^2} \right) \right)  .
%   \]
In the case $M < n \leq \frac{N}{2},$ Theorem~\ref{Theorem 1.1} implies
 $( \gamma ^N_n )^2 = O  \left( \frac{1}{N^4 M^2} \right) $ as well as
$ (\lambda ^N_{2n-1} - \lambda ^N_0)^{-1/2} = O(\frac{N}{n})$.
Indeed, as by \eqref{B.9uno}, 
$2 (1 - \cos(\frac{n\pi}{N}) ) \ge \frac{\pi n^2}{N^2}$ one gets
\[
\lambda ^N_{2n-1} - \lambda ^N_0 =  \frac{\pi n^2}{N^2} -
\frac{\lambda^-_0}{8N^2} + O\Big( (\frac{M^2}{N} +\frac{1}{M}) \frac{1}{N^2} \Big)
\] 
yielding in view of the inequality $\lambda^-_0 \le 0$ (cf Section \ref{spec}) the stated estimate.
Combining the above estimates leads to 
$I^N_n = O \left( \frac{1}{nM^2} \frac{1}{N^2} \right) $. 
In the case $\frac{N}{2} \leq n < N - M$ one argues in a similar fashion.
Going through the arguments of the proof one verifies that the claimed uniformity statement holds.
\hspace*{\fill }$\square $

\medskip

To prove the asymptotic estimates (\ref{1.6}) and  (\ref{1.7}) of Theorem~\ref{thm1.5}, 
we will use the following version of the above formulas of the actions (cf \cite{FM}, \cite{HK})
   \begin{equation}
   \label{99.1} I^N_n = \frac{1}{\pi} \int ^{\lambda ^N_{2n}}_{\lambda ^N_{2n-1}}
                 {\rm arcosh} \big( (-1)^{N-n} \frac{\Delta _N(\mu )}{2} \big) d\mu \quad
                 \forall \, 1 \leq n \leq N - 1 .
   \end{equation}
The action variables of KdV can be expressed
in a similar way (cf \cite{FM})

   \begin{equation}
   \label{99.2} I^\pm _n = \frac{2}{\pi } \int ^{\lambda ^\pm _{2n}}_{\lambda ^\pm
                 _{2n-1}} {\rm arcosh} \big( (-1)^n \frac{\Delta _\pm (\lambda )}{2}
                 \big) d\lambda \ \quad \forall \, n\geq 1.
   \end{equation}
In \cite{KP} one finds a detailed derivation of formula \eqref{99.2} for $1$-periodic
potentials. As indicated in \cite{FM}, the formula remains valid for potentials of any
given period.

\medskip

{\it Proof of Theorem~\ref{thm1.5}:} 
%Let $L=[F(M)]$.
In view of Proposition \ref{Proposition C.1} it remains to prove \eqref{1.6} and \eqref{1.7}.
As both estimates can be proven in a similar way, we concentrate on \eqref{1.6} only. 
So let $1 \le n \le M.$
By the change of variable $\mu = - 2 + \frac{1}{4N^2} \lambda $ in the integral
in the above formula for $I^N_n$ one gets
   \begin{equation}
   \label{99.3} 4N^2 I^N_n = \frac{1}{\pi} \int ^{\nu ^N_{2n}}_{\nu ^N_{2n-1}} {\rm arcosh}
               \big( \frac{(-1)^{N-n}}{2} \Delta _N(-2 + \frac{\lambda }{4N^2})\big)
               d\lambda
   \end{equation}
where for any $0 \leq j \leq 2N - 1$ we set
$\nu ^N_{j}= 4N^2(\lambda ^N_j + 2)$. By Theorem~\ref{Theorem 1.1}
one has for any $0 \leq j \leq 2M$,
   \begin{equation}
   \label{99.4} \nu ^N_j = \lambda ^-_j + O \big( \frac{M^2}{N} \big) .
   \end{equation}

Let us consider first the case where $\gamma^-_n = \lambda ^-_{2n} - \lambda ^-_{2n-1} > 0$. Then
with the change of variable
   \[ \lambda (x) = \nu ^N_{2n-1} + 
\frac{4N^2\gamma^N_n}{\gamma_n^-} \, x, \quad \,\, 0 \leq x \leq \gamma_n^-
   \]
one gets in view of the identity \eqref{99.1}
   \begin{equation}
   \label{99.5} 8N^2 I^N_n = 
		\frac{8N^2\gamma^N_n}{\gamma_n^-}I^-_n +
\frac{8N^2\gamma^N_n}{\gamma_n^-}\frac{1}{\pi } \int ^{\gamma_n^-}_0 f(x) dx
   \end{equation}
where 
\[
f(x) : = {\rm arcosh} \big( \frac{(-1)^{N-n}}{2} \Delta _N(-2 + \frac{\lambda(x) }{4N^2})\big)
- {\rm arcosh} \big( \frac{(-1)^n}{2} \Delta _- (\lambda ^- _{2n-1} + x) \big).
\]

Note that by Theorem~\ref{Theorem9.1} one has
\[
 \frac{(-1)^{N-n}}{2} \Delta _N(-2 + \frac{\lambda(x) }{4N^2})
 =  \frac{(-1)^n}{2} \Delta _-(\lambda (x)) +O\big( \frac{F(M)^2}{M} \big). 
\]
                  
By \eqref{99.4}, the first term on the right hand side of \eqref{99.5} can be estimated as
\[
\frac{8N^2\gamma^N_n}{\gamma_n^-} I^-_n = 
     I^-_n + O \big( \frac{I^-_n}{\gamma^-_n}  \frac{M^2}{N} \big) .
   \]
As $I^-_n = O( \frac{1}{n} (\gamma^-_{n})^2)$ (cf
\cite[Theorem 7.3]{KP}) one concludes that
   \begin{equation}
   \label{99.6} 
\frac{8N^2\gamma^N_n}{\gamma_n^-} I^-_n  = I^-_n + 
 O \big( \gamma^-_{n} \frac{M^2} {nN} \big) .
   \end{equation}
The second term on the right hand side of \eqref{99.5} is estimated by using
that ${\rm arcosh}$ is $\frac{1}{2}-$H\"older continuous (cf Lemma~\ref{Hoelder}).
Indeed, together with the fact that $\Delta _-(\lambda )$ is Lipschitz continuous 
and therefore
   \[ \Delta _-(\lambda (x)) - \Delta _-(\lambda ^-_{2n-1}+x)
       = O(\lambda (x) - \lambda ^-_{2n-1} - x) = O \big( \frac{M^2} {N} \big)
   \]
one gets
   \[ f(x) = O \big( \frac{F(M)}{M^{1/2}} + \frac{M}{N^{1/2}} \big) .
   \]
Hence
\begin{equation}
   \label{99.7}  
\frac{8N^2\gamma^N_n}{\gamma_n^-} \frac{1}{\pi } \int ^{\gamma_n^-} _0 f(x) dx 
   =O \big( (\gamma_n^- +  \frac {M^2}{N}) (\frac {F(M)}{M^{1/2}} + \frac{M}{N^{1/2}} )\big).
   \end{equation}
Substituting \eqref{99.6} - \eqref{99.7} into \eqref{99.5} leads to
   \[ 8N^2 I^N_n = I^-_n + O\big( (\frac{M^{2}}{N} +\gamma_n^-) 
     (\frac {F(M)}{M^{1/2}} + \frac{M}{N^{1/2}} ) \big)
   \]
as claimed. Since the estimates are uniform with respect to the size of the gap length
$\gamma^-_n$ one can use an approximation argument to extend the result to the case
where $\gamma^-_n = 0$.
Going through the arguments of the proof one verifies that the claimed uniformity statement holds.
\hspace*{\fill }$\square $

\bigskip

To prove the asymptotics of the Toda frequencies stated in Theorem~\ref{thm1.4} we will need in addition the
asymptotics of the quantities $J^N_n, 0 < n < N$, given by
   \[ J^N_n:= \begin{cases} I^N_n / \gamma ^N_n &\mbox{ if } \gamma ^N_n >
      0 \\ 0 &\mbox{ if } \gamma ^N_n = 0 \end{cases}
   \]
in terms of $J^\pm _n, n \geq 1$, given by
   \[ J^\pm_n:= \begin{cases} I^\pm_n / (2\gamma ^\pm_n) &\mbox{ if } \gamma ^
      \pm_n > 0 \\ 0 &\mbox{ if } \gamma ^\pm_n = 0 \end{cases} .
   \]

\begin{guess}
\label{Theorem 9.1} Uniformly in $1 \leq n \leq F(M)$ and uniformly on
bounded subsets of $\alpha , \beta $ in $C^2_0({\mathbb T})$,
  \[ J^N_n - J^-_n,\,\,\,  J^N_{N-n} - J^+_{n} = O \Big( \frac{M}{\sqrt{N}} +
     \frac{F(M)}{\sqrt{M}} \Big) .
  \]
\end{guess}

\begin{remark}
\label{Remark1} The proof of Theorem~\ref{Theorem 9.1} will show that both
$(J^N_n)_{1 \leq n \leq M}$ and $(J^-_n)_{1 \leq n \leq M}$ respectively
$(J^N_{N-n})_{1 \leq n \leq M}, (J^+_n)_{1 \leq n \leq M}$ are uniformly
bounded with respect to $n$ and with respect to bounded subsets of functions
$\alpha , \beta $ in $C^2_0({\mathbb T})$.
\end{remark}

\begin{remark}
\label{Remark2} In \cite{HK1} and \cite{KP} it is established that for any fixed $N$, $I^N_n =
O((\gamma ^N_n)^2)$  and $I^\pm_n = O((\gamma ^\pm _n)^2)$ uniformly on
closed bounded subsets of vectors $(b, a)$ in 
${\mathbb R}^N \times {\mathbb R}_{>0}^N$ respectively of
potentials in $L^2({\mathbb T})$. Theorem~\ref{Theorem 9.1} states
that the estimates for $(J^N_n)_{1 \leq n \leq M}, (J^N_{N-n})_{1 \leq n \leq M}$ 
are uniform in $N$.
\end{remark}

{\it Proof:} The estimates for $J^N_n - J^-_n$ and $J^N_{N-n} - J^+_n$
are proven in a similar way and so we concentrate on the first one only.
We argue similarly as in the proof of Theorem~\ref{Theorem 9.1}. Let
$1 \leq n \leq F(M)$. First consider the case where $\gamma ^-_n > 0$.
With the change of variable $\lambda = \lambda ^-_{2n-1} + t \gamma ^-_n$,
one gets
   \begin{align*} I^-_n &= \frac{2}{\pi } \int ^{\lambda ^-_{2n}}_{\lambda
                     ^-_{2n-1}} \mbox{arcosh} \big( (-1)^n \frac{\Delta _-(\lambda )}
                     {2} \big) d\lambda \\
                  &= 2 \gamma ^-_n \frac{1}{\pi } \int ^1_0 \mbox{arcosh} \big(
                     \frac{(-1)^n}{2} \Delta _-(\lambda ^-_{2n-1} + t\gamma ^-
                     _n) \big) dt .
   \end{align*}
Note that the assumption $\gamma ^-_n > 0$ together with the asymptotics
$\gamma ^N_n = \frac{\gamma ^-_n}{4N^2} + O \left( \frac{M^2}{N^3} \right) $
 implies that for $N$ sufficiently large one has
$\gamma ^N_n > 0$. Hence with the change of variable $\mu = \lambda
^N_{2n-1} + t \gamma ^N_n$ one gets
   \[ I^N_n = \gamma ^N_n \frac{1}{\pi } \int ^1_0 \mbox{arcosh} \big(
      \frac{(-1)^{N-n}}{2} \Delta _N(\lambda ^N_{2n-1} + t \gamma ^N_n)\big)dt .
   \]
We then conclude that
   \begin{align*} \frac{I^N_n}{\gamma ^N_n} - \frac{I^-_n}{2\gamma ^-_n} =
                     \frac{1}{\pi } \int ^1_0 \Big( &\mbox{ arcosh} \big(
                     \frac{(-1)^{N+n}}{2} \Delta _N \left( \lambda ^N_{2n-1}
                     + t \gamma ^N_n \right) \big)  \\
                  &- \mbox{arcosh} \left( \frac{(-1)^n}{2} \Delta _- \big(
                     \lambda ^-_{2n-1} + t \gamma ^-_n \big) \right) \Big)dt
   \end{align*}   
and as $\mbox{arcosh}$ is $\frac{1}{2}$-H\"older continuous (cf Lemma~\ref{Hoelder}) one then has
   \[  \frac{I^N_n}{\gamma ^N_n} - \frac{I^-_n}{2\gamma ^-_n} = O \big( \sup _{0 \leq t \leq 1} \big\arrowvert (-1)^N
      \Delta _N(\lambda ^N_{2n-1} + t\gamma ^N_n) - 
      \Delta _- (\lambda ^-_{2n-1} + t\gamma ^-_n) \big\arrowvert ^{1/2} \big).
   \]
As 
$\lambda ^N_{2n-1} + t\gamma ^N_n = -2 + 
(\lambda^-_{2n-1} + t\gamma_n^- + O( \frac{M^2}{N} )  )\frac{1}{4N^2}  $
Theorem~\ref{Theorem9.1} implies that
   \begin{align*} &\big\arrowvert (-1)^N \Delta _N(\lambda ^N_{2n-1} + t
                     \gamma ^N_n) - \Delta _-(\lambda ^-_{2n-1} + t\gamma^-_n) \big\arrowvert \\
                  &= \big\arrowvert \Delta _-(\lambda ^-_{2n-1} + t\gamma ^-_n
                     +O( \frac{M^2}{N} )) - \Delta _-(\lambda ^-
                     _{2n-1} + t \gamma ^-_n) \big\arrowvert 
                     + O \Big( \frac{F(M)^2}{M} + \frac{M^2}{N} \Big)
   \end{align*}
and by the Lipschitz continuity of $\Delta_-$ it then follows that
   \[ \Big\arrowvert \frac{I^N_n}{\gamma ^N_n} - \frac{I^-_n}{2\gamma ^-_n}
      \Big\arrowvert = O \Big(  \frac{F(M)}{\sqrt{M}}  + \frac{M}{\sqrt{N}} \Big).
   \]
Since the estimates are uniform with respect to the size of the gap length
$\gamma^-_n$ one can use an approximation argument to extend the result to the case
where $\gamma^-_n = 0$.
Going through the arguments of the proof one verifies that the claimed uniformity holds.
\hspace*{\fill }$\square $

%%%%%%%%%%%%%%%%%%%%%%%%%%%%%%%%%%%%%%%%%%%%%%%%%%%%%%%%%%%%%%%%%%%%%%%%%
%%%%%%%%%%%%%%%%%%%%%%%%%%%%%%%%%%%%%%%%%%%%%%%%%%%%%%%%%%%%%%%%%%%%%%%%%
%%%%%%%%%%%%%%%%%%%%%%%%%%%%%%%%%%%%%%%%%%%%%%%%%%%%%%%%%%%%%%%%%%%%%%%%%

\section{Formulas for the frequencies}
\label{frequencies}

In this section we prove Theorem~\ref{lemnouv} and review formulas expressing the frequencies of the periodic Toda chain and the KdV equation as periods of certain differentials of cycles on the corrsponding spectral curves.

\medskip

%%%%%%%%%%%%%%%%%%%%%%%%%%%%%%%%%%%%%%%%%%%%%%%%%%%%%%%%%%%%%%%%%%%%%%%%

{\em Spectral curves.}
%\label{form}
The Toda curve $\mathcal C_N \equiv \mathcal C_N^{(\alpha, \beta)}$ corrsponding to the
state $(b^N,a^N)$ is defined as the affine curve 
\[
\mathcal C_N = \{ (\mu, z) \in \mathbb C^2 : z^2 = \Delta^2_N (\mu) - 4\}.
\]
In case the spectrum of $Q_N^{\alpha, \beta}$ is simple, it is a two sheeted open Riemann surface
whose ramification points are the eigenvalues $(\lambda^N_n)_{0 \le n \le 2N-1}$
of $Q_N^{\alpha, \beta}$. 
If $\lambda^N_{2n}$ is a double eigenvalue of $Q_N^{\alpha, \beta}$, then 
$(\lambda^N_{2n}, 0)$  is a singular point of $\mathcal C_N$. 
On $\mathcal C_N$ we define a set of one forms $\frac{\psi^N_n(\mu)}{\sqrt{\Delta^2_N (\mu) - 4}}d\mu$
where $\psi^N_n(\mu)$ are polynomials of degree $N-2$, uniquely determined by the normalisation
conditions
\[ \frac{1}{2\pi } \int _{\Gamma ^N_k} \frac{\psi ^N_n(\mu )}{\sqrt[c]{\Delta^2
      _N(\mu ) - 4}} d\mu = \delta _{nk} \quad \forall 1 \leq k \leq N - 1 .
\]
Here $\sqrt[c]{\Delta _N(\mu )^2 - 4}$ denotes the canonical root of 
   \be\label{rn} \Delta ^2_N(\mu ) - 4 = \frak q^{-2}_N \chi_N(\mu ), \quad
\frak q_N:=\prod\limits_{i=1}^N{a^N_i},
   \ee
and $(\Gamma ^N_k)_{ 1 \leq k \leq N - 1}$ the counterclockwise
oriented contours
introduced at the end of Section \ref{intro}.
In view of \eqref{rn} it is convenient to write $\psi_n^N$ in the form
   \[ \psi ^N_n (\mu ) = M^N_n q^{-1}_N \varphi ^N_n(\mu ) , \qquad  \varphi ^N_n(\mu ) =
      \prod _{\underset{k\not=n}{0< k < N}} (\mu - \sigma ^{N,n}_k ) 
	= (-1)^N \prod _{\underset{k\not=n}{0< k < N}} (\sigma ^{N,n}_k - \mu)  .
   \]
Note that if the spectrum of $Q_N^{\alpha, \beta}$ is simple
$(\frac{\psi^N_n(\mu)}{\sqrt{\Delta^2_N (\mu) - 4}}d\mu)_{0 < n < N}$ are
linearly independent holomorphic differentials on the Riemann surface $\mathcal C_N.$ 
If $\lambda^N_{2n}$ is a double eigenvalue,
$\frac{\psi^N_n(\mu)}{\sqrt{\Delta^2_N (\mu) - 4}}d\mu$ has a pole of order $1$
at $(\lambda^N_{2n}, 0).$

For the KdV equation, similar objects have been introduced. Denote by $\mathcal C_{\pm}$
the curve corresponding to $q_{\pm}$,
\[
\mathcal C_\pm = \{ (\lambda, w) \in \mathbb C^2 : w^2 = \Delta^2_\pm (\lambda) - 4\}.
\]
If the spectrum of the Hill operator $H_{\pm}$ is simple, $\mathcal C_\pm$ is a two sheeted open Riemann surface of infinite genus whose ramification points are the eigenvalues 
$(\lambda^{\pm}_n)_{n \ge 0}$ of $H_{\pm}.$
On $\mathcal C_\pm$ we define a set of one forms $\frac{\psi^\pm_n(\lambda)}{\sqrt{\Delta^2_\pm (\lambda) - 4}}d\lambda$.
Note that $\Delta_\pm(\lambda)^2-4$ admits the product representation 
(cf \cite[Appendix A]{bkp2})
\begin{equation}
\label{Deltapro}
\Delta^2_\pm(\lambda)-4 = (\lambda_0^\pm -  \lambda)
\prod_{n\geq 1}\frac{(\lambda_{2n}^{\pm}-\lambda)(\lambda_{2n-1}^{\pm}-\lambda)}{(4\pi^2n^2)^2},
\end{equation}
whereas $\psi_n^\pm(\lambda),\ n\geq 1,$ are entire functions of $\lambda$ of the form
\[
\psi_n^\pm(\lambda)=\frac1 {2\pi n}\prod_{k\geq 1,k\neq n}\frac{\sigma_k^{\pm,n}-\lambda}
{4\pi^2 k^2}
\]
constructed in \cite[Proposition D.10]{KP} for potentials of period $1$ - see \cite[Appendix A]{bkp2} for the adjustments needed in the case of potentials of period $1/2.$
Its zeros $\sigma_k^{\pm,n}$ are determined by the conditions
\begin{equation}
\label{psi-}
\int_{\Gamma_k^\pm}\frac{\psi_n^\pm(\lambda)}{\sqrt[c]{\Delta^2_\pm
(\lambda)-4}} d\lambda=0 \qquad \forall \,\, k \geq 1,\ k\neq n
\end{equation}
and the normalisation factor $\frac 1 {2\pi n}$ by the requirement that
\begin{equation}
\label{psi+}
\int_{\Gamma_n^\pm}\frac{\psi_n^\pm(\lambda)}{\sqrt[c]{\Delta^2_\pm
(\lambda)-4}} d\lambda=1.
\end{equation}
Here $(\Gamma_k^\pm)_{ k\geq 1}$ denote pairwise disjoint counterclockwise oriented
contours chosen in such a way that $\lambda^\pm_{2k-1},\lambda^\pm_{2k}$ are
inside of $\Gamma^\pm_k$ whereas all other eigenvalues are outside and $\sqrt[c]{\Delta^2_\pm(\lambda)-4}$ is the canonical root of 
$\Delta_\pm(\lambda)^2-4$ determined by the sign condition
$i\sqrt[c]{\Delta^2_\pm(\lambda)-4} \,>\, 0$ for $\lambda_0^\pm < \lambda <\lambda^\pm_{1}$ or  
\[
\sqrt[c]{\Delta^2_\pm(\lambda +i0)-4} \,>\, 0 \qquad \forall \,\, \lambda <\lambda^\pm_{0}.
\]
%The version of the identity \eqref{psi+} for potentials $q$ of period $1$ is proven 
%in \cite[Proposition D.10]{KP}. Going through the arguments in \cite{KP} for potentials 
%$q$ of period $T$ one sees that
%the normalisation factor of the entire function $\psi_n(\lambda)$, satisfying 
%\eqref{psi-} and \eqref{psi+}  is independent of $q$. Denote it by $c_n^T$.
%The computations at the end of \cite[Appendix A]{bkp2} show that for the zero
%potential, viewed as $T$-periodic potential, $c_n^T=\frac{2T^2}{n\pi}$, leading 
%to the claimed value $\frac{1}{2n\pi}$ in the case $T=1/2$.

\medskip 

%%%%%%%%%%%%%%%%%%%%%%%%%%%%%%%%%%%%%%%%%%%%%%%%%%%%%%%%%%%%%%%%%%%%%%%

{\em Toda frequencies.} 
As a first step we prove Theorem \ref{lemnouv}.

 \medskip

\noindent {\it Proof of Theorem \ref{lemnouv}.} Let $0 < n < N.$
Using a continuity argument, it suffices to prove the claimed identity in the case
where $\lambda ^N_{2n-1} < \lambda ^N_{2n}.$
By the definitions of $\psi ^N_n$ and $\varphi ^N_n$ one 
has
\begin{equation}
   \label{B.2} \frac{\psi ^N_n(\mu )}{\sqrt{\Delta ^2_N(\mu ) - 4}} = M^N_n
               \frac{\varphi ^N_n(\mu )}{\sqrt{\chi_N(\mu )}},
   \end{equation}
implying that the normalisation factor $M^N_n$ is given by
\begin{equation}\label{normalisation}
 M^N_n \cdot \frac{1}{\pi } \int _{\lambda^{N}_{2n-1}}^{\lambda^{N}_{2n}} \frac{\varphi ^N_n(\mu )}
   {\sqrt[c]{\chi_N(\mu -i0 )}} d\mu = 1.
\end{equation}
To see that $M^N_n >0 $ note that 
$ \lambda ^N_{2k-1} \leq \sigma ^{N,n}_k \leq \lambda ^N_{2k}$, $0 < k < N,$ and thus
   \[ (-1)^{N - 1 - n} \varphi ^N_n(\mu ) = \prod _{k<n} (\mu - \sigma ^{N,n} _k) 
      \cdot \prod _{k>n} (\sigma ^{N,n}_k - \mu ) > 0 \quad
            \forall \,\, \lambda ^N_{2n-1} \leq \mu \leq \lambda ^N_{2n}.
   \]
On the other hand, if $\lambda ^N_{2n-1} < \lambda ^N_{2n}$ one has by the definition of the $c$-root
   \[ \sqrt[c]{\chi_N(\mu - i0)} = (-1)^{N+1-n} \sqrt[+]{\chi_N(\mu)} \qquad \forall \,\,\lambda
      ^N_{2n-1} < \mu < \lambda ^N_{2n}\ .
   \]
Altogether it follows that $M^N_n > 0$. Going through the arguments of \cite{HK4}, but without assuming that the trace of $L(b^N, a^N)$ vanishes, one concludes from \cite[Lemma 4.4, formulas (28), (33), (38)]{HK4} that 
$|\omega ^N_n| = M^N_n.$ In particular, it follows that $\omega ^N_n \ne 0$. Taking into account that at the equilibrium $(0_N, 1_N),$
$\omega ^N_n = 2 \sin(\frac{n\pi}{N}) > 0$, a deformation arguments yields the claimed identity
$\omega ^N_n = M^N_n$ and positivity $\omega ^N_n > 0.$ \qed

\medskip

\begin{remark}
\label{cauch}
Assume that $\lambda ^N_{2n-1} < \lambda ^N_{2n}$. 
 Applying the mean value theorem (cf Lemma~\ref{atif}) to \eqref{normalisation} it follows that there is
$\lambda ^N_{2n-1} \leq \mu _\ast \leq \lambda ^N_{2n}$ such that
   \begin{equation}
   \label{B.3} \omega ^N_n = \sqrt[+]{(\lambda ^N_{2N-1}-\mu _\ast )(\mu _\ast -
               \lambda ^N_0)} \prod _{\underset{1 \leq k < N}{k \not= n}}
               \frac{\sqrt[+]{(\lambda ^N_{2k} - \mu _\ast )(\lambda ^N_{2k - 1}
               - \mu _\ast )}}{|\sigma ^{N,n}_k - \mu _\ast |} .
   \end{equation}
%In particular, $\omega ^N_n > 0$.
By a continuity argument, \eqref{B.3} continues to hold in the case $\lambda ^N_{2n-1}
= \lambda ^N_{2n}$ with $\mu _\ast = \tau ^N_n$ where
  \be\label{tauuuu}  \tau^N_n = (\lambda ^N_{2n} + \lambda ^N_{2n-1}) / 2.
  \ee
\end{remark}
\begin{remark}
\label{lin}
Actually, formula \eqref{B.3} is valid for any state $(b,a) \in \mathbb R^N \times \mathbb R^N_{>0}.$ It allows to compute the frequencies corresponding to the equilibrium states
\[ b = r 1_N,  \quad \ a = s1_N \quad \mbox{with} \quad 1_N = (1, \ldots , 1), 
\,\,\, r \in \mathbb R,\,\,\, s>0.
\]
By \cite[Lemma 2.6]{HK1} one has $\lambda ^N_0 = - 2 s + r,$ $ \lambda ^N_{2N - 1} = 2s+ r$ and
   \[ \lambda ^N_{2n} = \lambda ^N_{2n - 1} = -
      2s \cos \frac{n\pi }{N} + r \quad  \forall \,\, 0 < n < N.
   \]
Hence $\sigma ^{N,n}_k = - 2s \cos \frac{k\pi }{N} + r$ for any $1 \leq n, k \leq
N - 1$ and thus in this case, for any $0 < n < N $,
   \[ \omega ^N_n = 2s \sqrt[+]{(1 + \cos \frac{n\pi  }{N})(- \cos \frac{n\pi }
      {N} + 1)} = 2s \cdot \sin \frac{n\pi }{N} .
   \]
\end{remark}

\medskip

Formula \eqref{B.3} is suited for deriving asymptotic estimates of
$\omega ^N_n$ as $N \rightarrow \infty $ for $n$ in the range $M < n <
N - M$. However to get asymptotic estimates for the remaining frequencies
we need alternative formulas, obtained with the help of Riemann
bilinear relations.
Such formulas were derived in \cite{HK4}, Section 4,
for trace free Jacobi matrices. In addition it is proved there that they are independent of the trace. 
Hence these formulas hold in general modulo a spectral shift. It is however more convenient
for us to have formulas which do not involve such a shift.

\medskip

\begin{proposition}\label{freqtod}
The frequencies $(\omega ^N_n)_{ \ 0 < n < N}$ of the state $(b^N,a^N)$ satisfy 
 \begin{equation}
   \label{B.5}
 iN\omega ^N_n = \sum ^n_{j=1} \int ^{\lambda ^N_{2j-1}}
                     _{\lambda ^N_{2j-2}} \frac{(\mu - \frac{1}{N}\frak p_N)\dot \Delta _N
                     (\mu )d\mu}{\sqrt[c]{\Delta^2_N(\mu ) - 4}} 
                  - \sum _{k \in \mathcal J_N} I^N_k \omega ^N_k \sum ^n_{j=1}
                     \int ^{\lambda ^N_{2j-1}}_{\lambda ^N_{2j-2}}
                     \frac{\varphi ^N_k(\mu )d\mu}{\sqrt[c]{\chi_N(\mu )}}
   \end{equation}
where $\mathcal J_N:= \{ 1 \leq k \leq N - 1 | I^N_k \not= 0 \}$ and 
$\frak p_N= \sum_{i=1}^{N}b^N_i.$
\end{proposition}

\begin{remark}\label{quotient}
In the formula \eqref{B.5}, the quotients, $\frac{\dot \Delta _N(\mu )}{\sqrt[c]{\Delta^2_N(\mu ) -
4}}$ and $\frac{\varphi ^N_k(\mu )}{\sqrt[c]{\chi_N(\mu )}},\ k\in \mathcal J_N$, are desingularized
in the sense that factors which appear both in the nominator and in the denominator
are canceled. Both denominators then consist of the square root of a polynomial
with simple roots.
\end{remark}

\proof In the case where all the eigenvalues
$(\lambda ^N_j)_{0 \leq j \leq 2N-1}$ are simple such formulas have been computed
in \cite[Theorem 4.5]{HK4}. Here we will use a slightly modified version, not
involving a translation of the vector $b^N$. 
Going through the arguments \cite[Section 4]{HK4} -- in particular Theorem 4.5 and Lemma 4.2 --
one sees that if $(\lambda ^N_j)_{0 \leq j \leq 2N -
1}$ are simple, one has for any $1 < n <  N$
   \begin{equation}
   \label{B.4} 
iN \omega ^N_n =\int ^{\lambda ^N_{2n-1}}_{\lambda ^N_0}
   \frac{(\mu - \frac{1}{N} \frak p_N) \dot \Delta _N(\mu )d\mu}{\sqrt[c]{\Delta^2_N (\mu - i0) - 4}}
  - \sum ^{N - 1}_{k=1} I^N_k \omega ^N_k
                  \int ^{\lambda ^N_{2n-1}}_{\lambda ^N_0} \frac{\varphi ^N_k(\mu )d \mu}
                  {\sqrt[c]{\chi_N(\mu - i0)}}.
   \end{equation}
Here we used Theorem~\ref{lemnouv} and formula \eqref{B.2} together with the fact that the one form
   \[ \frac{(\mu - \frac{1}{N} \frak p_N) \dot \Delta _N(\mu )}{\sqrt{\Delta^2_N(\mu ) - 4}}
      d\mu
   \]
has an expansion in $\mu = \frac{1}{z}$ at $\infty^{\pm} $ of the form $(\epsilon_{\pm} \frac{N}
{z^2} + O(1)) dz$ with $\epsilon_{\pm} \in \{+ , -\}$ -- see \cite[Lemma 4.1 and Appendix A]{HK4}. The integrals in \eqref{B.4} are
split up into integrals over bands and gaps,
   \begin{equation}
   \label{B.4bis} \int ^{\lambda ^N_{2n-1}}_{\lambda ^N_0} = \sum ^n_{j=1} \int
                  ^{\lambda ^N_{2j-1}}_{\lambda ^N_{2j-2}} + \sum ^{n-1}_{j=1}
                  \int ^{\lambda ^N_{2j}} _{\lambda ^N_{2j-1}} .
   \end{equation}
Recall that $\int ^{\lambda ^N_{2j}}_{\lambda ^N_{2j-1}} \frac{\omega ^N_k
\varphi ^N_k(\mu)}{\sqrt[c]{\chi_N(\mu - i0)}} d\mu = \pi \delta _{jk}$ and, as
$\int ^{\lambda ^N_{2j}}_{\lambda ^N_{2j-1}} \frac{\dot \Delta _N(\mu)}
{\sqrt{\Delta^2_N(\mu )-4}} d\mu = 0$,
   \[ \int ^{\lambda ^N_{2j}}_{\lambda ^N_{2j-1}} \frac{(\mu - \frac{1}{N} \frak p_N)\dot
      \Delta _N(\mu)}{\sqrt[c]{\Delta^2_N(\mu -i0)-4}}d\mu = \pi I^N_j .
   \]
Hence the integrals in \eqref{B.4bis} over the gaps cancel and $iN\omega ^N_n$ equals
\[
\sum ^n_{j=1} \int ^{\lambda ^N_{2j-1}}
                     _{\lambda ^N_{2j-2}} \frac{(\mu - \frac{1}{N}\frak p_N)\dot \Delta _N
                     (\mu )}{\sqrt[c]{\Delta^2_N(\mu ) - 4}} d\mu 
                  - \sum _{k=1}^{N-1} I^N_k \omega ^N_k \sum ^n_{j=1}
                     \int ^{\lambda ^N_{2j-1}}_{\lambda ^N_{2j-2}}
                     \frac{\varphi ^N_k(\mu )}{\sqrt[c]{\chi_N(\mu )}} d\mu.
\]
In case some of the eigenvalues $(\lambda ^N_j)_{0 \leq j \leq 2N-1}$ are double
a continuity argument shows that \cite[Theorem 4.5]{HK4} continues to hold. Going through
the arguments above and using that $I_k^N = O((\gamma_k^N)^2)$ by  \cite[Proposition 3.6]{HK1}, 
a limiting argument shows that the claimed formula for the frequencies holds also in this case.
\qed

%%%%%%%%%%%%%%%%%%%%%%%%%%%%%%%%%%%%%%%%%%%%%%%%%%%%%%%%%%%%%%%%%%%%%%%%

\medskip

{\em KdV Frequencies.}
It turns out that there are formulas for the KdV frequencies which are of a very
similar form as the ones of the Toda chains. 

\begin{lemma}\label{scan0}
The KdV frequencies  of a periodic potential of period $T$ and mean $0$
are $\omega_{T,n}=-\frac{24}TW_n,\  n\geq 1$, where
\be\label{scan2}
W_n=\int ^{\lambda_{2n-1}}_{\lambda_0}
                 \frac{\lambda \dot \Delta (\lambda ) d \lambda}{i\sqrt[c]{\Delta^2(\lambda - i0) - 4}}  
                   - \sum _{k\geq 1, I_k\neq 0}
                   \frac{I_k}{2}\int ^{\lambda _{2n-1}}_{\lambda_0} \frac{\psi_n(\lambda) d \lambda}
                  {i\sqrt[c]{\Delta^2(\lambda - i0) - 4}}
		  \ee
where $(\lambda_n)_{n\geq 0}$ denote the periodic eigenvalues of $-\partial_x^2+q$ on
$[0,2T]$, $\Delta$ the discriminant of $-\partial_x^2+q$ on $[0,T]$, and
$I_k$ the action variables.
\end{lemma}

{\it Proof:} In the case  $T=1$, formula \eqref{scan2} has been proven in \cite[Lemma 2.5]{KP}. For
an arbitrary period we argue similarly, using that $\omega_{T,n}=K_TW_n,\ n\geq 1$ for
some constant $K_T > 0$  which can be computed by following the
arguments in \cite{KP} for the case $T = 1$: according to the computations at
the end of \cite[Appendix A]{bkp2}, the zero potential, viewed as potential
of period $T$ has discriminant
$\Delta(\lambda)=2\cos({\sqrt\lambda T})$. By the definition of the c-root one
has for $\lambda \ge 0,$ $\sqrt[c]{\Delta^2(\lambda)-4}=-2i\sin({\sqrt[+]\lambda T})$
and the periodic eigenvalues are
  \[
  \lambda_0=0, \qquad 
\lambda_{2n}=\lambda_{2n-1}=\frac{n^2\pi^2}{T^2} \quad \forall \,\, n\geq 1
  \]
and, for $\lambda \ge 0$, $\dot\Delta(\lambda)=-\frac T{\sqrt[+]\lambda}\sin{\sqrt[+]\lambda T}$. As $I_k=0$ for
any $k\geq 1$ it then follows that
\[
W_n=\int ^{\lambda_{2n-1}}_{\lambda_0}
                 \frac{\lambda \dot \Delta (\lambda ) d \lambda}
                 {i\sqrt[c]{\Delta^2(\lambda - i0) - 4}}=-\frac T3
                 \lambda^{\frac 3 2}|_0^{(\frac{n\pi}T)^2}
		=-\frac T3(\frac{n\pi}T)^3. 
\]
Arguing as in \cite{KP}, the frequencies of the
zero potential on $[0,T]$ can be computed as $\omega_{T,n}=(\frac{2n\pi}T)^3$.
Hence $\omega_{T,n}=-\frac{24}TW_n$ as claimed.
\hspace*{\fill }$\square $

\medskip

Applying Lemma~\ref{scan0} in the case $T=\frac12$ and $q=q_\pm$ and arguing as 
in the proof of Proposition~\ref{freqtod} leads to the following formulas.

\medskip

\begin{proposition}\label{scan1}
The KdV frequencies $\omega ^\pm_n,\ n\geq 1$, of the potentials $q_\pm$ of period
$1/2$ and mean $0$ are 
\begin{equation}\label{B.5KdV}
%		  \omega ^\pm_n = 
-48\sum ^n_{j=1} \int ^{\lambda ^\pm_{2j-1}}
                     _{\lambda ^\pm_{2j-2}} \frac{\lambda\dot \Delta _\pm
                     (\lambda )d \lambda}{i \sqrt[c]{\Delta^2_\pm(\lambda ) - 4}} \quad
                  +24 \sum_{k=1,I^{\pm}_k\neq 0}^\infty I^\pm_k \sum ^n_{j=1}
                     \int ^{\lambda ^\pm_{2j-1}}_{\lambda ^\pm_{2j-2}}
                     \frac{\psi ^\pm_n(\lambda ) d\lambda}{i\sqrt[c]{\Delta^2_\pm(\lambda ) - 4}}
\end{equation}	
\end{proposition}

\vskip 1cm
%The fact that in formula \eqref{B.5}, the frequencies also
%appear on the right hand side makes it necessary to analyze the Toda frequencies
%in several steps.
%\hspace*{\fill }$\square $

%%%%%%%%%%%%%%%%%%%%%%%%%%%%%%%%%%%%%%%%%%%%%%%%%%%%%%%%%%%%%%%%%%%%%%%
%%%%%%%%%%%%%%%%%%%%%%%%%%%%%%%%%%%%%%%%%%%%%%%%%%%%%%%%%%%%%%%%%%%%%%%%
%%%%%%%%%%%%%%%%%%%%%%%%%%%%%%%%%%%%%%%%%%%%%%%%%%%%%%%%%%%%%%%%%%%%%%%%

\section{Asymptotics of the zeroes of $\varphi ^N_k$}
\label{zeroes}
In this section we compute the asymptotics of the zeroes $\sigma_\ell^{N,k}$ of the 
polynomials $\varphi ^N_\ell$, introduced in Section~\ref{intro}, in terms of the zeroes $\sigma_\ell^{\pm,k}$ of the entire functions
$\psi ^\pm_k$ (cf Section \ref{frequencies}). 
To define the edges of the spectrum of $Q_N^{\alpha, \beta},$ we use an arbitrary
 function
$F: \mathbb N \to \mathbb R_{\ge 1}$, satisfying $(F)$ with $\eta \le 1/3.$ 
For the remainder of the paper we set $M=[F(N)]$ and $L=[F(M)].$  
Then $L^3 \le F(M)^3 \le M$ and $M \le \frac{N}{M^2}$. 
Furthermore, as $L \ge 1$, $M^3 \le N L^2$ or $\frac{M}{L^2} \le \frac{N}{M^2}.$
We record for reference the following inequalities
\begin{equation}
\label{Fwith1/3}
\frac{M^2}{N} \le L^{-3} \le L^{-2}, \qquad \frac{M^2}{N} \le  \frac{L^2}{M}, \qquad 
\frac{M}{N^{1/2}} \le \frac{L}{M^{1/2}}
\end{equation}
%%%%%%%%%%%%%%%%%%%%%%%%%%%%%%%%%%%%%%%%%%%%%%%%%%%%%%%%%%%%%%%%%%%%%
\begin{guess}
\label{theorem2.1} 
 Let $F : {\mathbb N} \rightarrow {\mathbb R}_{\geq 1}$ 
satisfy (F) with $\eta \le 1/3$.
For $1 \leq k, \ell
\leq L=[F(M)]$ with $\ell \not= k$
\[ 
{\rm (i)}\quad \sigma ^{N,k}_\ell 
= - 2 + \frac{\sigma ^{-,k}_\ell }{4N^2} + O \Big( \frac{1}{N^2L^{5/2}} \Big) \qquad 
{\rm (ii)}\quad \sigma ^{N,k}_{N-\ell } = 2 - \frac{\sigma ^{+,k}_\ell }{4N^2} 
+ O \Big( \frac{1}{N^2L^{5/2}} \Big).
\]  
The error estimate is uniform in $k, \ell $ and uniform on bounded subsets of functions
$\alpha , \beta $ in $C^2_0 ({\mathbb T})$.
\end{guess}

To prove Theorem \ref{theorem2.1} we first need to introduce additional notation and
derive auxilary results. It turns out to be convenient to define the following roots:
for $0 < \ell < N$  
denote by $w^N_\ell(\mu)$ the standard root
\begin{equation}
\label{wN}
w^N_\ell(\mu) = \sqrt[c]{( \lambda^N_{2\ell} - \mu)( \lambda^N_{2\ell -1} - \mu )},
\qquad \mu \in \mathbb C \setminus [\lambda^N_{2\ell -1}, \lambda^N_{2\ell}],
\end{equation}
determined by the sign condition $w^N_\ell(\mu) < 0$ for any $ \mu > \lambda^N_{2\ell}.$
The standard root is related to the canoncial root, introduced in \eqref{cNroot} as follows
\[
\sqrt[c]{\Delta^2_N(\mu - i0) - 4} 
=  \frak q_N^{-1} \sqrt[+]{\mu - \lambda^N_0} \sqrt[+]{\lambda^N_{2N-1} - \mu} \,
(-1)^N \prod_{0 < \ell < N} w_\ell^N(\lambda ) \quad \forall\,\, \mu > \lambda ^N_{2N - 1},
\]
leading to the formula 
$\frac{\psi^N_k(\mu)}{\sqrt[c]{\Delta^2_N(\mu - i0)-4}}=
\omega^N_k \frac{\varphi^N_k(\mu)}{\sqrt[c]{\chi_N(\mu - i0)}}\,$ where
\begin{equation}
\label{phidelta}
\frac{\varphi^N_k(\mu)}{\sqrt[c]{\chi_N(\mu - i0)}}=
\frac{1}{\sqrt[+]{\mu - \lambda^N_0} \sqrt[+]{ \lambda^N_{2N-1} - \mu} }\frac{1}{w_k^N(\mu)}
\prod_{\ell\not=k}\frac{ \sigma^{N,k}_\ell -\mu}{w_\ell^N(\mu)}.
\end{equation}
For later reference we record that
\begin{equation}
\label{signChiN}
(-1)^N i \sqrt[c]{\chi_N(\mu)} > 0 \quad \forall \mu \,\, \in (\lambda^N_0, \lambda^N_1).
\end{equation}
Similarly, introduce the standard root
\begin{equation}
\label{w-}
w^\pm_\ell(\lambda):=\sqrt[s]{(\lambda^\pm_{2\ell} - \lambda)(\lambda^\pm_{2\ell-1} - \lambda)},
\qquad \lambda \in \mathbb C \setminus [\lambda^\pm_{2\ell -1}, \lambda^\pm_{2\ell}],
\end{equation}
with the sign condition
$w^\pm_\ell(\lambda) < 0$ for any $ \lambda > \lambda^N_{2\ell }.$
The canonical root $\sqrt[c]{\Delta^2_\pm(\lambda)-4}$ is related to the standard roots by
\begin{equation}
\label{wdelta}
\sqrt[c]{\Delta^2_\pm(\lambda + i0)-4}=\sqrt[+]{\lambda_0^\pm - \lambda}\prod_{\ell\geq
  1}\frac{w_\ell^\pm(\lambda )}{4\pi_\ell^2} \qquad \forall \,\, \lambda < \lambda^{\pm}_0
\end{equation}
where $\pi_\ell=\ell \pi$ for $\ell\neq 0$ and $\pi_0=1$, yielding for $\lambda < \lambda_0^\pm$
\begin{equation}
\label{psidelta}
\frac{\psi^\pm_k(\lambda)}{\sqrt[c]{\Delta^2_\pm(\lambda)-4}}=
\frac{2k\pi}{\sqrt[+]{\lambda_0^\pm - \lambda }} \frac{1}{w_k^\pm(\lambda)}
\prod_{\ell\not=k}\frac{\sigma^{\pm,k}_\ell-\lambda}{w_\ell^\pm(\lambda)}.
\end{equation}
Again for later reference we record that
\begin{equation}
\label{signDelta-}
 i \sqrt[c]{\Delta^2_-(\lambda) -4} > 0 \quad \forall \lambda \,\, \in (\lambda^-_0, \lambda^-_1).
\end{equation}
As items (i) and (ii) are proved in the same way we concentrate on (i)
only. 

For the sequel it turns out to be convenient to choose the contours $\Gamma
^-_n$ in \eqref{psi-} as
follows. Choose $\rho > 0$ in such a way that
   \[ \lambda ^-_{2\ell } + 3\rho < \lambda ^-_{2\ell + 1} - 3\rho \quad
      \forall \ell \geq 0
   \]
and define $\Gamma ^-_n$ to be the rectangle with top and bottom side given by
$[\lambda ^- _{2n-1} - 2\rho, \lambda ^-_{2n} + 2\rho ] \pm i\rho
$. (Recall from \cite[Proposition B.11]{KP} that the periodic eigenvalues of Hill operators
are compact functions of the potential on $L^2(\mathbb T)$. It then follows from the Counting
Lemma that for any $k \ge 0,$  the k'th band length 
$\lambda^-_{2k+1} - \lambda^-_{2k}$ can be estimated from below
-- and hence $\rho > 0$ chosen -- uniformly on bounded subsets of functions $\alpha, \beta \in C^2_0(\mathbb T).$) 
Furthermore, let
\be\label{tau}
\tau^-_\ell=\frac{\lambda^-_{2\ell-1}+\lambda^-_{2\ell}}2
\ee
and denote by $V$ the set of all sequences $s = (s_\ell )_{\ell \geq 1}$ in
$\ell ^2$ such that for any $\ell \geq 1,$ the sequence $ \sigma _\ell =
\tau ^-_\ell + s_\ell $ satisfies $\lambda ^-_{2\ell - 1} - \rho \le
\sigma _\ell \le \lambda ^-_{2\ell } + \rho $. In particular the
sequence $\sigma_{\ell}$ is strictly increasing. Note that $V$ is a closed,
connected subset of $\ell ^2$. Furthermore, for any $k \geq 1$,
let $\frak F^k\equiv (\frak F^k_n)_{n \geq 1}$ denote the map 
from  $\ell ^2$ to $\ell ^2$ whose n'th component is given by
   \[ \frak F^k_n : s = (s_\ell )_{\ell \geq 1} \mapsto \begin{cases} \int _{\Gamma ^-_n}
      \left( \prod _{\ell \not= k} \frac{\sigma _\ell - \lambda }{w^-_\ell
      (\lambda)} \right) \cdot \frac{\pi ^2_k - \pi ^2_n}{w^-_k(\lambda )} \cdot
      \frac{2n\pi}{\sqrt[+]{\lambda - \lambda^-_0}} d \lambda &\mbox{if } n \not=
      k \\
      s_k &\mbox{if } n = k \end{cases}
   \]
In \cite[Appendix D]{KP}  it is shown that this map is real analytic and its
differential at any point $s \in V$ a linear
isomorphism (cf \cite[Lemma D.6]{KP}). Furthermore,
% the inverses $(d_s \frak F^k)^{-1}$ are uniformly bounded for $s \in V$ 
%(cf \cite[Proposition D.8]{KP}) and
$s^{-,k} = (s^{-,k}_\ell )
_{\ell \geq 1}$, given by
 \begin{equation}
\label{s^-solution}
 s^{-,k}_\ell = \begin{cases} \sigma ^{-,k}_\ell - \tau ^-_\ell=0 &\ell \not=
      k \\ 0 &\ell = k \end{cases}
\end{equation} 

is the unique solution of $\frak F^k(s) = 0$ in $V$. Let 
$V^k:=\bigcup _{N \geq 3} V^{k,N}$ where
$(V^{k,N})_{N \ge 3}$ is a sequence of subsets of $V$ with $V^{k,N}$ given by
   \[ \{ s \in V : s_\ell = \sigma ^{-,k}_\ell - \tau ^- _\ell  \ \forall \ell
      > L; \, \ - \frac{\gamma ^-_\ell }{2} - C \frac{M^2}{N} \leq
       s_\ell \leq \frac{\gamma ^-_\ell }{2} + C \frac{M^2}{N}
       \ \forall \ell \leq L \}
   \]
where $C > 0$ denotes the constant in the error estimate 
of Theorem \ref{Theorem 1.1} for the
asymptotics of $\lambda ^N_\ell , \ell \leq M$.  Next we want to show that for any $k \geq 1$, the
restriction of $\frak F^k$ to $V^k$ is $1 - 1$.  
To this end we first verify that $V^k$ is relatively compact in $\ell ^2$. To see it note that as $F(N) \leq N^\eta $ with $\eta \leq \frac{1}{3}$ one
has $\sum_{\ell \leq F(M)} (\frac{M^2}{N})^2 \le N^{\eta^2} \cdot N^{4\eta - 2} = O(N^{-5/9})$
and hence $V^k$ is bounded in $\ell ^2$. To see that the sequences in $V^k$ are uniformly
summable, choose a strictly increasing sequence
$(\ell_j)_{j \ge 1}$ in $\mathbb N$ so that 
$\sum_{\ell > \ell_j} (\gamma^-_\ell)^2 \le \frac{1}{j}$.
By the first part of assumption (F),  there then exists another strictly increasing sequence $(N_j)_{j \ge 1}$ 
in $\mathbb N$ so that $F[F(N_j)] \le \ell_j < F[F(N_j + 1)].$ Then, for any sequence 
$s=(s_\ell)_{\ell \ge 1} \in V^{k,N}$ with $N \le N_j$ one has
$\sum_{\ell > \ell_j} s_\ell^2 \le  \sum_{\ell >\ell_j} (\gamma^-_\ell)^2 \le \frac{1}{j} $
whereas for sequences in $V^{k,N}$ with $N > N_j$ one has
\[
\sum_{\ell > \ell_j} s_\ell^2\le  \sum_{\ell > \ell_j} (\gamma^-_\ell)^2  
+ 2C^2 \sum_{\ell \le F(M)} \frac{M^4}{N^2}
\le \frac{1}{j} + 2C^2 N_j^{\eta^2} \cdot N_j^{4\eta - 2} 
\le \frac{1}{j}  + 2C^2 N_j^{-5/9}.
\]
 Altogether, we thus have proved that $V^k$ is a relatively compact subset of $V$
and hence, as $V$ is closed,  $\overline{V^k} \subset V$ is
compact. Further note that $V^k$ and hence $\overline{V^k}$ is star shaped 
and thus connected. The following Lemma is the key ingredient of the proof
of Theorem~\ref{theorem2.1}.

\begin{lemma}
\label{1-1}
For any $k\geq 1$, the map  $\frak F^k
\big\arrowvert _{V^k} : V^k \rightarrow \frak F^k(V^k)$ is one-to-one and
together with its inverse uniformly Lipschitz.
\end{lemma}
\proof 
Let
$ \tilde V^k:= \{ s \in \overline{V^k} : \sharp \left( (\frak F^k)^{-1}(\frak F^k
      (s)) \cap \overline{V^k} \right) = 1 \}.$
By the observation above, $(s^{-,k}_\ell)_{\ell \ge 1}$ is in $\tilde V^k$ and hence  
$\tilde V^k \not= \emptyset $. As $\frak F^k$ is a local
diffeomorphism, $\tilde V^k$ is closed. Furthermore by the compactness of
$\overline{V^k}$ and the local diffeomorphism property of $\frak F^k$ one
concludes that the complement of $\tilde V^k$ in $V^k$, $\overline {V^k} \backslash \tilde V^k,$ is also closed. As
$\overline{V^k}$ is connected we thus have proved that $\overline{V^k} =
\tilde V^k$, i.e. $\frak F^k \big\arrowvert _{\overline {V^k}}$ is $1 - 1$. 
As $\overline{V^k}$ is compact one can argue as in the proof of
\cite[Proposition D.8]{KP} that $d_s \frak F^k$ is uniformly
boundedly invertible for any $s \in \overline{V^k}$ and any $k \geq 1$, implying that 
$\frak F^k \big\arrowvert _{\overline
{V^k}}$ and its inverse are uniformly Lipschitz for any $k \geq 1$. \qed

\medskip

For any $1 \le  k \le L$ introduce
   \[ \tilde\sigma ^{N,k}_\ell = \begin{cases} 4N^2(\sigma ^{N,k}_\ell + 2)
      &1 \leq \ell \leq L, \ell \not= k \\ \tau ^-_k &\ell = k \\ \sigma ^{-,
      k}_\ell &\ell > L \end{cases} ,
   \]
and the corresponding shifted sequence
 $  \tilde s^{N,k} = \big( \tilde s^{N,k}_\ell = \tilde \sigma ^{N,k}_\ell - \tau ^-_\ell
      \big) _{\ell \geq 1} \in \ell ^2 \ .
  $
\begin{lemma}
\label{inV}
There exists $N_0 \ge 3$ so that $\tilde s^{N,k}
\in V^{k,N}$ for any $1 \leq k \leq L$ and $N \geq N_0 $. 
\end{lemma}
\proof 
By Theorem \ref{Theorem 1.1} that there exist $N_0 \geq 3$ and $C>0$
so that for any $ N \geq N_0$ and $1 \leq  \ell \leq M$
\[ \lambda ^-_{2\ell - 1} - C \frac{M^2}{N} \leq \nu^N_{2\ell -1},\nu^N_{2\ell}\leq \lambda
      ^-_{2\ell } + C \frac{M^2}{N} \quad \mbox{and} \quad C \frac{M^2}{N} < \rho
   \]
where $\nu^N_j = 4N^2 (\lambda^{N,k}_j  +2).$
As $\lambda^N_{2 \ell -1} \le \sigma^{N,k}_\ell \le \lambda^N_{2 \ell}$ for $0 < k,\ell < N$ 
the definition of $\tilde \sigma^{N,k}_\ell$ implies that for any $1 \leq  k \leq L$,
$N \geq N_0$ one has for $1 \leq \ell \leq L$
\begin{equation}
\label{estimatesigmatilde}
 \lambda ^-_{2\ell - 1} - C \frac{M^2}{N} \leq \tilde \sigma ^{N,k}_\ell
     \leq \lambda ^-_{2\ell } + C \frac{M^2}{N}
 \end{equation}
or, by subtracting $\tau _\ell ^-$
 \[ - \gamma ^-_\ell / 2 - C\frac{M^2}{N} \leq \tilde s^{N,k}_\ell \leq
      \gamma ^-_\ell / 2 + C \frac{M^2}{N}
\]
and for any $\ell \geq 1, $
 $ \lambda ^-_{2\ell - 1} - \rho \leq \tilde \sigma^{N,k}_\ell \leq \lambda
      ^-_{2\ell } + \rho. 
 $
Altogether we thus have shown that $\tilde s^{N,k} \in V^{k,N}$
for any $ N \geq N_0$ and $1 \leq  \ell \leq L$. \qed

\begin{lemma}
\label{z.2}
The first $L$ elements of the sequence  $\frak F^k (\tilde s ^{N,k})$ satisfy
   \begin{equation}
   \label{Z.12} \big( \sum ^{L}_{n=1} \big\arrowvert \frak F^k_n(\tilde s^{N,k})
                \big\arrowvert ^2 \big) ^{1/2} = O \big( \frac{1}{L^{5/2}}  \big)
   \end{equation}
uniformly for $1 \leq k \leq L$ and on bounded sets of potentials
$\alpha , \beta $ in $C^2_0 ({\mathbb T})$. 
\end{lemma}
\proof Recall that the $N-1$ zeroes $(\sigma ^{N,k}_\ell)_{0<\ell <N} $ of
$\varphi ^N_k$ are determined by the $N-1$ equations $(\ 0 < n < N, n \not= k)$
   \begin{equation}
   \label{Z.1} 0 = \int _{\Gamma ^N_n} \big( \prod _{\underset{\ell \not= k}
               {0 < \ell < N}} \frac{ \sigma ^{N,k}_\ell - \mu }{w^N_\ell (\mu )}
               \big) \frac{\pi ^2_k - \pi ^2_n}{4N^2 w^N_k(\mu )} \frac{2n\pi }
               {2N\sqrt[+]{\mu - \lambda ^N_0}} \cdot \frac{2}{\sqrt[+]{\lambda ^N
               _{2N-1} - \mu }} d\mu .
   \end{equation}

Introduce the products
   \begin{align}
\label{P.4bis}
&Q^{N,L}_k(\mu )= \frac{2}{\sqrt[+]{\lambda ^N_{2N-1} - \mu }} \cdot
                               \prod _{L < \ell < N} \frac{ \sigma ^{N,k}_\ell -\mu}{w^N_\ell (\mu)}
\\
\label{P.4bis.1} &Q^{-,L}_k(\lambda )= \prod _{\ell > L} \frac{\sigma ^{-,k}
                     _\ell - \lambda }{w^-_\ell (\lambda )}
   \end{align}
\begin{equation}
\label{gnk0}
 \frak G^{N,L}_k(\lambda ) =  \frac{Q^{-,L}_k(\lambda )}
      {Q^{N,L}_k(\mu )} \cdot \frac{2N\sqrt[+]{\mu - \lambda ^N_0}}{\sqrt[+]{\lambda -
      \lambda ^-_0}} \cdot \frac{4N^2 w^N_k(\mu )}{w^-_k(\lambda )}
\end{equation}
where $\mu\equiv \mu(\lambda) = -2 + \frac{\lambda}{4N^2}$. Then
   \begin{align*} \frak F^k_n (\tilde s^{N,k}) &= \int _{\Gamma ^-_n} \big(
                     \prod _{\underset{\ell \not= k}{1 \leq \ell \leq L}} \frac
                     {\tilde \sigma ^{N,k}_\ell - \lambda }{4N^2w^N_\ell (\mu ) }
                     \big) \cdot Q^{N,L}_k(\mu ) \frac{\pi ^2_k - \pi ^2_n}
                     {4N^2 w^N_k(\mu )} \frac{2n\pi }{2N \sqrt[+]{\mu -
                         \lambda ^{N}_0}} 
\\
&\cdot
\frac{2}{\sqrt[+]{\lambda ^N_{2N-1} - \mu }}\frak G^{N,L}_k(\lambda ) 
\prod _{\underset{\ell \not= k}{1 \leq \ell \leq L}} \frac{4N^2w^N_\ell (\mu )}{w^-_\ell (\lambda )} d\lambda.
   \end{align*}
The terms in the latter integral are estimated separately. 
By Theorem \ref{Theorem 1.1} (cf also Proposition \ref{basicEstimate}), one has for any $\lambda \in \Gamma ^-_n$,
 \begin{equation}
   \label{Z.10} \prod _{\underset{\ell \not= k}{1 \leq \ell \leq L}} \frac{4N^2
                   w^N_\ell (\mu )}{w^-_\ell (\lambda )} = \frac{4N^2w^N_n(\mu )}
                   {w^-_n(\lambda )} \cdot
                   \prod _{\underset{\ell \not= k,n}{1 \leq \ell \leq L}}
                   \frac{4N^2 w^N_\ell (\mu )}{w^-_\ell (\lambda )} 
\end{equation}
\[= \big( 1 + O \big( \frac{M^2}{N} \big) \big) \mbox{exp}
                   \Big( \sum _{\underset{\ell \not= k,n}{1 \leq \ell \leq L}}
                   \log \big( 1 + O \big( \frac{M^2}{N} \frac{1}{|\ell ^2 -
                   n^2|} \big) \big) \Big) \nonumber 
                  = 1 + O \big( \frac{M^2} {N} \big) .
\]
By Lemma \ref{G} one then has
$\frak G^{N,L}_k(\lambda ) \cdot \prod _{\underset{\ell \not= k}{1 \leq \ell \leq L}}
 \frac{4N^2w^N_\ell (\mu )}{w^-_\ell (\lambda )}  = 1 + O \big(\frac{1}{L^3} \big)$
 uniformly for $1 \leq n, k \leq L, n \not= k$ and $\lambda \in \Gamma ^-_n$.
By \eqref{Z.1} it then follows that
   \begin{align*} &\frak F ^k_N(\tilde s^{N,k}) = \int _{\Gamma ^-_n} \big(
                     \prod _{\underset{\ell \not= k}{1 \leq \ell \leq L}} \frac
                     { \tilde \sigma ^{N,k}_\ell - \lambda}{4N^2 w^N_\ell (\mu )}
                     \big) \cdot Q^{N,L}_k(\mu ) \\
                  &\cdot \frac{\pi ^2_k - \pi ^2_n}{4N^2 w^N_k(\mu )} \cdot \frac{
                     2n\pi }{2N \sqrt[+]{\mu - \lambda ^N_0}} \cdot \frac{2}{\sqrt[+]{
                     \lambda ^N_{2N-1} - \mu }} \cdot O \left(\frac{1}{L^3} \right) d\lambda .
   \end{align*}
Furthermore, again by Theorem \ref{Theorem 1.1} one has  uniformly for $1 \leq n, k \leq L$ and $\lambda \in \Gamma ^-_n$
\[
\frac{2n\pi }{2N \sqrt[+]{\mu - \lambda ^N_0}}, \quad
  \frac{2}{\sqrt[+]{\lambda ^N_{2N-1} - \mu }},\quad
\frac{ \tilde \sigma ^{N,k}_n - \lambda }{4N^2 w^N_n (\mu )} \,=\, O(1) .
\]
It remains to estimate the term $\prod _{\underset
{\ell \not= k}{1 \leq \ell \leq L}} \frac{ \tilde \sigma ^{N,k}_\ell - \lambda}
{4N^2 w^N_\ell (\mu )}$. 
By Theorem \ref{Theorem 1.1} (cf also Proposition \ref{basicEstimate}),
 $\tilde \sigma ^{N,k}_\ell - \nu ^N_j = O \left( \gamma_\ell^- +\frac{M^2}{N}  \right)$
for $j \in \{ 2\ell , 2\ell - 1 \}$ and $1 \le \ell \le L$
where $\nu^N_j = 4N^2 (\lambda^{N,k}_j  +2)$.
Hence, for 
$\lambda \in \Gamma ^-_n$ and $1 \leq \ell \leq L$ with $\ell \not= k, n$
   \[ \frac{ \tilde \sigma ^{N,k}_\ell - \lambda }{4N^2 w^N_\ell (\mu ) } 
     = \Big( \frac{(\tilde \sigma^{N,k}_\ell - \lambda )(\tilde \sigma ^{N,k}_\ell -
      \lambda )}{(\nu ^N_{2\ell } - \lambda )(\nu ^N_{2\ell - 1} - \lambda )} \Big) ^{1/2}
      = 1 + O \Big( ( \gamma_\ell^- +\frac{M^2}{N} ) \frac{1}{|\ell ^2 - n^2|} \Big)
   \]
and thus
\[
 \prod _{\underset{\ell \not= k,n}{1 \leq \ell \leq L}} 
\frac{\lambda - \tilde \sigma ^{N,k}_\ell }{4N^2 w^N_\ell (\mu )} =
                     \exp \big( \sum _{\underset{\ell \not= k,n}{1 \leq \ell \leq
                     L}} \log \big( 1 + O \big( (( \gamma_\ell^- +\frac{M^2}{N} )) \frac{1}
                     {|\ell ^2 - n ^2|} \big) \big) \big) = O \big( 1 \big).
\]
In view of \eqref{Fwith1/3} we thus proved that for any $1 \leq n \leq L, n \not= k$,
$\frak F^k_n(\tilde s^{N,k}) = O \Big( \frac{1}{L^3} \Big) $. 
Going through the arguments of the proof one verifies the claimed uniformity
statement.\qed

\medskip

\begin{lemma}
\label{z.2.1}
The tail $(\frak F^k_n (\tilde s ^{N,k}))_{n > L}$ of $\frak F^k (\tilde s ^{N,k})$ satisfies the estimate
   \begin{equation}
   \label{Z.15} \big( \sum _{n > L} |\frak F^k_n (\tilde s ^{N,k})|^2 \big) ^{1/2} =
                O \big(\frac{1}{L^4} \big)
   \end{equation}
uniformly in $1 \leq k \leq L$ and on bounded sets of potentials
$\alpha , \beta $ in $C^2_0 ({\mathbb T})$. 
\end{lemma}
\proof By the definition of $\tilde s^{N,k}$ and $\frak F^k_n$ one has
   \[ \frak F^k_n(\tilde s^{N,k}) = \int _{\Gamma ^-_n} \Big( \prod _{\underset{\ell
      \not= k}{1 \leq \ell \leq L}}\frac{ \tilde \sigma ^{N,k}_\ell - \lambda }
      {\sigma ^{-,k}_\ell - \lambda } \Big) \zeta ^{-,k}_n (\lambda ) \frac
      {\sigma ^{-,k}_n - \lambda }{w^-_n(\lambda )}d\lambda
   \]
where
   \[ \zeta ^{-,k}_n(\lambda ) := \Big( \prod _{\underset{\ell
      \not= k,n}{\ell \geq 1}}  \frac
      {\sigma ^{-,k}_\ell - \lambda }{w^-_n(\lambda )} \Big) \frac{\pi ^2_k
      - \pi ^2_n}{w^-_k(\lambda )} \frac{2n\pi }{\sqrt{\lambda - \lambda ^-_0}} .
   \]
We argue as in \eqref{Z.10} and write for $\lambda \in \Gamma^-_n$ with $n > L$
\[
 \prod _{\underset{\ell \not= k}{1 \leq \ell \leq L}}
                     \frac{ \tilde \sigma ^{N,k}_\ell - \lambda }{\sigma ^{-,k}_ \ell - \lambda } 
     = \exp \big( \sum _{\underset{\ell \not= k}{1 \leq \ell \leq L}}
                     \log \big( 1 + \frac{\tilde \sigma ^{N,k} _\ell -\sigma ^{-,k}_\ell}
                    {\sigma ^{-,k}_\ell - \lambda} \big) \big) .
 \]
Note that by \eqref{estimatesigmatilde},
$\big\arrowvert \sigma ^{N,k}_\ell - \tau ^-_\ell \big\arrowvert \leq
      \gamma ^-_\ell/2 + C \frac{M^2}{N}$.
Taking into account that $| \sigma ^{-,k}_\ell - \tau ^-_\ell | \leq \gamma ^- _\ell / 2$ it then follows that
   \[ \big\arrowvert \tilde \sigma ^{N,k} _\ell -\sigma ^{-,k}_\ell
      \big\arrowvert \leq \gamma ^-_\ell + C \frac{M^2}{N}
      = O \left( \frac{1}{\ell ^2} + \frac{M^2}{N} \right) .
   \]
As $(\sigma ^{-,k}_\ell - \lambda)^{-1} = O \left( \frac{1}{L^2} \right) $
for $\ell \leq L/2$ and $n > L,$ one concludes from \eqref{Fwith1/3}
   \[ \sum _{\underset{\ell \not= k}{1 \leq \ell \leq L/2}} \big\arrowvert 
    \frac{\tilde \sigma ^{N,k} _\ell -\sigma ^{-,k}_\ell }{\sigma ^{-,k}_\ell - \lambda}
      \big\arrowvert = O \big( \frac{1}{L^2} (\sum_{\ell \ge 1} \frac{1}{\ell^2} + \frac{M^2}{N} L) \big)
      = O(\frac{1}{L^2}).
   \]
Furthermore
   \[ \sum _{\underset{\ell \not= k}{L/2 < \ell \leq L}}
      \big\arrowvert \frac{ \tilde \sigma ^{N,k} _\ell -\sigma ^{-,k}_\ell }{\sigma ^{-,k}_\ell - \lambda}
      \big\arrowvert = O \big( \frac{\log L}{L^3} \big).
   \]
Altogether we thus have proved that for $n > L$
   \[ \prod _{\underset{\ell \not= k}{1 \leq \ell \leq L}} \frac{ \tilde \sigma ^{N,k}
      - \lambda }{\sigma ^{-,k}_ \ell - \lambda } - 1 = O \big( \frac{1}{L^2}\big) .
   \]
Actually, these asymptotic estimates hold not only for $\lambda \in \Gamma ^-_n$, but also for
any $\lambda $ in the interior of $\Gamma ^-_n$, uniformly in $1 \leq k \leq L$
and $n > L$. Furthermore, one verifies in a straightforward way that 
$\zeta ^{-,k}_n (\lambda ) = O(1)$ for $\lambda \in \Gamma ^-_n$ 
and in the interior of $\Gamma ^-_n.$ 
To prove the claimed estimate note that $\int _{\Gamma ^-_n} \zeta ^{-,k}_n(\lambda ) 
\frac{\sigma ^{-,k}_n - \lambda }{w^-_n(\lambda )} d\lambda = 0$
 by the definition of $\big( \sigma ^{-,k}_\ell \big)_{\ell \ne k},$ implying that
   \[\frak F^k_n(\tilde s^{N,k}) = \int _{\Gamma ^-_n} \zeta ^{-,k}_n(\lambda ) 
     \big( \prod _{\underset{\ell \not= k}{1 \leq \ell \leq L}} \frac{ \tilde \sigma ^{N,k}
      - \lambda }{\sigma ^{-,k}_ \ell - \lambda } - 1)
     \frac{\sigma ^{-,k}_n - \lambda }{w^-_n(\lambda )} d\lambda .
   \]
In the case $\gamma^-_n = 0, $ one has $\frac{\sigma ^{-,k}_n - \lambda }{w^-_n(\lambda )} =1$
and as the integrand is analytic in $\lambda$ the latter integral vanishes by Cauchy's theorem.
In the case $\gamma^-_n > 0$ we deform the contour $\Gamma ^-_n$ to the interval 
$[\lambda ^-_{2n-1}, \lambda ^-_{2n}]$ to conclude by the mean value theorem that
   \[ \frak F^k_n(\tilde s^{N,k}) = O \big(\frac{1}{L^2} \gamma ^-_n \big) 
= O \big(  \frac{1}{L^4}  n^2\gamma^-_n \big).
   \]
The claimed estimate then follows from the assumption that 
$\alpha, \beta \in C_0^2(\mathbb T).$
Going through the arguments of the proof one verifies the claimed uniformity
statement. \qed

\bigskip

{\it Proof of Theorem \ref{theorem2.1}.}  Combining
\eqref{Z.12} and \eqref{Z.15} yields
   \[ \|\frak F^k(\tilde s^{N,k}) \| _{\ell ^2} = O \big(\frac{1}{L^{5/2}} \big)
   \]
uniformly for $1 \le k \le L$. 
As for any $1 \le k \le L,$ $\frak F^k(\tilde s^{-,k}) = 0$ (cf \eqref{s^-solution}), $\tilde s^{N,k} \in V^k$ (Lemma \ref{inV}),
$\frak F^k$ is $1-1$ and together with its inverse uniformly Lipschitz on $V^k$ (Lemma \ref{1-1})
 we then get the claimed estimate
   \[ \| \tilde \sigma ^{N,k} - \sigma ^{-,k} \| _{\ell ^2} = 
        \| \tilde s ^{N,k} - s ^{-,k} \| _{\ell ^2}
      = O \big(\frac{1}{L^{5/2} } \big) .
 %     \eqno \square
   \]
Going through the arguments of the proof one verifies the claimed uniformity
statement. 
\hspace*{\fill }$\square $

\medskip

%%%%%%%%%%%%%%%%%%%%%%%%%%%%%%%%%%%%%%%%%%%%%%%%%%%%%%%%%%%%%%%%%%%%%%%%%
%%%%%%%%%%%%%%%%%%%%%%%%%%%%%%%%%%%%%%%%%%%%%%%%%%%%%%%%%%%%%%%%%%%%%%%%%%
%%%%%%%%%%%%%%%%%%%%%%%%%%%%%%%%%%%%%%%%%%%%%%%%%%%%%%%%%%%%%%%%%%%%%%%%%%

\section{Leading order asymptotics}\label{priori}

In this section we compute the leading order asymptotics of the Toda frequencies $\omega^N_n$.
In particular we prove the asymptotics \eqref{freq3inSec6} in the bulk.
In addition we identify the principal contributions in the formulas \eqref{B.5} 
 to the asymptotics of $\omega ^N_n$ at the two edges.

\begin{remark}
\label{actionedge}
We recall that for the remainder of the paper we set $M=[F(N)]$ and $L=[F(M)]$ 
where $F : {\mathbb N} \rightarrow {\mathbb R}_{\geq 1}$
satisfies $(F)$ with $\eta \le 1/3.$ Then \eqref{Fwith1/3} holds and
the asymptotics of Theorem \ref{thm1.5} and Theorem \ref{Theorem 9.1} at the right and left edges read as follows:
for any $1 \le n \le L$ 
\[8N^2 I^N_n = I^-_n + O((\frac{M^2}{N} + \gamma_n^-)\frac{L}{M^{1/2}}) \quad
8N^2 I^N_{N-n} = I^+_n + O((\frac{M^2}{N} + \gamma_n^+)\frac{L}{M^{1/2}})
\]
\[
J^N_n - J^-_n,\,\,\, J^N_{N-n} - J^+_n = O(\frac{L}{M^{1/2}})
\]
\end{remark}

\begin{proposition}
\label{Proposition B.1} Uniformly for $M < n < N - M$ and on bounded
subsets of functions $\alpha , \beta $ in $C_0^2(\mathbb T)$
   \[ \omega ^N_n = 2 \big( \sin \frac{n\pi }{N}\big) 
     \big( 1 + O \big( \frac{\log M}{M^2} \big) \big) .
   \] 
\end{proposition}

{\it Proof:} By \eqref{B.3} there exists $\mu _\ast \in [\lambda ^N_{2n-1}, \lambda ^N_{2n}]$ so that
   \[ \omega ^N_n = \sqrt[+]{(\lambda ^N_{2N-1} - \mu _\ast )(\mu _\ast - \lambda ^N
      _0)} \prod _{\underset{1 \leq k < N}{k \not= n}} \frac{\sqrt[+]{(\lambda^N
      _{2k}-\mu _\ast )(\lambda ^N_{2k-1} - \mu _\ast )}}{|\sigma ^{N,n}_k -
      \mu _\ast|}
   \]
By Theorem~\ref{Theorem 1.1}, 
$\mu _\ast = -2 \cos \frac{n\pi }{N} + O \big( \frac{1}{N^2 M} \big)$, 
$\lambda ^N_0 = - 2 + O(N^{-2})$, and $\lambda ^N_{2N-1} = 2 +O(N^{-2})$. 
Hence
   \[ (\lambda ^N_{2N-1} - \mu _\ast)(\mu _\ast - \lambda ^N_0) 
     = 4 \sin ^2 \frac{n \pi }{N} + O \big( \frac{1}{N^2} \big) .
   \]
As $\sin \frac{n\pi}{N} \geq \sin \frac{M\pi }{N} = O \big( \frac{M}{N} \big)$
 for $M < n < N - M$,
   \begin{equation}
   \label{B.61} \sqrt[+]{(\lambda ^N_{2N-1} - \mu _\ast )(\mu _\ast - \lambda ^N_0)}
               = (2\sin \frac{n\pi }{N}) \big( 1 + O \big( \frac{1}{M^2} \big) \big) .
   \end{equation}
Next we write
   \[ \prod _{\underset{0 < k < N}{k \not= n}} \frac{\sqrt{(\lambda ^N_{2k} -
      \mu _\ast )(\lambda ^N_{2k-1} - \mu _\ast)}}{|\sigma ^{N,n}_k - \mu _\ast|} =
      \Big( \prod _{\underset{0 < k < N}{k \not= n}} \frac{|\lambda ^N_{2k}
      - \mu _\ast|}{|\sigma ^{N,n}_k - \mu _\ast |}
       \prod _{\underset{0 < k < N}{k \not= n}} \frac{|\lambda ^N_{2k-1}
      - \mu _\ast|}{|\sigma ^{N,n}_k - \mu_\ast|} \Big) ^{1/2}.
   \]
The two products are estimated in the same way, so we concentrate on
the first one.  Note that $0 < \frac{\lambda ^N_{2k} - \mu _\ast }
        {\sigma ^{N,n}_k - \mu _\ast} 
= 1 + \frac{\lambda ^N_{2k} - \sigma ^{N,n}_k}{\sigma ^{N,n}_k - \mu _\ast}$ 
for any $k \ne n.$ It will turn out that 
$\frac{\lambda ^N_{2k} - \sigma ^{N,n}_k}{\sigma ^{N,n}_k - \mu _\ast}$ 
is very small so that we can write
   \begin{equation}
   \label{B.71} 
   \prod _{\underset{0 < k < N}{k \not= n}} \frac{\lambda ^N_{2k} - \mu _\ast }
        {\sigma ^{N,n}_k - \mu _\ast} 
   = \exp \Big( \sum _{\underset{0 < k < N}{k \not= n}} 
      \log \big( 1 + \frac{\lambda ^N_{2k} - \sigma ^{N,n}_k}{\sigma ^{N,n}_k - \mu _\ast}\big) \Big)
   \end{equation}
and get the estimate 
 \[
\exp \Big( \sum _{\underset{0 < k < N}{k \not= n}} 
      \log \big( 1 + \frac{\lambda ^N_{2k} - \sigma ^{N,n}_k}{\sigma ^{N,n}_k - \mu _\ast}\big) \Big)
-1 = O( \sum _{\underset{0 < k < N}{k \not= n}} \frac{|\lambda ^N_{2k} - \sigma
      ^{N,n}_k|}{|\sigma ^{N,n}_k - \mu _\ast |}) .
 \]
The latter sum is split up and as 
$\lambda ^{N}_{2k-1} \le \sigma ^{N,n}_k \le \lambda ^{N}_{2k}$ can be estimated as follows
\[
\sum _{\underset{0 < k < N}{k \not= n}} 
\frac{|\lambda ^N_{2k} - \sigma^{N,n}_k|}{|\sigma ^{N,n}_k - \mu _\ast |}
\le  \sum _{0 < k < n} \frac{\gamma^N_k} {\lambda^N_{2n-1} - \lambda^N_{2k}}
+ \sum _{n < k < N} \frac{\gamma^N_k} {\lambda^N_{2k-1} - \lambda^N_{2n}}.
\]
Taking into account that $\eta \le 1/3$ and hence $\frac{M^3}{N} = o(1),$ the right hand
side of the latter inequality is $O(\frac{\log M}{M^2})$ by Proposition \ref{basicEstimate} (ii). 
Altogether we thus have proved that
\[ \prod _{\underset{0 < k < N}{k \not= n}} \frac{\sqrt{(\lambda ^N_{2k} -
      \mu _\ast )(\lambda ^N_{2k-1} - \mu _\ast)}}{|\sigma ^{N,n}_k - \mu _\ast|} =
1 +O(\frac{\log M}{M^2}).
\]
Combined with \eqref{B.61}, we then obtain the stated estimate.
Going through the arguments of the proof, one verifies the claimed uniformity statement.
\hspace*{\fill }$\square $	

\begin{proposition}
\label{Proposition B.2} Uniformly for $1 \le n \le M$ and on bounded
subsets of functions $\alpha , \beta $ in $C_0^2(\mathbb T)$ one has
$ {\mbox{(i)}}\,\,\, \omega ^N_n = O \left( \frac{n}{N} \right)$  and $ {\mbox{(ii)}}
\,\, \, \omega ^N_{N-n} = O \left( \frac{n}{N} \right).$
\end{proposition}

{\it Proof:} The estimates (i) and (ii) are proved in the same way and
so we concentrate on (i) only. We use again formula \eqref{B.3}.  As by
Theorem~\ref{Theorem 1.1}
   \[\lambda ^N_{2N-1} - \mu _\ast \le 4 - \frac{\lambda ^+_0 + \lambda ^-
      _{2n-1}}{4N^2} + O \left( \frac{M^2}{N^3} \right), \,\,\,
\mu _\ast - \lambda ^N_0 \leq 
      \frac{\lambda ^-_{2n} - \lambda ^-_0}{4N^2} + O \left( \frac{M^2}{N^3} \right) 
\]
one gets
   \begin{equation}
   \label{B.15}
      \sqrt[+]{(\lambda ^N_{2N-1} - \mu _\ast )(\mu _\ast - \lambda ^N_0)} = O
      \big( \frac{\sqrt{\lambda ^-_{2n} - \lambda ^-_0}}{N} \big) = O
      \big( \frac{n}{N} \big) .
   \end{equation}
Next we need to estimate the product
 \[
 \prod _{k\not= n} \frac{\sqrt[+]{(\lambda ^N_{2k} - \mu _\ast )(\lambda ^N_{2k-1}
         - \mu _\ast)}}{|\sigma ^{N,n}_k - \mu _\ast |} 
     = \Big( \prod
         _{\underset{1 \leq k < N}{k\not= n}} \frac{|\lambda ^N_{2k} - \mu _\ast |}{|
         \sigma ^{N,n}_k - \mu _\ast |} \Big) ^{1/2}  \Big( \prod
         _{\underset{1 \leq k < N}{k\not= n}} \frac{|\lambda ^N_{2k-1} - \mu _\ast |}{|
         \sigma ^{N,n}_k - \mu _\ast |} \Big) ^{1/2} .
 \]
The two products in the latter expression are estimated in the same way and thus 
we concentrate again on the first one. Argue as in the proof of Proposition \ref{Proposition B.1}
and use that by Proposition \ref{basicEstimate}(i)
\[
\sum _{\underset{0 < k < N}{k \not= n}} 
\frac{|\lambda ^N_{2k} - \sigma^{N,n}_k|}{|\sigma ^{N,n}_k - \mu _\ast |}
\le  \sum _{0 < k < n} \frac{\gamma^N_k} {\lambda^N_{2n-1} - \lambda^N_{2k}}
+ \sum _{n < k < N} \frac{\gamma^N_k} {\lambda^N_{2k-1} - \lambda^N_{2n}} = O(1)
\]
to conclude that $\big( \prod
         _{\underset{1 \leq k < N}{k\not= n}} \frac{|\lambda ^N_{2k} - \mu _\ast |}{|
         \sigma ^{N,n}_k - \mu _\ast |} \big) ^{1/2} = O(1).$
Substituting the obtained estimates into formula \eqref{B.3} yields (i). Going through
 the arguments of the proof one verifies the claimed uniformity.
\hspace*{\fill }$\square $	

\medskip

%%%%%%%%%%%%%%%%%%%%%%%%%%%%%%%%%%%%%%%%%%%%%%%%%%%%%%%%%%%%%%%%%%%%%%%%%

%\section{Principal contributions}
%\label{first}
We now use the rough estimates of  the Toda frequencies $\omega ^N_n$ obtained above to identify the principal contributions in the formulas \eqref{B.5}, 
\[
iN\omega ^N_n = \sum ^n_{j=1} \int ^{\lambda ^N_{2j-1}}
                     _{\lambda ^N_{2j-2}} \frac{(\mu - \frac{1}{N}\frak p_N)\dot \Delta _N
                     (\mu )d\mu}{\sqrt[c]{\Delta _N(\mu )^2 - 4}} 
                  - \sum _{k \in \mathcal J_N} I^N_k \omega ^N_k \sum ^n_{j=1}
                     \int ^{\lambda ^N_{2j-1}}_{\lambda ^N_{2j-2}}
                     \frac{\varphi ^N_k(\mu )d\mu}{\sqrt[c]{\chi_N(\mu )}},
\]
to their asymptotics at the two edges.
Recall that 
%$\mathcal J_N:= \{ 0 < k < N | I^N_k \not= 0 \}$, 
$\frak p_N= \sum_{i=1}^{N}b^N_i$ and $\chi_N$ is the characteristic polynomial, defined by \eqref{rn}.

\medskip

\begin{proposition}
\label{Da.Proposition.1} 
Uniformly for $1 \le n \le M$, $n \le M_0 \le M,$ and on bounded
subsets of functions $\alpha , \beta $ in $C_0^2(\mathbb T)$
\begin{equation}
\label{da.1.1}
\omega_n^N=\frac{2\pi n}{N} -\frac{\pi}{4N^3}\sum_{j=1}^{n}\lambda^-_{2j-2} 
- (P1) - (P2) - (P3) + Error \quad
\end{equation}
where
\begin{equation}
\label{da.1.2}
(P1) := \frac{1}{4N^3} \sum ^n_{j=1} \int ^{\nu ^N_{2j-1}}_{\nu
      ^N_{2j-2}} \arccos \big( \frac{(-1)^{N+j}}{2} \Delta _N \big( - 2 +
      \frac{\lambda }{4N^2} \big) \big) d\lambda \qquad \qquad 
\nonumber
\end{equation}
\[
(P2) := \frac{1}{N} \sum _{\underset{I^N_k \ne 0}{1 \leq k \leq M_0}} {I^N_k}\omega^N_k  \int _{\I_n^N}
           \frac{\varphi ^N_k(\mu )}{i \sqrt[c]{\chi_N(\mu )}} d\mu,
\quad
(P3) := \frac{1}{N} \sum _{M_0 < k \leq M} {I^N_k}\omega^N_k  \int _{\I_n^N}
           \frac{\varphi ^N_k(\mu )}{i \sqrt[c]{\chi_N(\mu )}} d\mu 
\]
with $\nu^N_j = 4N^2(2 + \lambda ^N_j)$, $\I^N_n:=[\lambda_0^N,\lambda^N_{2n-1}]\backslash \big(\cup_{j=1}^{n-1} [\lambda_{2j-1}^N,\lambda_{2j}^N] \big)$, and
\begin{equation}
\label{remain}
Error = 
O\Big(\frac{1}{N^{3}} \big(\frac{nM^2}{N} +\frac n{M^3}\big) \Big)
\end{equation}
Furthermore, 
\begin{equation}
\label{da.1.4}
(P3) = O\big(\frac{\log(n+1)}{N^3} (\frac{1}{M_0^5}+\frac{M^4\log M}{N^2}) \big).
\end{equation}
\end{proposition}
\begin{remark}
The integer $M_0$ in the statement of Proposition \ref{Da.Proposition.1} is a free parameter
for which different choices will be made in the course of the proof of Theorem \ref{thm1.5inSec6}:
in the proof of Proposition \ref{Proposition B.3A}, Proposition \ref{Da.Proposition.1} will be applied
for $1 \le n \le M$ with $M_0 = M$ whereas in the proof of Theorem \ref{theo.fre}
we apply it for $1 \le n \le L$ with $M_0 = L.$
\end{remark}
\proof For any $1 \leq n \leq M,$ write the formula \eqref{B.5} for the nth frequency, recalled
above, as a sum
 $\omega_n^N  = (T1) + (T2) - (T3) - (T4)$ where
\begin{equation}
\label{o.1}
(T1) := - \frac{2}{N} \sum ^n_{j=1} \int ^{\lambda ^N_{2j-1}}_{\lambda ^N_{2j-2}} 
\frac{\dot \Delta _N(\mu )d\mu }{i \sqrt[c]{\Delta ^2_N(\mu ) - 4}}  
- \frac{\frak p_N}{N^2} \sum ^n_{j=1} \int ^{\lambda ^N_{2j-1}}_{\lambda ^N_{2j-2}} 
\frac{\dot \Delta _N(\mu )d\mu }{i \sqrt[c]{\Delta ^2_N(\mu ) - 4}} 
\end{equation}
\begin{equation}
\label{o.3}
(T2) :=  \frac{1}{N} \sum ^n_{j=1} \int ^{\lambda ^N_{2j-1}}
                     _{\lambda ^N_{2j-2}} \frac{(\mu + 2)\dot \Delta _N(\mu )d\mu }
                     {i \sqrt[c]{\Delta ^2_N(\mu ) - 4}}  \qquad
\end{equation}
\begin{equation}
\label{o.6}
 (T3) :=   \frac{1}{N} \sum _{\underset{I^N_k \ne 0}{1 \leq k \leq M_0}} I^N_k \omega^N_k
 \sum^n_{j=1}\int ^{\lambda ^N_{2j-1}}_{\lambda ^N_{2j-2}} 
\frac{\varphi ^N_k(\mu )d\mu } {i \sqrt[c]{\chi_N(\mu )}} .
\end{equation}
\begin{equation}
\label{o.5}
 (T4) :=  \frac{1}{N} \sum _{\underset{I^N_k \ne 0}{M_0 < k < N}} I^N_k \omega ^N_k\sum^n_{j=1}
\int ^{\lambda ^N_{2j-1}} _{\lambda ^N_{2j-2}} 
\frac{\varphi ^N_k(\mu )d\mu }{i \sqrt[c]{\chi_N(\mu )}} 
\end{equation}
The various terms are treated individually. We begin with the term $(T1)$.
Notice that 
$|\Delta _N(\mu )| < 2 \,\, \forall \mu \in (\lambda ^N_{2j-2}, \lambda ^N_{2j-1})\,$ and
 \[\Delta _N(\lambda^N _{2j-1}) = 2(-1)^{N+j}, \quad \Delta_N(\lambda^N _{2j-2}) = 2(-1)^{N+j+1} 
\]
and denote by $\arccos (x)$  the principal branch of the inverse of cosine, defined on
$[-1,1]$ with $\arccos (-1) = \pi $ and $\arccos (1) = 0$. Then for $s \in \{ \pm 1 \}$
   \[ \frac{d}{d\mu } \arccos \left( \frac{(-1)^s \Delta _N(\mu )}{2} \right) =
      - \frac{\dot \Delta _N(\mu )(-1)^s}{\sqrt[+]{4 - \Delta _N(\mu )^2}}
	\quad \forall \,\, \mu \in (\lambda _{2j-2}, \lambda _{2j-1}).
   \]
As by the definition of the $c$-root,
$ \sqrt[+]{4 - \Delta _N(\mu )^2} = (-1) ^{N + 1 + j} i \sqrt[c]{\Delta _N (\mu )^2 - 4}
$
one concludes that
$ \frac{\dot \Delta _N(\mu )}{i \sqrt[c]{\Delta _N(\mu )^2 - 4}} = \frac{d}
      {d \mu } \arccos  \left( (-1)^{N +j} \frac{\Delta _N(\mu )}{2} \right)
$
and hence
   \[ \int ^{\lambda ^N_{2j-1}} _{\lambda ^N_{2j-2}} \frac{\dot \Delta _N(\mu )}{i
      \sqrt[c]{\Delta ^2_N(\mu ) - 4}} d\mu = - \pi .
   \]
Finally, as $\frak p_N = O(N^{-3})$ by \cite[Proposition 8.1]{bkp2},
   \begin{equation}
     (T1) = \frac{2\pi n}{N}+O \big( \frac{n}{N^5} \big).
   \end{equation}
To estimate the term $(T2)$ integrate by parts to get for $(T2)$
\[
\sum ^n_{j=1} \frac{\mu + 2}{N} \arccos ((-1)^{N+j}
 \frac {\Delta _N(\mu )}{2}) \Big\arrowvert ^{\lambda
   ^N_{2j-1}}_{\lambda ^N_{2j-2}} -\frac{1}{N} \sum ^n_{j=1} \int
 ^{\lambda ^N_{2j-1}}_{\lambda ^N_{2j-2}} \arccos ((-1)^{N+j} \frac{\Delta _N(\mu )}{2}) d\mu
\]
\[
\quad  =- \frac {\pi }{4N^3} 
\sum^n_{j=1}\nu ^N_{2j-2} - \frac{1}{4N^3} \sum ^n_{j=1} \int ^{\nu ^N_{2j-1}}_{\nu ^N_{2j-2}} 
\arccos \big(\frac{(-1)^{N+j}}{2} \Delta _N \big( - 2 + \frac{\lambda}{4N^2} \big) \big) d\lambda
\] 
where we made the change of variable of integration  $\mu = - 2 + \lambda / 4N^2$. 
As by Theorem \ref{Theorem 1.1} $\nu ^N_j =\lambda^-_j+O(M^2/N)$
\begin{equation}
\label{esti2}
 (T2) = - \frac {\pi }{4N^3} \sum^n_{j=1}\lambda ^N_{2j-2} + O\big(\frac{nM^2}{N^4}\big) - (P1).
\end{equation}
Next observe that $(T3) = (P2)$. Concerning $(T4)$ note that 
the summation indices $k, j$ in the definition of $(T4)$ satisfy $n < k$ and $j \le n$,
implying that the improper Riemann integral $\int ^{\lambda ^N_{2j-1}} _{\lambda ^N_{2j-2}} 
\frac{\varphi ^N_k(\mu )d\mu }{i \sqrt[c]{\chi_N(\mu )}} $ exists even if $I_k^N =0$. Hence
splitting up $(T4)$ we get $(T4) =  (P3) + (E1) + (E2)+ (E3)$ where
\[
(E1) := \sum _{M < k \le \frac{N}{2}} \dots \quad
(E2) := \sum _{\frac{N}{2} < k < N - M} \dots  \quad
(E3) := \sum _{N-M \le k < N} \dots
\]
In fact, $\int ^{\lambda ^N_{2j-1}} _{\lambda ^N_{2j-2}} 
\frac{\varphi ^N_k(\mu )} {\sqrt[c]{\chi_N(\mu )}}d\mu =
O\big( \frac{1}{\tau^N_k - \tau^N_j} \frac{1}{N}\big)$
for such $k,j$ by Lemma \ref{integrale} yielding
\[
 (T4) = O \Big(  \sum _{n < k < N} I^N_k \omega ^N_k
 \sum ^n_{j = 1}\frac{1}{\tau^N_k - \tau^N_j}\frac{1}{N^2} \Big) .
\]
Concerning the estimate \eqref{da.1.4}, note that for $1 \le j \le n$ and $n< k \leq M$,
$\frac{1}{\tau^N_k - \tau^N_j }= O \big( \frac{N^2}{k^2 - j^2} \big)$ (cf Theorem \ref{Theorem 1.1}),
$I^N_k = O\big(\frac{1}{kN^2} ((\gamma^-_k)^2 + \frac{M^4}{N^2}) \big)$ 
(Proposition \ref{Proposition C.1}), 
and $\omega^N_k = O(\frac{k}{N})$ (Proposition \ref{Proposition B.2}).
As $( (k^2\gamma^-_k)^2)_{k \ge 1}$ is summable and 
$\sum_{j=1}^{n} \frac{1}{k^2 - j^2} = O(\frac{\log(n+1)}{k})$  (Lemma \ref{Lemma C.2}) 
one has for any $n \le M_0 \le M$,
\begin{align*}
%\label{esti3} 
(P3) =& O \Big( \sum_{\underset{1 \le j \le n}{M_0 < k \le M}} 
\frac{1}{N^3}\big(\frac{(k^2\gamma^-_k)^2}{k^4}\frac{1}{k^2-j^2}+\frac{M^4}{N^2}\frac{1}{k^2-j^2}\big) \Big) \\
= &O\big(\frac{\log(n+1)}{N^3} (\frac{1}{M_0^5}+\frac{M^4\log M}{N^2}) \big).
\end{align*}
To estimate $(E1)$ note that for $M < k \leq N/2$,
$\omega ^N_k = O \big( \sin \frac{k\pi }{N} \big) $ (Proposition~\ref{Proposition B.1}),
$I^N_k = O \big( \frac{1}{M^2 N^2}\frac{1}{k} \big) $ (Proposition~\ref{Proposition C.1}),
and
$ \frac{1}{\tau^N_k - \tau^N_j} = O \big( \frac{1}{\cos \frac{k\pi }{N} - \cos \frac{M\pi}{N}} \big)$
(Theorem \ref{Theorem 1.1}, Lemma \ref{DeltaEV}).  By Lemma \ref{cos1} one then concludes that
   \be\label{esti4} 
(E1) = \sum _{M < k \leq \frac{N}{2}} \sum_{1 \leq j \leq n}
  O \Big( \frac{1}{M^2N^2} \frac{1}{k} \frac{ \sin\frac{k\pi }{N}}{\cos \frac{k\pi }{N}
 - \cos \frac{M\pi } {N}}\frac{1}{N^2}  \Big) = O \big( \frac{n}{M^3 N^3} \big).
   \ee
To estimate $(E2)$, observe that  
$\omega ^N_k = O(\sin\frac{k\pi }{N})$ (Proposition~\ref{Proposition B.2}),
  $  I^N_k = O \big( \frac{1}{M^2 N^2} \frac{1}{N - k} \big)$ (Proposition~\ref{Proposition C.1}),
and $\frac{1}{\tau^N_k - \tau^N_j} = O(1)$ (Theorem \ref{Theorem 1.1}, cf  Lemma \ref{DeltaEV}))
and hence  
\be\label{esti5} 
(E2) = O\Big( \sum _{\frac{N}{2} \leq k < N - M} \sum_{1 \leq j \leq n}
 \frac{1}{M^2 N^2} \frac{\sin\frac{k\pi }{N}}{N-k} \frac{1}{N^2} \Big) = O \big( \frac{n}{M^2 N^4} \big)
   \ee
where  we used that
$ \sum _{\frac{N}{2} \leq k < N-M} \frac{1}{N-k} \sin \frac{k\pi}{N} =
      \sum _{M < \ell \leq \frac{N}{2}} \frac{1}{\ell } \sin \frac{(N-\ell )\pi }{N} $
 is bounded as  $\sin \frac{(N-\ell )\pi }{N} = \sin \frac{\ell \pi }{N} \leq \frac{\ell \pi }{N}$.
Finally we consider $(E3)$. Changing the index of summation, $\ell = N - k$
and using that $\omega^N_{N-\ell} = O(\frac{\ell}{N})$ (Proposition~\ref{Proposition B.2}),
  $  I^N_{N - \ell} = O \big( \frac{1}{\ell N^2} ((\gamma^+_\ell )^2 + \frac{M^4}{N^2}) \big) \big)$ (Proposition~\ref{Proposition C.1}),
and $\frac{1}{\tau^N_{N - \ell} - \tau^N_j} = O(1)$ (Theorem \ref{Theorem 1.1}, cf  Lemma \ref{DeltaEV}))
\be
\label{esti6} 
(E3) = O\big( \sum _{1 \leq \ell \leq M} \sum_{1 \leq j \leq n} I^N_{N-\ell }
      \omega ^N_{N - \ell } \frac{1}{N^2} \big) = O \big( \frac{n}{N^5} \big) .
\ee
Combining the estimates obtained for $(E1), (E2),$ and $(E3)$ yields 
the claimed bound \eqref{remain}.
Going through the arguments of the proof one verifies that the claimed uniformity
statement holds.
\qed

\section{End of proof of Theorem \ref{thm1.5inSec6}}
\label{ProofFrequ}
 In this section we complete the proof of Theorem \ref{thm1.5inSec6} using the results 
obtained so far together with auxilary estimates obtained in Appendix \ref{tkdv}. 
First we prove the asymptotics \eqref{freq4inSec6} for the frequencies near the edges.
%%%%%%%%%%%%%%%%%%%%%%%%%%%%%%%%%%%%%%%%%%%%%%%%%%%%%%%%%%%%%%%%%%%%%%%%%%
\begin{proposition}
\label{Proposition B.3A} 
Uniformly  for $1 \le  n  \le M$ and on bounded
subsets of functions $\alpha , \beta $ in $C_0^2(\mathbb T)$ one has 
$\omega ^N_n ,\,\,\omega ^N_{N-n}= \frac{2\pi n}{N} + O \big( \frac{n^3}{N^3} \big).$
\end{proposition}
%%%%%%%%%%%%%%%%%%%%%%%%%%%%%%%%%%%%%%%%%%%%%%%%%%%%%%%%%%%%%%%%%%%%%%%%%%
\proof The asymptotics of $\omega ^N_n $ and $\omega ^N_{N-n}$ are shown in
a similar way so we concentrate on $\omega ^N_n $.  Starting point is formula \eqref{da.1.1}
with the choice $M_0 = M,$
\[
\omega_n^N  - \frac{2\pi n}{N} 
= -\frac{\pi}{4N^3}\sum_{j=1}^{n}\lambda^-_{2j-2}  - (P1) - (P2) + Error.
\]
Note that 
$\sum_{j=1}^{n}\lambda^-_{2j-2} = O\big( n^3 \big)$ (asymptotics of $(\lambda^-_j)_{j \ge 1}$)
and $(P1) = O\big( \frac{n^2}{N^3}   \big) $ (Theorem \ref{Theorem 1.1}) where
$(P1) = \frac{1}{4N^3} \sum ^n_{j=1} \int ^{\nu ^N_{2j-1}}_{\nu^N_{2j-2}} 
\arccos \big( \frac{(-1)^{N+j}}{2} \Delta _N \big( - 2 +\frac{\lambda }{4N^2} \big) \big) d\lambda.$
Concerning
$(P2) = \frac{1}{N}\sum _{I^N_k \ne 0, 1 \leq k \leq M}
  I^N_k\omega^N_k \int _{\I_n^N} \frac{\varphi ^N_k(\mu )}{i \sqrt[c]{\chi_N(\mu )}} d\mu$
note that $N^2  I^N_k = O\big( \frac{1}{k} ( (\gamma_k^-)^2 + \frac{M^4}{N^2} )\big)$ (Proposition \ref{Proposition C.1}), $\sum_{k=1}^{\infty} k^4 (\gamma_k^-)^2 < \infty$ 
($q_- \in C^2_0(\mathbb T)$), $\frac{I^N_k}{ \gamma^N_k} = O(1)$ 
(Remark \ref{Remark1}, Theorem \ref{Theorem 9.1}), and 
$\omega^N_k = O \big( \frac{k}{N} \big)$ (Proposition \ref{Da.Proposition.1}).
It then follows from  Corollary \ref{c.n.1}, \ref{c.n.2}, \ref{c.n.21}, \ref{c.n.3}, \ref{c.n.4},
and Corollary \ref{c.n.5} that $(P2) = O(\frac{1}{N^3}).$
As 
$Error = O\Big(\frac{1}{N^{3}} \big(\frac{nM^2}{N} +\frac n{M^3}\big) \Big)$ 
(Proposition \ref{Da.Proposition.1})
and thus $Error = O\big(\frac{1}{N^3}\big),$ the stated asymptotics follow.
Going through the arguments of the proof one verifies that the claimed uniformity statement holds.
\qed

\medskip

To prove the asymptotics \eqref{freq1inSec6} for the frequencies at the left edge.
the starting point is formula \eqref{da.1.1} with the choice $M_0 = L,$
\begin{equation}
\label{FORMULA}
\omega_n^N  - \frac{2\pi n}{N} 
= -\frac{\pi}{4N^3}\sum_{j=1}^{n}\lambda^-_{2j-2}  - (P1) - (P2) - (P3) + Error.
\end{equation}
First consider
$(P1) = \frac{1}{4N^3} \sum ^n_{j=1} \int ^{\nu ^N_{2j-1}}_{\nu
      ^N_{2j-2}} \arccos \big( \frac{(-1)^{N+j}}{2} \Delta _N \big( - 2 +
      \frac{\lambda }{4N^2} \big) \big) d\lambda$. 
%%%%%%%%%%%%%%%%%%%%%%%%%%%%%%%%%%%%%%%%%%%%%%%%%%%%%%%%%%%%%%%%%%%%%%%%%%%
\begin{lemma}
\label{da.l.1}
Uniformly for $1 \le n \le L$ and on bounded
subsets of functions $\alpha , \beta $ in $C_0^2(\mathbb T)$
\begin{equation}
\label{s.1}
-\pi\sum_{j=1}^{n}\lambda^-_{2j-2}  - 4N^3 (P1) = \underset{1\leq j \leq n}{\sum }
		\int ^{\lambda ^-_{2j-1}}_{\lambda ^-_{2j-2}} \frac{\lambda \dot \Delta _-(\lambda )}
                     {i \sqrt[c]{\Delta ^2_-(\lambda ) - 4}} d\lambda + O  \big(  \frac{n^2 L}{M^{1/2}} \big) .
\end{equation}
%\[ use: \forall \,\, 1 \le n \le L, \,\,\,\frac{M}{L^2} \le \frac{M}{nL} \le \frac{N}{nM^{2}} \]
\end{lemma}
%%%%%%%%%%%%%%%%%%%%%%%%%%%%%%%%%%%%%%%%%%%%%%%%%%%%%%%%%%%%%%%%%%%%%%%%%%%
\proof Integrating by parts, one gets for any $1 \le j \le n$
$$
\int^{\lambda ^-_{2j-1}}_{\lambda ^-_{2j-2}}
                     \frac{\lambda \dot \Delta _-(\lambda )}{i \sqrt[c]{\Delta _-
                     (\lambda )^2 - 4} } d\lambda 
                  = - \pi \lambda ^-_{2j-2} - 
                     \int^{\lambda ^-_{2j-1}}_{\lambda ^-_{2j-2}} \arccos \left(
                     \frac{(-1)^j}{2} \Delta (\lambda )  \right) d\lambda .
$$
To estimate the difference of 
$(Q^N_j) = \int ^{\nu ^N_{2j-1}}_{\nu ^N_{2j-2}} 
\arccos \left(\frac{(-1)^{N+j}}{2} \Delta _N \left( -2 + \frac{\lambda }{4N^2} \right) \right) d \lambda$
with the corresponding KdV integral,
$(Q^-_j) = \int ^{\lambda ^-_{2j-1}}_{\lambda ^-_{2j-2}} 
\arccos \left(\frac{(-1)^{j}}{2} \Delta _- \left(\lambda \right) \right)d \lambda$,
we use that $\arccos$ is $\frac{1}{2}-$H\"older continuous. 
Indeed, making in the integral $(Q^N_j)$ the change of variable
$[0,\delta^-_j] \to [\lambda ^-_{2j-1} ,\lambda ^-_{2j-2}], x \mapsto 
\lambda (x) =\nu ^N_{2j-2}  + h^N_j x$
with $\delta^-_j = \lambda ^-_{2j-1} -\lambda ^-_{2j-2} $
 and $h^N_j = \frac{\nu ^N_{2j-1} - \nu ^N_{2j-2}}{\lambda ^-_{2j-1} - \lambda ^-_{2j-2}}$, the difference
$(Q^N_j) - (Q^-_j)$ takes the form
\begin{equation}
\label{10.50} 
h^N_j  \int _0^{\delta^-_j} f(x) dx \,\,+\,\,
( h^N_j - 1 ) \int _0^{\delta^-_j} 
\arccos \big( \frac{(-1)^j}{2} \Delta _- (\lambda ^-_{2j-2} + x) \big) dx
\end{equation}
where 
$f(x) : = \arccos \big( \frac{(-1)^j}{2} \Delta _N(\lambda (x)) \big)
- \arccos \big( \frac{(-1)^j}{2} \Delta _- (\lambda ^-_{2j-2} + x) \big)$ 
satisfies by Theorem \ref{Theorem9.1}
\[
f(x) = \arccos \big( \frac{(-1)^j}{2} \Delta _-(\lambda (x)) + O \big( \frac{L^{2}}{M} \big) \big) 
       - \arccos \big( \frac{(-1)^j}{2} \Delta _- (\lambda ^-_{2j-2} + x) \big).
\]
Let us first consider the second term on the right hand side of \eqref{10.50}. 
As the length $\delta^-_j$ of the j'th band satisfies $1/\delta^-_j = O(1/j)$, one concludes from
Theorem \ref{Theorem 1.1} that
$h^N_j = 1 + O \big( \frac{M^2}{jN} \big)$. Hence
\[ ( h^N_j - 1) \int _0^{\delta^-_j } 
    \arccos \big( \frac{(-1)^j}{2} \Delta _- (\lambda ^-_{2j-2} + x) \big) dx
      = O \big( \frac{M^2}{N} \big) .
   \]
Towards the first term on the right hand side of \eqref{10.50}, note that
by the Lipschitz continuity of $\Delta _-$,
$\Delta _-(\lambda (x)) - \Delta _- (\lambda ^-_{2j-2} + x) =
                     O\big(\lambda (x) - \lambda ^-_{2j-2} - x \big).$
As
$
\lambda (x) - (\lambda ^-_{2j-2} + x)
= \nu ^N_{2j-2} - \lambda ^-_{2j-2} + ( h^N_j - 1) x 
$
it follows that 
\[
 \Delta _-(\lambda (x)) - \Delta _- (\lambda ^-_{2j-2} + x) 
= O \big( j \frac{M^2} {N} \big),
\]
and hence, by the $\frac{1}{2}-$H\"older continuity of arccos,
$f(x) = O(\frac{L}{M^{1/2}})$ (use that by \eqref{Fwith1/3},
$j \frac{M^2}{N} \le \frac{L^2}{M}$ ),
implying that
   \[ h^N_j \int ^{\delta^-_j }_0 f(x)dx \,\,= \,\,
     O \big( j \frac{L}{M^{1/2}} \big).
   \]
Combining the obtained estimates, the stated asymptotics follow.
Going through the arguments of the proof one verifies that the claimed uniformity statement holds.\qed

\medskip

Next we consider the term
$(P2) = \frac{1}{N} \sum _{\underset{I^N_k \ne 0}{1 \leq k \leq L}} {I^N_k}\omega^N_k  
\int _{\I_n^N}\frac{\varphi ^N_k(\mu )}{i \sqrt[c]{\chi_N(\mu )}} d\mu$ 
in \eqref{FORMULA}.
In Appendix \ref{tkdv} we show that the analogous expression for KdV 
plays an important role in describing the asymptotics of $(P2)$.
The asymptotics proved in the lemmas of Appendix \ref{tkdv} show that
$\gamma^N_k \int _{\I_n^N}\frac{\varphi ^N_k(\mu )}{i \sqrt[c]{\chi_N(\mu )}} d\mu$
stays bounded if $I^N_k$ approaches $0$. By Theorem \ref{Theorem 9.1},
$(P2)$ can thus be written  alternatively as
$(P2) = \frac{1}{N} \sum _{1 \leq k \leq L} {J^N_k}\omega^N_k  
\cdot \gamma^N_k \int _{\I_n^N}\frac{\varphi ^N_k(\mu )}{i \sqrt[c]{\chi_N(\mu )}} d\mu$. 

%%%%%%%%%%%%%%%%%%%%%%%%%%%%%%%%%%%%%%%%%%%%%%%%%%%%%%%%%%%%%%%%%%%%%%%%%%
\begin{lemma}
\label{difsum}
Uniformly for $1 \le n \le L$ and on bounded
subsets of functions $\alpha , \beta $ in $C_0^2(\mathbb T)$
\be\label{da.1.8}
(P2) = \frac{ 1}{8N^3} \sum_{\underset{I^-_k \ne 0}{1 \leq k \leq L}} I^-_k 
\int_{\I_n^-} \frac{\psi ^-_k(\lambda )}{i\sqrt[c]{\Delta^2_-(\lambda ) - 4}} d\lambda
+ O\big( \frac{1}{N^3}(\frac{L}{M^{1/2}} + \frac{1}{L^{5/2}}) \big).
\ee
\end{lemma}
%%%%%%%%%%%%%%%%%%%%%%%%%%%%%%%%%%%%%%%%%%%%%%%%%%%%%%%%%%%%%%%%%%%%%%%%%%%
\proof Using the terminology of Appendix \ref{tkdv}, decompose $(P2)$ into three parts,
\[
(P2) = \frac{1}{N} \sum _{\underset{I^N_k \ne 0, k \ne n}{1 \leq k \leq L}} {I^N_k}\omega^N_k  
\int _{\I_n^N}\frac{\varphi ^N_k(\mu )}{i \sqrt[c]{\chi_N(\mu )}} d\mu \qquad \qquad \qquad
\]
\[
\qquad +\,\, \frac{1}{N} {I^N_n}\omega^N_n  
\int _{\I_n^N \setminus C^N_n}\frac{\varphi ^N_n(\mu )}{i \sqrt[c]{\chi_N(\mu )}} d\mu\,\,
+\,\, \frac{1}{4N^3} {J^N_n}\omega^N_n  
\tilde \gamma^N_n \int _{C^N_n}\frac{\varphi ^N_n(\mu )}{i \sqrt[c]{\chi_N(\mu )}} d\mu
\]
where we used that $I^N_n = J^N_n \tilde \gamma^N_n$ and 
$\tilde \gamma^N_n = 4N^2  \gamma^N_n.$
The three parts are treated separately, but in a similar way, using the asymptotics of $I^N_k$
(Appendix \ref{actions}), $\omega^N_k$ (Proposition \ref{Proposition B.3A}), 
and the integrals (Section \ref{tkdv}). Each of them yields 
$O\big( \frac{1}{N^3}(\frac{L}{M^{1/2}} + \frac{1}{L^{5/2}}) \big)$ as an error term. 
Let us compute the asymptotics for the third part 
in more detail. Write it as
\[
{J^N_n}\omega^N_n \tilde \gamma^N_n 
\int _{C^N_n}\frac{\varphi ^N_n(\mu )}{i \sqrt[c]{\chi_N(\mu )}} d\mu
= J^-_n \gamma^-_n 
\int_{C_n^-} \frac{\psi ^-_n(\lambda )}{i\sqrt[c]{\Delta^2_-(\lambda ) - 4}} d\lambda
+ (T1) + (T2) + (T3)
\]
where $(T1) := (J^N_n - J^-_n) \omega^N_n \tilde \gamma^N_n 
\int _{C^N_n}\frac{\varphi ^N_n(\mu )d\mu}{i \sqrt[c]{\chi_N(\mu )}} $, 
$(T2) := {J^-_n} (\omega^N_n - \frac{2\pi n}{N} ) \tilde \gamma^N_n 
\int _{C^N_n}\frac{\varphi ^N_n(\mu )d\mu}{i \sqrt[c]{\chi_N(\mu )}}$,
\[
(T3) := {J^-_n} \Big(\frac{2\pi n}{N} \tilde \gamma^N_n 
\int _{C^N_n}\frac{\varphi ^N_n(\mu )}{i \sqrt[c]{\chi_N(\mu )}} d\mu -\gamma^-_n 
\int_{C_n^-} \frac{\psi ^-_n(\lambda )}{i\sqrt[c]{\Delta^2 _-(\lambda ) - 4}} d\lambda \Big).
\]
For any $1 \le k \le L,$ one has $J^N_k - J^-_k = O( \frac{L}{M^{1/2}})$
(Theorem \ref{Theorem 9.1}, Remark \ref{actionedge}), 
$\omega ^N_n = \frac{2\pi n}{N} + O \big( \frac{n^3}{N^3} \big)$ (Proposition \ref{Proposition B.3A})
and by Lemma \ref{l.n.4}, 
\[
\frac{2\pi n}{N}\tilde \gamma^N_n 
\int _{C_n^N}\frac{\varphi^N_n(\mu )}{i \sqrt[c]{\chi_N(\mu)}} d\mu=
\gamma^-_n \int_{C_n^-} \frac{\psi ^-_n(\lambda )}{i\sqrt[c]{\Delta _-(\lambda )^2 - 4}} d\lambda
+ O ( \frac{L}{M^{1/2}} + \frac{1}{L^{5/2}}).
\]
Furthermore, $J^-_n =O(n^{-5})$ (\cite{KP}), $\omega^N_n = O(\frac{n}{N})$ 
(Proposition \ref{Proposition B.3A}), and finally
$\tilde \gamma^N_n \int _{C_n^N}\frac{\varphi^N_n(\mu )}{i \sqrt[c]{\chi_N(\mu)}} d\mu
= O(\frac{N}{n})$ (Corollary \ref{c.n.4}). Combining all these estimates and using that 
$J^-_n \gamma^-_n = I^-_n /2$ one then concludes that
\[
 {J^N_n}\omega^N_n \tilde \gamma^N_n 
\int _{C^N_n}\frac{\varphi ^N_n(\mu )}{i \sqrt[c]{\chi_N(\mu )}} d\mu
= \frac{1}{2} I^-_n \int_{C_n^-} 
\frac{\psi ^-_n(\lambda )}{i\sqrt[c]{\Delta _-(\lambda )^2 - 4}} d\lambda
+ O (\frac{L}{M^{1/2}} + \frac{1}{L^{5/2}}).
\]
As already mentioned above, the first second part in the decomposition of $(P2)$ are estimated
in a similar way.
Going through the arguments of the proof one verifies that the claimed uniformity 
statement holds.
\qed
\medskip 

Combining the results obtained so far we get the following asymtptotics for the
frequencies at the left edge. 
%%%%%%%%%%%%%%%%%%%%%%%%%%%%%%%%%%%%%%%%%%%%%%%%%%%%%%%%%%%%%%%%%%%%%%
\begin{theorem}
\label{theo.fre}
Uniformly for $1 \le n \le L$ and on bounded
subsets of functions $\alpha , \beta $ in $C_0^2(\mathbb T)$, 
\begin{equation}
\label{for.fre}
\omega_n^N=\frac{2\pi n}{N}-\frac{1}{24}\frac{1}{(2N)^3}\omega^-_n
+ O\big( \frac{1}{N^3}(\frac{n^2L}{M^{1/2}} + \frac{1}{L^{5/2}}) \big)
\end{equation}
\end{theorem}  
\begin{remark}
Correspondingly, the asymptotics of $\omega_{N-n}^N$ are given by
$\omega_{N-n}^N=\frac{2\pi n}{N}-\frac{1}{24}\frac{1}{(2N)^3}\omega^+_n
+ O\big( \frac{1}{N^3}(\frac{n^2 L}{M^{1/2}} + \frac{1}{L^{5/2}}) \big).$
\end{remark}
%%%%%%%%%%%%%%%%%%%%%%%%%%%%%%%%%%%%%%%%%%%%%%%%%%%%%%%%%%%%%%%%%%%%%%%%%
\proof
By \eqref{FORMULA}, Lemma \ref{da.l.1} and Lemma \ref{difsum}, one has for any
$1 \le n \le L,$
\[
\omega_n^N  = \frac{2\pi n}{N} +   \frac{1}{4N^3}\underset{1\leq j \leq n}{\sum }
\int ^{\lambda ^-_{2j-1}}_{\lambda ^-_{2j-2}} \frac{\lambda \dot \Delta _-(\lambda )}
{i \sqrt[c]{\Delta ^2_-(\lambda ) - 4}} d\lambda  + O  \big(\frac{1}{4N^3}  \frac{n^2 L}{M^{1/2}} \big)
\]
\[
- \frac{ 1}{8N^3} \sum_{\underset{I^-_k \ne 0}{1 \leq k \leq L}} I^-_k 
\int_{\I_n^-} \frac{\psi ^-_k(\lambda )}{i\sqrt[c]{\Delta^2_-(\lambda ) - 4}} d\lambda
+ O\big( \frac{1}{N^3}(\frac{L}{M^{1/2}} + \frac{1}{L^{5/2}}) \big)
- (P3) + Error.
\]
By the formula for the KdV frequencies of Proposition \ref{scan1}, we thus get
\begin{align*}
\omega_n^N  = &\frac{2\pi n}{N} - \frac{1}{8N^3} \frac{1}{24} \omega^-_n
+ \frac{ 1}{8N^3} \frac{1}{24} \sum_{\underset{I^-_k \ne 0}{k > L}} I^-_k 
\int_{\I_n^-} \frac{\psi ^-_k(\lambda )}{i\sqrt[c]{\Delta^2_-(\lambda ) - 4}} d\lambda\\
&+ O\big( \frac{1}{N^3}(\frac{n^2 L}{M^{1/2}} + \frac{1}{L^{5/2}}) \big) - (P3) + Error.
\end{align*}
By Proposition \ref{Da.Proposition.1},
$Error = O\big( \frac{1}{N^3} (\frac{nM^2}{N} + \frac{n}{M^3}) \big)$. As
 by \eqref{Fwith1/3}, $\frac{nM^2}{N} \le \frac{n^2 L}{M^{1/2}}$ and 
$\frac{n}{M^3} \le \frac{1}{L^{5/2}}$, one concludes that
$Error = O\big( \frac{1}{N^3}(\frac{n^2 L}{M^{1/2}} + \frac{1}{L^{5/2}}) \big)$.
Taking into account in addition that $M_0 = L,$
Proposition \ref{Da.Proposition.1}
also implies that 
$(P3) = O\big( \frac{1}{N^3}(\frac{n^2 L}{M^{1/2}} + \frac{1}{L^{5/2}}) \big)$.
Arguing as in Appendix \ref{Appendix C} one sees that for any $k > L,$
$\int_{\I_n^-} \frac{\psi ^-_k(\lambda )}{i\sqrt[c]{\Delta^2_-(\lambda ) - 4}} d\lambda
= O\big( \frac{\log (n+1)}{k^2} \big)$. As $\alpha, \beta \in C^2_0(\mathbb T),$
one has $\sum_{k=1}^{\infty} k^5 I^-_k < \infty$ (\cite{KST}), implying that 
 $ \sum_{\underset{I^-_k \ne 0}{k> L}} I^-_k 
\int_{\I_n^-} \frac{\psi ^-_k(\lambda )}{i\sqrt[c]{\Delta^2_-(\lambda ) - 4}} d\lambda
= O \big( \frac{\log L}{L^5} \big)$. Hence the stated asymptotics are proved.
Going through the arguments of the proof one verifies that the claimed uniformity
statement holds.
\qed

\medskip

{\em Proof of Theorem \ref{thm1.5inSec6}}
By Theorem \ref{theo.fre} and the definitions of $\mathcal H^N_{KdV}$ and $\omega_n^-$, 
the claimed asymtotics \eqref{freq1inSec6} of the frequencies at the left edge hold. 
The asymptotics \eqref{freq2inSec6} of the frequencies at the right edge are proved 
using the symmetry of Toda chains, discussed in Appendix \ref{TodaSymmetry}.
Indeed, let
\[
\tilde \beta (x) := - \beta (-x); \quad \tilde \alpha (x) : = \alpha (-x); \quad
\tilde q_-(x) = - 2\tilde \alpha (2x) + \tilde\beta( 2x)
\]
and note that
$\tilde q_-(-x) = q_+(x)$.  Furthermore, one verifies in a straightforward way that the periodic eigenvalues
of $-\partial_x^2 + \tilde q_-(-x)$ coincide with the ones of $-\partial_x^2 + \tilde q_-(x)$. 
Hence $q_+$ and $ \tilde q_-$ have the same actions and the same KdV frequencies.
In view of Appendix \ref{TodaSymmetry} one has
$\tilde b^N_n = - b^N_{N - n} = \frac{1}{4N^2} \tilde \beta(\frac{n}{N})$ and
\[
\tilde a^N_n = a^N_{N-n-1} = 1 + \frac{1}{4N^2} \tilde \alpha(\frac{n+1}{N})
= 1 + \frac{1}{4N^2} \tilde \alpha(\frac{n}{N}) + O\big(\frac{1}{N^3}\big).
\]
By \eqref{sym7},  the claimed asymptotics \eqref{freq2inSec6} then follow from \eqref{freq1inSec6}, applied
to $\tilde \beta$ and $\tilde \alpha$.
Alternatively, one can derive the asymptotics \eqref{freq2inSec6} from \eqref{freq1inSec6}  using Corollary \ref{cor.sym}. The asymptotics \eqref{freq4inSec6} of the frequencies near the left and right edges 
are established in Proposition \ref{Proposition B.3A}. Finally the asymptotics 
\eqref{freq3inSec6}  of the frequencies in the bulk are given in \ref{Proposition B.1}. 
\qed

%%%%%%%%%%%%%%%%%%%%%%%%%%%%%%%%%%%%%%%%%%%%%%%%%%%%%%%%%%%%%%%%%%%%%%%%%%%%
%%%%%%%%%%%%%%%%%%%%%%%%%%%%%%%%%%%%%%%%%%%%%%%%%%%%%%%%%%%%%%%%%%%%%%%%%%%
%%%%%%%%%%%%%%%%%%%%%%%%%%%%%%%%%%%%%%%%%%%%%%%%%%%%%%%%%%%%%%%%%%%%%%%%%%%

\appendix

\section{Spectral estimates}
\label{AppendixSpec}

In this appendix we apply Theorem \ref{Theorem 1.1} to prove formulas involving, the eigenvalues
of the Jacobi matrices $Q^{\alpha, \beta}_N$, used at various places of the paper.

\begin{proposition}
\label{basicEstimate}
Under the same assumptions as in Theorem \ref{Theorem 1.1} the following estimates
hold uniformly in $1 \le n \le N/2$ and on bounded subsets of 
functions $\alpha, \beta \in C^2_0(\mathbb T)$:\\
(i) If $1 \le n \le M,$ then
\[
\sum_{0 < k <n} \frac{\gamma^N_k}{\lambda^N_{2n-1} - \lambda^N_{2k}},
\,
\sum_{n < k \le M} \frac{\gamma^N_k}{\lambda^N_{2k-1} - \lambda^N_{2n}}= O(\frac{1}{n}), 
\,\,\,
\sum_{M < k \le N/2} \frac{\gamma^N_k}{\lambda^N_{2k-1} - \lambda^N_{2n}} 
= O(\frac{\log M}{M^2}),
\]
\[
\sum_{N/2 < k < N - M } \frac{\gamma^N_k}{\lambda^N_{2k-1} - \lambda^N_{2n}} = O(\frac{1}{NM}),
\quad
\sum_{N-M \le k < N } \frac{\gamma^N_k}{\lambda^N_{2k-1} - \lambda^N_{2n}} 
= O(\frac{1}{N^2} + \frac{M^3}{N^3}).
\]
(ii) If $M <  n \le N/2,$ then
\[
\sum_{0 < k \le M} \frac{\gamma^N_k}{\lambda^N_{2n-1} - \lambda^N_{2k}}
= O(\frac{1}{M^2} + \frac{M\log M}{N}) 
\]
\[
\sum_{M < k < n} \frac{\gamma^N_k}{\lambda^N_{2n-1} - \lambda^N_{2k}}, \quad
\sum_{n < k \le \frac{N}{2} } \frac{\gamma^N_k}{\lambda^N_{2k-1} - \lambda^N_{2n}} 
= O(\frac{\log M}{M^2}),
\]
\[
\sum_{\frac{N}{2} < k < N-M} \frac{\gamma^N_k}{\lambda^N_{2k-1} - \lambda^N_{2n}} 
= O(\frac{\log M}{M^2}),
\quad
\sum_{N-M \le  k < N } \frac{\gamma^N_k}{\lambda^N_{2k-1} - \lambda^N_{2n}} 
= O(\frac{1}{N^2} + \frac{M^3}{N^3}) .
\]
Analogous results hold for $N/2 < n < N$.
\end{proposition}

To prove Proposition \ref{basicEstimate} we first establish the following
auxilary result.
\begin{lemma}
\label{DeltaEV}
Under the same assumptions as in Theorem \ref{Theorem 1.1} there exist $N_0 \ge 3$
and $C \ge 1$ so that
\[
\lambda^N_{2n-1}-\lambda^N_{2k} \ge \frac{1}{C}\frac{n^2 - k^2}{N^2} 
\quad \forall \,\,0 \le k < n \le \frac{N}{2}
\]
\[
\lambda^N_{2k-1}-\lambda^N_{2n} \ge \frac{1}{C}\frac{k^2 - n^2}{N^2} 
\quad \forall \,\,0 \le n < k \le \frac{N}{2}
\]
where $C$ and $N_0$ can be chosen uniformly for $k,n$ and on bounded subsets of 
$\alpha, \beta \in C^2_0(\mathbb T)$.
Similar estimates hold for $\frac{N}{2}\le k < n \le N$ respectively $\frac{N}{2} \le n < k \le N$.
\end{lemma}
{\em Proof of Lemma \ref{DeltaEV}.}
Recall from \cite[Proposition B.11]{KP} that the periodic eigenvalues of Hill operators
are compact functions of the potential on $L^2(\mathbb T)$. It then follows from the Counting
Lemma that for any $k \ge 0,$  the k'th band length 
$\lambda^-_{2k+1} - \lambda^-_{2k}$ can be estimated from below uniformly on bounded subsets of functions $\alpha, \beta \in C^2_0(\mathbb T).$ Furthermore, $\lambda^-_{2\ell}, \lambda^-_{2\ell -1} = \frac{\ell^2\pi^2}{N^2} + O(\frac{1}{\ell^2})$ as $\ell \to \infty$ where the error term
is uniform on bounded subsets of functions $\alpha, \beta \in C^2_0(\mathbb T)$ 
(cf e.g. \cite{KST} ). After these preparations we can prove the claimed estimates.
For the following we choose $N \ge N_0$ with $N_0 \ge 3$ sufficienly large.
Let us first consider the case where $0 < n \le M.$ If $0 \le k < n$, it follows from 
the above lower bound of the spectral bands and Theorem 2.1
that $\lambda^N_{2n-1}-\lambda^N_{2k} \ge \frac{1}{C}\frac{n^2 - k^2}{N^2}$ for some $C \ge 1.$ 
The case where $n < k \le M$ can be treated in the same way. If $M < k \le \frac{N}{2}$, write
\[
\lambda^N_{2k-1}-\lambda^N_{2n} =
(\lambda^N_{2k-1} + 2\cos\frac{M\pi}{N}) +(-2\cos\frac{M\pi}{N} -\lambda^N_{2M})
+ (\lambda^N_{2M} - \lambda^N_{2n})
\]
and estimate each expression in a bracket separately. We have already seen that
$\lambda^N_{2M} - \lambda^N_{2n} \ge \frac{1}{C}\frac{M^2 - n^2}{N^2}$. 
As by Theorem 2.1, $\lambda^N_{2k-1} = -2\cos\frac{k\pi}{N} +O(\frac{1}{N^2M})$
it follows from \eqref{B.91} that 
\[
\lambda^N_{2k-1} + 2\cos\frac{M\pi}{N} = (2\cos\frac{M\pi}{N} - 2\cos\frac{k\pi}{N}) 
(1 + O(\frac{1}{M^2}) ). 
\]
Hence by \eqref{B.9uno} there exists $ C \ge 1$ so that
$\lambda^N_{2k-1} + 2\cos\frac{M\pi}{N} \ge \frac{1}{C} \frac{k^2 - M^2}{N^2}$.
Finally, by Theorem \ref{Theorem 1.1} and the asymptotics of the eigenvalues $\lambda^-_\ell$
reviewed above, $\lambda^N_{2M} = -2 + \frac{M^2 \pi^2}{N^2} +O(\frac{1}{N^2M^2})$
whereas by a straightforward Taylor expansion, 
$2\cos\frac{M\pi}{N} = 2 - \frac{M^2 \pi^2}{N^2} +O(\frac{M^4}{N^4})$. 
Hence we get
$-2\cos\frac{M\pi}{N} -\lambda^N_{2M} = O(\frac{1}{N^2}(\frac{1}{M^2} + \frac{M^4}{N^2}) ).$
In view of the assumption $\eta < 1/2$ in $(F)$, 
we then conclude that for $C \ge 1$ sufficiently large,
$\lambda^N_{2k-1}-\lambda^N_{2n} \ge \frac{1}{C}\frac{k^2 - n^2}{N^2} .$
The case where $M < n \le \frac{N}{2}$ is treated in a similar fashion.
Going through the arguments of the proof one verifies that $C\ge 1$ and $N_0 \ge 3$
can be chosen uniformly for $k,n$ and for bounded subsets of functions $\alpha, \beta \in C^2_0(\mathbb T).$ The estimates for $\frac{N}{2}\le k < n \le N$ respectively $\frac{N}{2} \le n < k \le N$ are proven in the same way.
\qed

\medskip

{\em Proof of Proposition \ref{basicEstimate}.} (i) By Lemma \ref{DeltaEV} and
Theorem \ref{Theorem 1.1}
\[
\sum_{0 < k <n} \frac{\gamma^N_k}{\lambda^N_{2n-1} - \lambda^N_{2k}} 
=  O ( \sum_{0 < k <n} \frac{\gamma^-_k}{n^2 - k^2} + 
\frac{M^2}{N} \sum_{0 < k <n} \frac{1}{n^2 - k^2})
= O(\frac{1}{n})
\]
Similarly, $\sum_{n < k \le M} \frac{\gamma^N_k}{\lambda^N_{2k-1} - \lambda^N_{2n}} = O(\frac{1}{n})$ and in view of Lemma \ref{Lemma C.2} one gets
\[
\sum_{M < k \le N/2} \frac{\gamma^N_k}{\lambda^N_{2k-1} - \lambda^N_{2n}} 
\le O(\frac{1}{M} \sum_{M < k \le N/2}\frac{1}{k^2 - n^2} )
= O(\frac{\log M}{M^2}).
\]
As for any $\frac{N}{2} < k < N -M$ and $N$ sufficiently large, 
$\lambda^N_{2k-1} - \lambda^N_{2n} \ge \frac{1}{C}$ and $\gamma^N_k = O(\frac{1}{N^2M})$ 
one concludes
\[
\sum_{N/2 < k < N -M } \frac{\gamma^N_k}{\lambda^N_{2k-1} - \lambda^N_{2n}} = O(\frac{1}{NM})
\]
Finally, by similar arguments and taking into account that 
$\gamma^N_k = O(\frac{\gamma^+_{N-k}}{N^2} + \frac{M^2}{N^3})$ 
for $N-M \le k < N $ one gets
\[
\sum_{N - M \le  k < N } \frac{\gamma^N_k}{\lambda^N_{2k-1} - \lambda^N_{2n}} =
O(\frac{1}{N^2} \sum_{1 \le  \ell \le M } \gamma^+_\ell + \frac{M^3}{N^3} ) 
= O(\frac{1}{N^2} + \frac{M^3}{N^3}).
\]
(ii) The claimed estimates are proved in a similar way as the ones of item (i). We only mention
that to prove $\sum_{0 < k \le M} \frac{\gamma^N_k}{\lambda^N_{2n-1} - \lambda^N_{2k}}
= O(\frac{1}{M^2} + \frac{M\log M}{N}) $ we split the sum 
$\sum_{0 < k \le M} = \sum_{0 < k \le \frac{M}{2}} + \sum_{\frac{M}{2}< k \le M}$ and use that 
$\sum_{\ell >0} (\ell^2\gamma_\ell^-) ^2 < \infty$ and for the estimate 
$\sum_{\frac{N}{2} < k < N-M} \frac{\gamma^N_k}{\lambda^N_{2k-1} - \lambda^N_{2n}} 
= O(\frac{\log M}{M^2})$ we use that 
$\lambda^N_{2k-1} - \lambda^N_{2n}  =
 -2\cos \frac{k\pi}{N} + 2 \cos \frac{n\pi}{N}+O(\frac{1}{N^2M})$, $\cos \frac{n\pi}{N} \ge 0$
and $-\cos \frac{k\pi}{N} = \cos \frac{(N-k)\pi}{N}$
to conclude that by Lemma \ref{cos}
\[
\sum_{\frac{N}{2} < k < N-M} \frac{\gamma^N_k}{\lambda^N_{2k-1} - \lambda^N_{2n}} 
= O ( \frac{1}{N^2M} \sum_{M < k < \frac{N}{2}} \frac{1}{ \cos \frac{k\pi}{N}}) = O(\frac{\log M}{M^2}).
\]
Going through the arguments of the proofs of item (i) and (ii) one verifies that 
the estimates hold uniformly for $1 \le n \le \frac{N}{2}$ and on
 bounded subsets of functions $\alpha, \beta \in C^2_0(\mathbb T).$ The estimates for $\frac{N}{2}  < n \le N$ are proven in the same way.
\qed

%%%%%%%%%%%%%%%%%%%%%%%%%%%%%%%%%%%%%%%%%%%%%%%%%%%%%%%%%%%%%%%%%%%%%%%%%%%%%%
%%%%%%%%%%%%%%%%%%%%%%%%%%%%%%%%%%%%%%%%%%%%%%%%%%%%%%%%%%%%%%%%%%%%%%%%%%%%%%
%%%%%%%%%%%%%%%%%%%%%%%%%%%%%%%%%%%%%%%%%%%%%%%%%%%%%%%%%%%%%%%%%%%%%%%%%%%%%%

\section{Estimates of products}
\label{Appendix C}

We prove estimates on products needed at various places of the paper.
Throughout this appendix, $F : {\mathbb N} \rightarrow {\mathbb R}_{\geq 1}$ denotes a function satisfying (F) with $\eta \le 1/3$, $M = [F(N)], L = [F(M)]$, and
$w^N_k(\mu)$, $w^\pm_k(\lambda)$ the standard roots introduced in Section \ref{zeroes}.

\medskip 

In Section \ref{zeroes} and Section \ref{priori} we need the following lemmas.

\begin{lemma}
\label{integrale}
For any $1 \le j \le  M$ and $0 < k < N$, $k\not=j, j-1$ 
\begin{equation}
\label{int.fi}
 \int ^{\lambda ^N_{2j-1}} _{\lambda ^N_{2j-2}} 
\frac{\varphi ^N_k(\mu )} {\sqrt[c]{\chi_N(\mu )}}d\mu =
O\big( \frac{1}{\tau^N_k - \tau^N_j} \frac{1}{N}\big).
\end{equation}
The estimate holds uniformly in $j, k$ and uniformly on bounded subsets of 
functions $\alpha, \beta$ in $C^2_0(\mathbb T).$
\end{lemma}

\proof 
%By \eqref{phidelta}, for any $\lambda^N_{2j-2} < \mu < \lambda^N_{2j-1},$
%\[
%\frac{\varphi ^N_k(\mu )} {\sqrt[c]{\chi_N(\mu)}} =
%\frac{1}{i(-1)^{N-j}\sqrt[+]{\lambda^N_{2N-1} - \mu} \sqrt[+]{\mu - \lambda^N_0}} 
%\frac{1}{w_k^N(\mu)} \prod_{l\not=k}\frac{\mu - \sigma^{N,k}_l}{w_l^N(\mu)}.
%\]
Let us first consider the case $j = 1$. Then by assumption, $k \ge 2$. As
$\lambda_{2\ell -1}^N \le \sigma^{N,k}_\ell \le \lambda_{2\ell}^N$ one sees that
for $\lambda^N_0 \le \mu \le \lambda^N_1$, $\big\arrowvert \frac{\varphi ^N_k(\mu )} {\sqrt{\chi_N(\mu )}}\big\arrowvert $ 
is bounded by
\[  \frac{1}{\sqrt[+]{(\lambda ^N_1-\mu)(\mu - \lambda ^N_0)}}
   \frac{\sqrt[+]{\lambda ^N_2 - \lambda^N_0}} {\sqrt[+]{\lambda ^N_{2N-1}-\lambda^N_1}} \cdot
\frac{1}{\lambda^N_{2k-1} - \lambda_1^N } \cdot \prod _{\underset{\ell \not= k}{1 < \ell < N}} 
    \big( \frac{\lambda ^N_{2\ell }-\mu }{\lambda ^N_{2\ell -1} - \mu } \big) ^{1/2}.
 \]
By Lemma \ref{atif} there exists $\lambda ^N_0 <
\rho \equiv \rho^{N,k}_1 < \lambda ^N_1$ so that
\[
\int ^{\lambda ^N_1} _{\lambda ^N_0} 
\big\arrowvert \frac{\varphi ^N_k(\mu)} {\sqrt{\chi_N(\mu )}}  \big\arrowvert d\mu  
\leq \pi 
\frac{\sqrt[+]{\lambda ^N_2 - \lambda_0^N}} {\sqrt[+]{\lambda ^N_{2N-1}-\lambda_1^N}} \cdot
\frac{1}{\lambda^N_{2k-1} - \lambda_1^N } \cdot \prod _{\underset{\ell \not= k}{1 <\ell <N}} 
    \big( \frac{\lambda ^N_{2\ell }-\rho }{\lambda ^N_{2\ell -1} - \rho } \big) ^{1/2}.
\]
As
   $\frac{\lambda ^N_{2\ell} - \rho}{ \lambda ^N_{2\ell -1} -\rho } = 1
      + \frac{\gamma ^N_\ell }{\lambda ^N_{2\ell -1 } - \rho } 
     \leq \exp \big( \frac{\gamma ^N_\ell }{\lambda ^N_{2\ell -1} - \rho } \big) 
      \leq \exp
      \big( \frac{\gamma ^N_\ell }{\lambda ^N_{2\ell -1} - \lambda ^N_2} \big)$
one has
\[ \prod _{\underset{\ell > 1}{\ell \not= k}} \big( \frac{\lambda ^N_{2\ell }
      - \rho}{\lambda ^N_{2\ell -1} - \rho} \big) ^{1/2} \leq
      \exp \big( \frac{1}{2} \sum _{\underset{\ell > 1}{\ell \not= k}}
      \frac{\gamma ^N _\ell }{\lambda ^N_{2\ell -1 } - \lambda ^N_2} \big) = O (1)
\]
where for the latter estimate we used Proposition \ref{basicEstimate}. 
Furthermore, by Theorem~\ref{Theorem 1.1},  $\lambda ^N_{2N-1} - \lambda ^N_{1} =
4 + O \big( \frac{1}{N^2} \big) $, $\sqrt[+]{\lambda ^N_{2} - \lambda ^N_0}
= O \big( \frac{1}{N} \big) $, and 
   $(\lambda ^N_{2k-1} - \lambda ^N_1)^{-1} = O \left( (\tau^N_{k} - \tau^N_{1})^{-1} \right).$
This proves the claim for $j = 1.$ The case where $2 \le j \le M$ is treated in a similar way.
Going through the arguments of the proof one verifies that 
the claimed uniformity statement holds. \qed

\medskip

In Section \ref{zeroes} we introduced the rectangles $\Gamma^-_n$ with top and bottom side 
given by $[\lambda ^- _{2n-1} - 2\rho, \lambda ^-_{2n} + 2\rho ] \pm i\rho$
where $\rho > 0$ satisfies
$ \lambda ^-_{2\ell } + 3\rho < \lambda ^-_{2\ell + 1} - 3\rho \quad \forall \ell \geq 0.$
Furthermore, we chose $N_0 \ge 3$ in such a way that the error term in Theorem \ref{Theorem 1.1}
is smaller than $\rho$ for any $N \ge N_0$ (cf Lemma \ref{inV}). Then
\begin{equation}
\label{rhoBound}
\lambda^-_{2\ell-1} - \rho <   \nu^N_{2\ell -1} \le \nu^N_{2\ell} < \lambda^-_{2\ell} + \rho
\quad \forall \,\,1 \le  \ell \le M .
\end{equation}
%%%%%%%%%%%%%%%%%%%%%
\begin{lemma}
\label{Q-Estimate}
For any $\mu \equiv \mu (\lambda) = -2 + \frac{\lambda}{4N^2}$ with $\lambda \in \Gamma^-_n$
and $N \ge N_0,$
\[
\frac{\sigma ^{N,k}_\ell - \mu }{w ^N_\ell (\mu )} = 
1 + O(\frac{\gamma^N_\ell}{\tau^N_{\ell} - \tau^N_{n}}) \quad 
\]
uniformly in $1 \leq n, k \leq L$, $L < \ell < N$, 
and on bounded subsets of functions $\alpha, \beta$ in $C^2_0(\mathbb T).$
\end{lemma}
\proof Given $\mu \equiv \mu (\lambda)$ with $\lambda \in \Gamma^-_n$, write
$\frac{\sigma ^{N,k}_\ell - \mu }{w ^N_\ell (\mu )} - 1 
= \frac{\sigma ^{N,k}_\ell - \mu - w ^N_\ell (\mu ) }{w ^N_\ell (\mu )}$. We estimate
denumerator and numerator separately. In view of the definition of $\Gamma^-_n$
and the lower bound $\lambda^N_{2n} -(-2 + \frac{\lambda^-_{2n}}{4N^2}) \ge - \frac{\rho}{4N^2}$
from \eqref{rhoBound} it follows that
\[
|w ^N_\ell (\mu ) | \ge \lambda^N_{2\ell -1} - (-2 + \frac{\lambda^-_{2n} + \rho}{4N^2})
\ge  \lambda^N_{2\ell -1} - \lambda^N_{2n} -\frac{2\rho}{4N^2}.
\]
Hence $\frac{1}{w ^N_\ell (\mu )} = O (\frac{1}{\tau^N_\ell -\tau^N_n }). $ Concerning the numerator,
it follows from the definition of the standard root that
\[
w ^N_\ell (\mu ) = (\tau^N_\ell - \mu) 
\sqrt[+]{1  - \big(\frac{\gamma^N_\ell /2}{\tau^N_\ell - \mu}\big)^2}
= (\tau^N_\ell - \mu) + (\tau^N_\ell - \mu) 
\big(\sqrt[+]{1  - \big(\frac{\gamma^N_\ell /2}{\tau^N_\ell - \mu}\big)^2} - 1\big).
\]
Note that by Theorem \ref{Theorem 1.1}, $\gamma^N_\ell  = O(\frac{M^2}{N^3} + \frac{1}{MN^2})$
and by the definition of $\rho,$ $|\tau^N_\ell - \mu| \ge \frac{\rho}{N^2}$
Thus by choosing $N_0$ larger if necessary one can assume that for any $N \ge N_0$
\[
\big|\frac{\gamma^N_\ell /2}{\tau^N_\ell - \mu}\big| \le \frac{1}{2}  \qquad 
\forall \,\, L < \ell < N, \quad  \forall \,\, \mu \in \Gamma^-_n \mbox{ with } 1 \le n \le L
\]
and hence 
\[
\big|\tau^N_\ell - \mu \big| 
\big| 1 - \sqrt[+]{1  - \big(\frac{\gamma^N_\ell /2}{\tau^N_\ell - \mu}\big)^2} \,\, \big|
\le \frac{|\gamma^N_\ell|}{4} .
\]
Finally, as $\big|\sigma ^{N,k}_\ell - \tau^N_\ell \big| \le \big| \gamma^N_\ell/2 \big|$,
one has
$\sigma ^{N,k}_\ell - \mu - w ^N_\ell (\mu ) = O(\gamma^N_\ell)$. Altogether we thus proved
that 
\[
\frac{\sigma ^{N,k}_\ell - \mu }{w ^N_\ell (\mu )} = 
1 + O(\frac{\gamma^N_\ell}{\tau^N_{\ell} - \tau^N_{n}}) .
\]
Going through the arguments of the proof one verifies the claimed uniformity statement.
\qed

\medskip
%%%%%%%%%%%%%%%%%%%%%%%%%%%%%%%%%%%%%%%%%%%%%%%%%%%%%%%%%%%%%%%%%%%%%%%%
Lemma \ref{Q-Estimate} is used to get an estimate for
\[
Q^{N,L}_k(\mu )=\frac{2}{\sqrt[+]{\lambda ^N_{2N-1} - \mu }} \cdot
  \prod _{L < \ell < N} \frac{\sigma ^{N,k}_\ell - \mu }{w^N_\ell (\mu)}
\]
introduced in \eqref{P.4bis}. With the notation introduced above we have
\begin{lemma}
\label{le.peel}
For any $\mu \equiv \mu (\lambda) = -2 + \frac{\lambda}{4N^2}$ with $\lambda \in \Gamma^-_n,$
$1 \le n, k \le L$, and $N \ge N_0,$
   \begin{equation}
   \label{P.6} Q^{N,L}_k(\mu ) = 1 + O \big( \frac{1}{L^3} \big)
   \end{equation}
uniformly in $1 \leq n, k \leq L$
and on bounded subsets of functions $\alpha, \beta$ in $C^2_0(\mathbb T).$
\end{lemma}
%%%%%%%%%%%%%%%%%%%%%%%%%%%%%%%%%%%%%%%%%%%%%%%%%%%%%%%%%%%%%%%%%%%%%
\begin{remark}
\label{le.peelR}
Similarly one shows that \eqref{P.6} holds for 
$\mu \le (\lambda^N_{2L} + \lambda^N_{2L+1})/2$
uniformly in $1 \le k \leq L$
and on bounded subsets of functions $\alpha, \beta$ in $C^2_0(\mathbb T).$
\end{remark}
\proof By Theorem \ref{Theorem 1.1}, 
   $\frac{2}{\sqrt[+]{\lambda ^N_{2N-1} - \mu }} = 1 + O \left( \frac{L^2}{N^2} \right)$ 
and by Lemma \ref{Q-Estimate} 
\[
\prod _{L < \ell < N} \frac{\sigma ^{N,k}_\ell - \mu }{w ^N_\ell (\mu )} = 
\prod _{L < \ell < N} \big( 1 + O \big( \frac{\gamma ^N_\ell }{\tau^N_{\ell} - \tau^N_{n}} \big) \big)
= \exp \big( \sum _{L < \ell < N} 
\log \big( 1 + O \big( \frac{\gamma ^N_\ell }{\tau^N_{\ell} - \tau^N_{n}} \big) \big) \big)
\]
yielding the estimate
\[
\prod _{L < \ell < N} \frac{ \sigma ^{N,k}_\ell - \mu}{w ^N_\ell (\mu )}
= 1 + O \big( \sum _{L < \ell < N} \frac{\gamma ^N_\ell } {\tau^N_{\ell} - \tau^N_{n}} \big).
\]
Hence it suffices to show that
 $ \sum _{L < \ell < N} \frac{\gamma ^N_\ell }{\tau^N_{\ell} - \tau^N_{n}} 
= O \big(\frac{1}{L^3} \big) .
 $
The latter sum is split up into four parts
 $ \sum _{L < \ell \leq M} + \sum _{M < \ell \leq \frac{N}{2}} + \sum
      _{\frac{N}{2} < \ell < N - M} + \sum _{N - M \leq \ell < N} .
  $
We argue as in the proof of Proposition \ref{basicEstimate}. In particular there exists $C \ge 1$, independent of $n,$ so that
\[ \frac{\gamma ^N_\ell }{\tau^N_{\ell} - \tau^N_{n}} \leq 
 C \frac{\gamma ^-_\ell +  \frac{M^2}{N} }{\ell ^2 - L^2}, \,L < \ell \leq M,
 \quad
 \frac{\gamma ^N_\ell }{\tau^N_{\ell} - \tau^N_{n}} \leq \frac{C}{N^2} (\gamma ^+_{N-\ell }
     + \frac{M^2}{N}), \, N - M \leq \ell < N.
 \]
By Lemma \ref{Lemma C.2}  and as
$(\ell ^2 \gamma ^-_\ell )_{ \ell \geq 1}  \in \ell ^2$ due to $\alpha, \beta \in C^2_0(\mathbb T)$ one then has
   \begin{align*} \sum _{L < \ell \leq M} \frac{\gamma ^N_\ell }{\tau^N_{\ell} - \tau^N_{n}}
                     &\leq C \big( \sum _{L < \ell \leq M} \frac{\ell ^2
                     \gamma ^-_\ell }{L^3 (\ell - L)} + \sum _{L < \ell \leq
                     M} \frac{M^2}{N} \frac{1}{L} \frac{1}{\ell - L} \big) \\
                  &= O \big( \frac{1}{L^3} + \frac{M^2}{N} \frac{\log L}{L}\big) 
= O \big( \frac{1}{L^3} \big)
   \end{align*}
and as by Theorem \ref{Theorem 1.1}, $\frac{1 }{\tau^N_{\ell} - \tau^N_{n}} = O(1)$
for any $N - M \leq \ell < N$ and $\eta \le 1/3,$
   \[ \sum _{N - M \leq \ell < N} \frac{\gamma ^N_\ell }{\tau^N_{\ell} - \tau^N_{n}} =
      O \left( \frac{1}{N^2} \right).
   \]
Finally, for $M < \ell \leq \frac{N}{2}$ one has $\gamma ^N_\ell = O(N^{-2} M^{-1})$ and
   \[ \tau^N_{\ell} - \tau^N_{n} \geq \tau^N_{\ell} - \tau^N_{M}
        = - 2\cos \frac{\ell \pi }{N} + 2\cos \frac{M\pi }{N} + O \left( N^{-2} M^{-1} \right) .
   \]
Using 
   $ - \cos \frac{\ell \pi }{N} + \cos \frac{M\pi }{N} \geq \pi  \frac{\ell ^2 - M^2}{N^2}
   $
( cf \eqref{B.9uno}) and Lemma \ref{atif} one sees that
   \[ \sum _{M < \ell \leq \frac{N}{2}} \frac{\gamma ^N_\ell }{\tau^N_{\ell} - \tau^N_{n}}
      = O \big( \frac{1}{N^2M} \sum _{\ell > M} \frac{N^2}{\ell ^2 -M^2} \big)= 
     O \big( \frac{\log M}{M^2} \big) .
   \]
Note that $\frac{\log M}{M^2} =O(\frac{1}{L^3})$. Finally, as $(\tau^N_{\ell} - \tau^N_{n})^{-1} = O(1)$ for $\frac{N}{2} < \ell < N - M$, one has
 $ \sum _{\frac{N}{2} < \ell < N - M} \frac{\gamma ^N_\ell }{\tau^N_{\ell} - \tau^N_{n}}
      = O\big( \frac{1}{NM} \big) .
  $
Going through the arguments of the proof one verifies the claimed uniformity statement.
\qed

\medskip
In a straightforward way one estimates  $Q^{-,L}_k(\lambda ) =
\prod _{\ell > L} \frac{\sigma ^{-,k} _\ell - \lambda }{w^-_\ell (\lambda )},$
defined by \eqref{P.4bis.1} and obtains the following result.
%%%%%%%%%%%%%%%%%%%%%%%%%%%%%%%%%%%%%%%%%%%%%%%%%%%%%%%%%%%%%%%%%%%%%%%
\begin{lemma}
\label{le.peel.1} Uniformly for $1 \le n, k \le L$, $\lambda \in \Gamma^-_n$,
and uniformly on bounded subsets of functions $\alpha, \beta$ in $C^2_0(\mathbb T)$
   \begin{equation}
    \label{1.10bis} Q^{-,L}_k(\lambda ) = 1 + O\big( \sum _{\ell > L}
                    \frac{\gamma ^-_\ell }{\ell ^2- L^2} \big) = 1 + O \big( \frac{1}{L^3} \big) .
   \end{equation}
 \end{lemma}
\begin{remark}
\label{le.peel.1R}
Similarly, one shows that \eqref{1.10bis} holds for 
$\lambda \le (\lambda^-_{2L} + \lambda^-_{2L+1})/2$
uniformly in $1 \le k \leq L$
and on bounded subsets of functions $\alpha, \beta$ in $C^2_0(\mathbb T).$
\end{remark}

Using the lemmas above we now estimate $\frak G^{N,L}_k(\lambda )$, introduced in \eqref{gnk0},
\[
\frak G^{N,L}_k(\lambda ) =  \frac{Q^{-,L}_k(\lambda )}
      {Q^{N,L}_k(\mu )} \cdot \frac{2N\sqrt[+]{\mu - \lambda ^N_0}}{\sqrt[+]{\lambda -
      \lambda ^-_0}} \cdot \frac{4N^2 w^N_k(\mu )}{w^-_k(\lambda )}
\] 
where $\mu \equiv \mu(\lambda)= -2 + \frac{\lambda}{4N^2}$.

\begin{lemma}
\label{G}
Uniformly for $1 \le n, k \le L$ with $n \ne k,$ $\lambda\in \Gamma^-_n$, and
uniformly on bounded subsets of functions $\alpha, \beta$ in $C^2_0(\mathbb T)$, one has 
\begin{equation}
\label{gnk}
\frak G^{N,L}_k(\lambda)=1+O\big(\frac{1}{L^3}\big).
\end{equation}
\end{lemma}
%%%%%%%%%%%%%%%%%%%%%%%%%%%%%%%%%%%%%%%%%%%%%%%%%%%%%%%%%%%%%%%%%%
\begin{remark}
Similarly, one shows that \eqref{gnk} holds for 
$\lambda \in  [\lambda^-_{1}, \lambda^-_{2k-2}] \cup  [\lambda^-_{2k+1}, \lambda^-_{2L}]$
uniformly in $1 \le k \leq L$
and on bounded subsets of functions $\alpha, \beta$ in $C^2_0(\mathbb T).$
\end{remark}
\proof 
By Theorem \ref{Theorem 1.1} one has, with $\nu^N_j = 4N^2(\lambda^N_j + 2),$
\[
\frac{4N^2 w^N_k(\mu )}{w^-_k(\lambda )} 
  = \sqrt[+]{ \frac{(\nu ^N_{2k} - \lambda )(\nu ^N_{2k-1}-\lambda )}{(\lambda ^-_{2k}
   - \lambda )(\lambda ^-_{2k-1}-\lambda )}} = 1 + O \big( \frac{M^2}{N} \frac{1}{|k^2 - n^2|} \big) 
\]
and
\[
\frac{2N\sqrt[+]{\mu - \lambda ^N_0}}{\sqrt[+]{\lambda - \lambda ^-_0}}
= \sqrt[+]{1 + \frac{\lambda ^-_0 - \nu ^N_0}{\lambda -\lambda ^-_0}}
= 1 + O \big( \frac{M^2}{N} \frac{1}{n^2} \big).
\]
The claimed estimate then follows from Lemma \ref{le.peel}, Lemma \ref{le.peel.1},
and \eqref{Fwith1/3}.
Going through the arguments of the proof one verifies the uniformity statement.\qed

%%%%%%%%%%%%%%%%%%%%%%%%%%%%%%%%%%%%%%%%%%%%%%%%%%%%%%%%%%%%%%%%%%%%%%%%%%
%%%%%%%%%%%%%%%%%%%%%%%%%%%%%%%%%%%%%%%%%%%%%%%%%%%%%%%%%%%%%%%%%%%%%%%%%%
%%%%%%%%%%%%%%%%%%%%%%%%%%%%%%%%%%%%%%%%%%%%%%%%%%%%%%%%%%%%%%%%%%%%%%%%%%%

\section{Auxilary lemmas}
\label{peeling}

For the convenience of the reader, we collect in this appendix elementary lemmas
which are used at various places of the paper. The first two lemmas concern uniform bounds
for the H\"older continuity of some special functions. 

\begin{lemma}
\label{Hoelder} 
{\rm(i)} $0 \le \sqrt{x} - \sqrt{y} \le \sqrt{x-y} \quad \forall \,\, x \ge y \ge 0;$

{\rm(ii)} $0 \le \log x - \log y \le  x - y \quad \forall \,\,x \ge  y \ge 1;$

{\rm(iii)} $0  \le \log (x + \sqrt{x^2 -1}) - \log (y + \sqrt{y^2 -1}) \le 2 \sqrt{x+y} \sqrt{x-y} \,\,\,
\forall \,\ x \ge  y \ge 1.$
\end{lemma}

\begin{remark}
Note that for $x \ge 1,$ the principal branch of $arcos(x)$ is given by 
$arcosh(x) = \log (x + \sqrt{x^2 -1}).$
\end{remark}

{\it Proof:} (i) The claimed estimate is obtained from combining the following 
inequalities, valid for any $x \ge y \ge 0,$
\[
x - y = (\sqrt{x} - \sqrt{y})(\sqrt{x} + \sqrt{y}), \quad x -y \le \sqrt{x-y}\sqrt{x+y},
\quad \sqrt{x+y} \le \sqrt{x} + \sqrt{y}.
\]
 (ii) For any $x \ge y \ge 1$ one has
\[
0 \le \log x - \log y = \int_0^1 \partial_t \log(y + t(x-y)) dt \le x - y .
\]
(iii) By applying first (ii) and then (i) one sees that for any $x \ge y \ge 1,$
\[
\log (x + \sqrt{x^2 -1}) - \log (y + \sqrt{y^2 -1})  \le (x-y) + \sqrt{x^2 - y^2}.
\]
As $\sqrt{x^2 - y^2} = \sqrt{x-y} \sqrt{x+y},$ the claimed estimate follows.
\hspace*{\fill }$\square $

\begin{lemma}
\label{Lemma C.3bis} For any $a, b> 0$ and $ x,y \geq 0$,
   \[ \Big\arrowvert \frac{a}{\sqrt{a + x}} - \frac{b}{\sqrt{b + y}}
      \Big\arrowvert \leq |a - b|^{1/2} +  |x - y|^{1/2} .
   \]
By a limiting argument, the inequality continues to hold for $a,b \ge 0.$
\end{lemma}

{\it Proof:} For any $a, b > 0$ and $x,y \ge 0$ one has
   \begin{align*} \Big\arrowvert \frac{a}{\sqrt{a + x}} - \frac{b}{\sqrt{b + y}}
      \Big\arrowvert ^2 &\leq \Big\arrowvert \left( \frac{a}{\sqrt{a + x}} -
      \frac{b}{\sqrt{b + y}} \right) \left( \frac{a}{\sqrt{a + x}} +
      \frac{b}{\sqrt{b + y}} \right) \Big\arrowvert\\
  &= \Big\arrowvert \frac{a^2}{a + x} - \frac{b^2}{b + y} \Big\arrowvert =
      \Big\arrowvert \frac{a^2b - ab^2 + a^2y - b^2x}{ab + xb + ay + xy} \Big\arrowvert 
\end{align*}
As $a^2b - ab^2 + a^2y - b^2x =  (ab + xb + ay)(a - b) + ab(y - x)$
   one gets 
$ \arrowvert \frac{a}{\sqrt{a + x}} - \frac{b}{\sqrt{b + y}} \arrowvert ^2 \leq |a - b| + |y - x|$ 
or
   \[ \Big\arrowvert \frac{a}{\sqrt{a + x}} - \frac{b}{\sqrt{b + y}}
      \Big\arrowvert \leq  |a - b|^{1/2} + |y - x|^{1/2} .
   \]
\hspace*{\fill }$\square $	
\begin{lemma}
\label{cos} 
For any
$M < n,k < N - M, k \not= n$
   \begin{equation}
   \label{B.91}  |\cos \frac{n\pi}{N} - \cos \frac{k\pi}{N}| \geq \frac{M}{N} \cdot \frac{\pi }{N}
   \end{equation}
whereas for any $0 < \ell \leq \frac{N}{2}$ and $0 \le j < \ell$
   \begin{equation}
   \label{B.9uno}  
		\cos \frac{j\pi}{N} - \cos \frac{\ell \pi}{N} \ge \frac{(\ell^2 - j^2)\pi }{N^2}.
   \end{equation}
Furthermore, for any $M < n < N-M$ one has
   \begin{equation}
   \label{B.9bis} \sum _{\underset{M < k < N - M}{k \not= n}} \frac{1}{|\cos \frac{n
                  \pi }{N} - \cos \frac{k\pi }{N}|} = O(N^2 \frac{\log M}{M})
   \end{equation}
\end{lemma}
\proof
The difference $\cos \frac{n\pi }{N} - \cos \frac{k\pi }{N}$ is bounded in absolute
value from below by $\min _\pm |\cos \frac{n\pi }{N} - \cos \frac{(n\pm 1)\pi }
{N}|$. Writing $|\cos x - \cos y| = |\int ^y_x \sin t dt|$ and using that
   \[ \sin x \geq \sin \frac{M\pi }{N} \geq \frac{2}{\pi } \cdot \frac{M\pi }{N}
      \qquad  \forall  \,\,\,x \in 
      [\frac{M\pi }{N} , \pi - \frac{M\pi }{N}]
   \]
\eqref{B.91} follows. 
Next, for any 
$0 < \ell \leq \frac{N}{2}$ and $0 \le j < \ell$ one has
 \[\cos \frac{j\pi}{N} - \cos \frac{\ell \pi}{N}  = 
     \int ^{\frac{\ell\pi}{N}}_{\frac{j\pi }{N}} \sin x dx \geq  
     \frac{x^2}{\pi } \big\arrowvert^{\frac{\ell \pi }{N}}_{\frac{j\pi }{N}}  
     = \frac{(\ell^2 - j^2)\pi }{N^2}
 \]
which proves \eqref{B.9uno}.
Concerning \eqref{B.9bis} let us concentrate on the
case where $M < n \leq \frac{N}{2}$.  (The case where $\frac{N}{2}
\leq n < N - M$ is treated in the same way.) Then one sees in a
straightforward way that
   $\sum_{M < k < N - M, \,\,k\not= n} \cdots 
  \leq \ 3 \sum_{M < k \leq N/2,\,\, k \not= n} \cdots .$
Hence \eqref{B.9bis} follows from \eqref{B.9uno} and Lemma \ref{Lemma C.2}.
\qed

\begin{lemma}
\label{cos1}
The following estimate holds
\[ \sum _{M < k \leq \frac{N}{2}} \frac{\sin \frac{k\pi}{N}}{\cos \frac{M\pi }
      {N} - \cos \frac{k\pi }{N}} = O(N).
   \]
\end{lemma}
\proof
Note that
   \[ \sum _{M < k \leq \frac{N}{2}} \frac{\sin \frac{k\pi}{N}}{\cos \frac{M\pi }
      {N} - \cos \frac{k\pi }{N}} \leq \frac{\sin \frac{(M+1)\pi}{N}}{\cos \frac{M\pi }
      {N} - \cos \frac{(M+1)\pi }{N}} + \int ^{\frac{\pi}{2}}_{\frac{(M+1)\pi}{N}}
      \frac{\sin x}{\cos \frac{M\pi }{N} - \cos x} dx 
   \]
and
\[
      \int ^{\pi /2}_{\frac{(M+1)\pi}{N}} \frac{\sin x}{\cos \frac{M\pi }{N} - \cos x}
         dx = \log \left( \cos \frac{M\pi }{N} - \cos x \right) \Big\arrowvert ^{\pi /2}
         _{\frac{(M + 1)\pi }{N}} 
\]
\[
\le  - \log \big( \cos \frac{M\pi } {N} - 
     \cos \frac{(M + 1)\pi}{N} \big) .
 \]

By \eqref{B.9uno} one has
$\cos \frac{M\pi }{N} - \cos \frac{(M + 1) \pi}{N}  \ge \frac{(2M+1)\pi}{N^2}$ and hence
   \[ - \log \big( \cos \frac{M\pi }{N} - \cos \frac{(M+1)\pi }{N} \big) \leq 2 \log N .
   \]
Furthermore, together with the estimate $ 0 \le \sin \frac{(M + 1)\pi }{N} \le \frac{(M+1)\pi}{N}  $ 
one gets
\[\frac{\sin ((M + 1)\pi / N)}{\cos \frac{M\pi }{N} - \cos \frac{(M + 1)\pi }{N}}  \le N.
\]
Combining the two latter estimates the claim follows.
\qed

\begin{lemma}\label{atif} (Mean value theorem in integral form)
Let $f$ be continuous and $g\geq 0$ integrable on the interval $[a,b]$ with $a < b$. Then 
there exists $a < y < b$ so that $\int_a^bf(x)g(x)dx=f(y)\int_a^bg(x)dx.$
\end{lemma}

\begin{lemma}
\label{Lemma C.2} For any $k \in \mathbb Z_{\ge 1}$ and any $1 \le n < k,$
$
\sum\limits _{0 \le j \le  n } \frac{1}{k^2 - j^2}  = O \big( \frac{\log( n +1)} {k} \big).$
Similarly, for any $n > k,$
$\sum _{j \ge n } \frac{1}{j^2 - k^2} = O \big( \frac{\log( k +1)} {k} \big).$
\end{lemma}
\proof
Note that for any $n > k$,
   $ \sum\limits _{j \geq n} \frac{1}{j^2 - k^2} \leq \frac{1}{2k} 
+ \int ^\infty _{n} \frac{1}{x^2 - k^2} dx .$
   Using that 
 $\frac{1}{x^2 - k^2} = \big( \frac{1}{x - k} - \frac{1}{x + k} \big)
      \frac{1}{2k} = \frac{1}{2k} \frac{d}{dx} \log \frac{x - k}{x + k}
= \frac{1}{2k} \frac{d}{dx} \log (1 - \frac{2k}{x + k})$
one gets
   \begin{align*} \int ^N_{n} \frac{1}{x^2-k^2} dx 
		   &= \frac{1}{2k} \log \big( 1 - \frac{2k}{N + k} \big) -
                     \frac{1}{2k} \log \big( 1 - \frac{2k}{n+k} \big)  \\
                  &= \frac{1}{2k} \log \big( 1 - \frac{2k}{N + k} \big) +
                     \frac{1}{2k} \log (1 +\frac{2k}{n-k}  ) .
   \end{align*}
Hence $\lim _{N \rightarrow \infty } \int ^N_{n} \frac{1}{x^2-k^2}dx = 
\frac{1}{2k} \log (1 +\frac{2k}{n-k})  \le \frac{1}{2k} \log (1 +2k) $. 
The sum $\sum\limits _{0 \le j \le n } \frac{1}{k^2 - j^2}$
can be estimated in a similar fashion.
\qed
%\hspace*{\fill }$\square $
%\end{proof}	

%%%%%%%%%%%%%%%%%%%%%%%%%%%%%%%%%%%%%%%%%%%%%%%%%%%%%%%%%%%%%%%%%%%%%%%%%%%%%%
%%%%%%%%%%%%%%%%%%%%%%%%%%%%%%%%%%%%%%%%%%%%%%%%%%%%%%%%%%%%%%%%%%%%%%%%%%%%%%%
%%%%%%%%%%%%%%%%%%%%%%%%%%%%%%%%%%%%%%%%%%%%%%%%%%%%%%%%%%%%%%%%%%%%%%%%%%%%%%%

\section{Auxilary estimates}\label{tkdv}

In this appendix we prove auxilary estimates needed in Section \ref{ProofFrequ} to derive the claimed asymptotics of $\omega_n^N$
for $1 \le n \le L$. In view of \eqref{remain}, \eqref{da.1.4}, and the assumptions about $F,$ formula \eqref{da.1.1} reads
in the case $1 \le n \le L$ and $M_0 = L$ as
$\omega_n^N=
\frac{2\pi n}{N} -\frac{\pi}{4N^3}\sum_{j=1}^{n}\lambda^-_{2j-2} - (P1) - (P2) + Error_L$ 
where
\begin{equation}
\label{da.1.3}
(P2) = \frac{1}{N} \sum _{\underset{I^N_k \ne 0}{1 \leq k \leq L}} {I^N_k}\omega^N_k  \int _{\I_n^N}
           \frac{\varphi ^N_k(\mu )}{i \sqrt[c]{\chi_N(\mu )}} d\mu,
\end{equation}
\[
Error_L = O\Big(\frac{1}{N^{3}} \big(\frac{\log L}{L^5} + \frac{nM^2}{N} +\frac n{M^3}\big) \Big).
\]
We will need estimates for the integral in \eqref{da.1.3}, relating it
to a corresponding one of KdV. The integral in \eqref{da.1.3} is 
a sum of convergent improper Riemann integrals on intervals at the endpoints of which 
the integrand might be singular. More specifically, the improper Riemann integrals 
$\int _{\lambda^N_{2k - 2}}^{\lambda^N_{2k-1}}
\frac{\varphi ^N_k(\mu )}{i \sqrt[c]{\chi_N(\mu )}} d\mu$ and
$\int _{\lambda^N_{2k}}^{\lambda^N_{2k+1}}
      \frac{\varphi ^N_k(\mu )}{i \sqrt[c]{\chi_N(\mu )}} d\mu$ 
can be shown to be of the order $O(\frac{N}{k} \log\gamma^N_k)$
and to diverge if $I^N_k = 0.$ On the other hand, $i \sqrt[c]{\chi_N(\mu )}$
is real valued on both intervals $[\lambda^N_{2k - 2}, \lambda^N_{2k-1}]$
and $[\lambda^N_{2k} , \lambda^N_{2k+1}]$ and has opposite signs 
which will lead to cancelations. 
%Indeed for $\mu \in [\lambda^N_{2k} , \lambda^N_{2k+1}]$,
%$i \sqrt[c]{\chi_N(\mu )} = (-1)^{N+k} \sqrt[c]{-\chi_N(\mu )}$.
Hence a careful analysis is needed to obtain estimates which are uniform on bounded subsets of functions $\alpha , \beta $ in $C_0^2(\mathbb T)$. We split the
domain $\I^N_n$ of the integral in \eqref{da.1.3} as follows. Fix 
$\rho > 0$  in such a way that  
$\lambda ^- _{2j} + 3\rho < \lambda ^- _{2j+1} - 3\rho$ for any $ j \geq 0$ (cf \eqref{rhoBound}),
set $\rho ^N = \frac{\rho }{4N^2}$ and introduce
\begin{align*}
A_j^N:=&[\lambda_{2j-1}^N-\rho^N, \lambda_{2j-1}^N]\cup
[\lambda_{2j}^N, \lambda_{2j}^N+\rho^N],   \,\, \quad 0 <  j < n
\\
B_j^N:=&[\lambda_{2j}^N+\rho^N, \lambda_{2j+1}^N-\rho^N], \,\,\quad 0 < j < n
\\
A_0^N :=&[\lambda_0^N, \lambda_{1}^N-\rho^N], \,\, \quad
C_{n}^N :=[\lambda_{2n-1}^N - \rho^N, \lambda_{2n-1}^N] 
\end{align*}
so that $\I^N_n=\big(\cup_{j=0}^{n-1} A_j^N\big)\cup
\big(\cup_{j=1}^{n-1} B_j^N\big) \cup C_{n}^N.$
Correspondingly  define
\begin{align*}
A_j^-:=&[\lambda_{2j-1}^--\rho, \lambda_{2j-1}^-]\cup
[\lambda_{2j}^-, \lambda_{2j}^-+\rho],   \ ,\quad 0 < j < n
\\
B_j^- :=&[\lambda_{2j}^-+\rho, \lambda_{2j+1}^--\rho], \,\, \quad 0 < j < n
\\
A_0^- :=&[\lambda_0^-,\lambda_1^- -\rho]\ ,\quad
C^-_{n} :=[\lambda^-_{2n-1}-\rho,\lambda^-_{2n-1}]
\end{align*}
and set $\I_n^-= \big(\bigcup_{j=0}^{n-1} A_j^-\big)\cup
\big(\cup_{j=1}^{n-1} B_j^-\big) \cup C^-_{n}.$
In the lemmas below we will derive the asymptotics of
the integrals over the intervals $(A_j^N)_{0 \le  j < n}$, $(B_j^N)_{0 < j < n}$, and $C_{n}^N$
   in terms of corresponding integrals for KdV over the intervals 
$(A_j^-)_{0 \le j < n}$, $(B_j^-)_{0 < j < n},$ and $C^-_{n}$. We begin with the analysis of the integral over the interval $A^M_0= [\lambda_0^N, \lambda_{1}^N-\rho^N]$.

%%%%%%%%%%%%%%%%%%%%%%%%%%%%%%%%%%%%%%%%%%%%%%%%%%%%%%%%%%%%%%%%%%%
\begin{lemma}
\label{l.n.1} Uniformly for $1 \le k \le L$ and on bounded
subsets of functions $\alpha , \beta $ in $C_0^2(\mathbb T)$
\begin{equation}
\label{da.i.1}
\int _{A_0^N}
      \frac{\varphi ^N_k(\mu )}{i \sqrt[c]{\chi_N(\mu )}} d\mu = 
\frac{N}{2\pi k} \int_{A_0^-} \frac{\psi^-_k(\lambda)}{i\sqrt[c]{\Delta^2_-(\lambda)-4 }}d\lambda 
+ O\big(\frac{N}{k^2}\frac{1}{L^{5/2}} \big).
\end{equation}
\end{lemma}
\proof For $\mu\in A^N_0$, and $1 \le k \le L$ one has by  \eqref{wN}, \eqref{phidelta},
\eqref{P.4bis}, \eqref{P.4bis.1} 
\[
 \frac{\varphi ^N_k(\mu )}{i \sqrt[c]{\chi_N(\mu )}}=\frac{1}{2}
 \Big( \prod _{\underset{\ell \not= k}{1 \leq \ell \leq L}} \frac{\sigma ^{N,k}_\ell - \mu }
  {w^N_\ell (\mu )} \Big) \frac{1}{\sqrt[+]{\mu - \lambda ^N_0}}\frac{1}{w^N_k(\mu )}  Q^{N,L}_k(\mu).
\]
Here we used that  
$\varphi^N_k(\mu ) = (-1)^{N-2} \prod _{\underset{\ell \not= k}{1 \leq \ell \leq L}} 
(\sigma ^{N,k}_\ell - \mu )$ and that by \eqref{signChiN}
$(-1)^N i \sqrt[c]{\chi_N(\mu )} > 0 $ for any $\mu \in (\lambda^N_0 , \lambda^N_1)$
to get the correct sign.
Setting $\mu \equiv \mu(\lambda) = -2 + \frac{\lambda}{4N^2}$ and writing $\frac{1}{w^N_k(\mu )} \prod _{\underset{\ell \not= k}{1 \leq \ell \leq L}} 
\frac{\sigma ^{N,k}_\ell - \mu } {w^N_\ell (\mu )} $ as
\[
\frac{4N^2}{w^-_k(\lambda )} \prod _{\underset{\ell \not= k}{1 \leq \ell \leq L}} 
\frac{\sigma ^{-,k}_\ell - \lambda } {w^-_\ell (\lambda )}  \cdot
\prod _{\underset{\ell \not= k}{1 \leq \ell \leq L}} 
\frac{4N^2(\sigma ^{N,k}_\ell - \mu )} {\sigma^{-,k}_\ell - \lambda} \cdot
\prod _{1 \leq \ell \leq L} 
\frac{w ^-_\ell(\lambda) } {4N^2w^N_\ell (\mu )} 
\]
one argues as in the proof of \eqref{Z.10} to conclude that for $\lambda \le \lambda^-_1 - \rho/2,$
$\prod _{1 \leq \ell \leq L} \frac{w ^-_\ell(\lambda) } {4N^2w^N_\ell (\mu )} = 1 + O(\frac{M^2}{N})
$ 
and, taking into account also Theorem \ref{theorem2.1},
$\prod _{\underset{\ell \not= k}{1 \leq \ell \leq L}} 
\frac{4N^2(\sigma ^{N,k}_\ell - \mu )} {\sigma^{-,k}_\ell - \lambda} = 1 + O(\frac{1}{L^{5/2}}).$
By Lemma \ref{le.peel} (cf Remark \ref{le.peelR}) 
and Lemma \ref{le.peel.1} (cf Remark \ref{le.peel.1R}),
$Q^{N,L}_k(\mu) = 1 + O(\frac{1}{L^3})$ respectively
$Q^{-,L}_k(\lambda)^{-1} = 1 + O(\frac{1}{L^3} )$. By \eqref{Fwith1/3} one then gets uniformly for 
$\lambda \le \lambda^-_1 - \rho/2$,
\[
\frac{\varphi ^N_k(\mu )}{i \sqrt[c]{\chi_N(\mu )}}
=  \frac{1}{w^-_k(\lambda )} \big( \prod _{\underset{\ell \not=k}{0 <\ell <\infty}} 
\frac{\sigma ^{-,k}_\ell - \lambda }  {w^-_\ell (\lambda )}  \big)
\frac{2N^2}{\sqrt[+]{\mu - \lambda ^N_0}}  \Big(1+O\big(\frac{1}{L^{5/2}} \big)\Big).
\]
Note that uniformly for $\lambda\in A^-_0$,
$ \prod _{\underset{\ell \not=k}{0 <\ell <\infty}} 
\frac{\sigma ^{-,k}_\ell - \lambda }  {w^-_\ell (\lambda )} = O(1)$ 
and, by the asymptotics of the periodic eigenvalues of $H_-,$ $\tau^-_k = (\lambda^-_{2k-1} + \lambda^-_{2k})/2 = 4k^2\pi^2 + O(1),$  whence
$w^-_k(\lambda)^{-1} = O(k^{-2})$. Thus for $\lambda \le \lambda^-_1 - \rho/2,$
\be\label{fcor}
f_0^k(\lambda) := \frac{1}{w^-_k(\lambda )} \prod _{\underset{\ell \not=k}{0 <\ell <\infty}} 
\frac{\sigma ^{-,k}_\ell - \lambda }  {w^-_\ell (\lambda )}  =O(k^{-2}).
\ee
To prove the claimed asymptotics of the integral 
$\int _{A_0^N} \frac{\varphi ^N_k(\mu )}{i \sqrt[c]{\chi_N(\mu )}} d\mu$
introduce the change of variable 
$[\lambda^-_0, \lambda^-_1 - \rho] \to [\lambda^N_0, \lambda^N_1 - \rho^N]$,
given by
\begin{align}
\label{ch.va}
x \mapsto \mu=\lambda_0^N+\frac{h^N_0}{4N^2}(x-\lambda^-_0),
\qquad
h^N_0:= 4N^2 \frac{\lambda_1^N-\rho^N-\lambda^N_0}{\lambda_1^--\rho-\lambda^-_0} > 0.
\end{align}
Then $h^N_0=1+O\big(\frac{M^2}{N}\big)$
(cf Theorem \ref{Theorem 1.1})
and for any $x \in [\lambda^-_0, \lambda^-_1 - \rho]$,
\begin{align*}
\lambda(x) := (\mu(x) +2)4N^2 = x + O\big(\frac{M^2}{N}\big),\qquad
\frac{d\mu}{\sqrt[+]{\mu-\lambda^N_0}}
=\frac{\sqrt[+]{h^N_0}}{2N}\frac{dx}{\sqrt[+]{x-\lambda_0^-}}
\ .
\end{align*}
Using \eqref{fcor} and \eqref{Fwith1/3} we get
\begin{align}
\int _{A_0^N} \frac{\varphi ^N_k(\mu )}{i \sqrt[c]{\chi_N(\mu )}} d\mu
\label{da.i.n}
=N\int_{\lambda_0^-}^{\lambda_1^--\rho}f_0^k(\lambda(x))
\frac{\sqrt[+]{h_0}}{\sqrt[+]{x-\lambda_0^-}} dx +
O\Big(\frac{N}{k^2}\frac{1}{L^{5/2}}\Big).
\end{align}
As by \eqref{psidelta} and \eqref{signDelta-},
$ \frac{1}{2\pi k}\frac{\psi^-_k(\lambda)} {i\sqrt[c]{\Delta^2_-(\lambda)-4 }} = \frac{1}{\sqrt[+]{\lambda - \lambda^-_0}} f^k_0(\lambda) $ 
and as $\lambda (x) = h^N_0 x + c_N$ with $c_N = O(\frac{M^2}{N})$  the stated asymptotics then follow by using once more \eqref{fcor}.
Going through the arguments of the proof one verifies that the claimed uniformity statement holds.
\qed

\medskip 
For the proof of Theorem \ref{thm1.5inSec6} we need bounds for 
$\int _{A_0^N} \frac{\varphi ^N_k(\mu )}{i \sqrt[c]{\chi_N(\mu )}} d\mu$. 
In particular, in the proof of Proposition \ref{Proposition B.3A} we need a bound for 
these integrals also for $k$ in the range $L < k \le M.$
By similar arguments as in the proof of Lemma \ref{l.n.1} one gets the following result.
\begin{corollary}
\label{c.n.1}
Uniformly for $1 \le k \le M$ and on bounded
subsets of functions $\alpha , \beta $ in $C_0^2(\mathbb T)$
\begin{equation}
\label{c.n.e1}
\int _{A_0^N} \frac{\varphi ^N_k(\mu )}{i \sqrt[c]{\chi_N(\mu )}} d\mu
=O\left(\frac{N}{k^2} \right).
\end{equation}
\end{corollary}
%%%%%%%%%%%%%%%%%%%%%%%%%%%%%%%%%%%%%%%%%%%%%%%%%%%%%%%%%%%%%%%%%%%%
\proof Going through the arguments of the proof of Lemma \ref{l.n.1} one sees that
$\prod _{1 \leq \ell \leq M} \frac{w ^-_\ell(\lambda) } {4N^2w^N_\ell (\mu )} = 1 + O(\frac{M^2}{N})$
and 
$\prod _{\underset{\ell \not= k}{1 \leq \ell \leq M}} \frac{4N^2(\sigma ^{N,k}_\ell - \mu )} {\sigma^{-,k}_\ell - \lambda} = O(1)$
(Theorem \ref{Theorem 1.1}, $\sigma^N_\ell \in [\lambda^N_{2\ell -1}, \lambda^N_{2\ell}]$).
As $Q^{N,M}_k(\mu), \,Q^{-,M}_k(\lambda)^{-1} = O(1 )$ 
the claimed estimate then follows from  \eqref{fcor}.
\qed

\medskip

Let us now analyze the integrals over $A^N_j$ in the case $0 < j < n$, $j \ne k.$
\begin{lemma}
\label{l.n.2} Uniformly for $1 \le n, k \le L$, $0 < j < n$ with $j \ne k$, and on bounded
sets of functions $\alpha , \beta $ in $C_0^2(\mathbb T)$, $\int _{A_j^N} \frac{\varphi ^N_k(\mu )}{i \sqrt[c]{\chi_N(\mu )}} d\mu$ has the asymptotics
\begin{equation}
\label{da.i.2}
\frac{N}{2\pi k } \int_{A_j^-}
\frac{\psi^-_k(\lambda)}{i\sqrt[c]{\Delta^2_-(\lambda)-4}}d\lambda
+O\Big(\big(\frac{M}{N^{1/2}}+\frac{1}{L^{5/2}}\big)\frac{N} {j(k^2-j^2)} \Big).
\end{equation}      
\end{lemma}
%%%%%%%%%%%%%%%%%%%%%%%%%%%%%%%%%%%%%%%%%%%%%%%%%%%%%%%%%%%%%%%%%%%%%
\proof We present a proof which can be easily adapted to the case where $j = k$.
Similarly as in the proof of Lemma \ref{l.n.1}, write for $\mu \in A^N_j$
\[
 \frac{\varphi ^N_k(\mu )}{i \sqrt[c]{\chi_N(\mu )}}=
  \frac{1}{2\sqrt[+]{\mu - \lambda ^N_0}}
\Big(\frac{\sigma ^{N,k}_j - \mu}{w^N_k(\mu )}\prod _{\underset{\ell \not= j, k}{1 \leq \ell \leq L}} 
\frac{\sigma ^{N,k}_\ell - \mu } {w^N_\ell (\mu )} \Big) 
 Q^{N,L}_k(\mu) \frac{1}{w^N_j(\mu )}.
\]
Set $\mu \equiv \mu(\lambda) = -2 + \frac{\lambda}{4N^2}$ and write 
$\frac{\sigma ^{N,k}_j - \mu}{w^N_k(\mu )} 
\prod _{\underset{\ell \not= j,k}{1 \leq \ell \leq L}} \frac{\sigma ^{N,k}_\ell - \mu } {w^N_\ell (\mu )}$ as
\[
\frac{4N^2(\sigma^{N,k}_j-\mu)}{w^-_k(\lambda )} \prod _{\underset{\ell \not=j, k}{1 \leq \ell \leq L}} 
\frac{\sigma ^{-,k}_\ell - \lambda } {w^-_\ell (\lambda )}  \cdot
\prod _{\underset{\ell \not= j,k}{1 \leq \ell \leq L}} 
\frac{4N^2(\sigma ^{N,k}_\ell - \mu )} {\sigma^{-,k}_\ell - \lambda} \cdot
\prod _{\underset{\ell \not= j}{1 \leq \ell \leq L}}
\frac{w ^-_\ell(\lambda) } {4N^2w^N_\ell (\mu )}.
\]
Arguing as in the proof of Lemma \ref{l.n.1} and taking into account that 
by Theorem \ref{Theorem 1.1},
$4N^2(\mu - \lambda ^N_0) = (\lambda - \lambda ^-_0) (1 + O(\frac{1}{j^2}\frac{M^2}{N})),$
and hence $\frac{1}{2\sqrt[+]{\mu - \lambda ^N_0}} = 
\frac{N}{\sqrt[+]{\lambda - \lambda ^-_0}}(1 + O(\frac{1}{j^2} \frac{M^2}{N}))$,
 one sees that for 
$\lambda \in \tilde A^-_j := [\lambda^-_{2j-1} - 2 \rho, \lambda^-_{2j} + 2 \rho]$
% $[\frac{\lambda^-_{2j-1} + \lambda^-_{2j+2}}{2},  \frac{\lambda^-_{2j+1} +\lambda^-_{2j}}{2}]$
$$
\frac{\varphi ^N_k(\mu )}{i \sqrt[c]{\chi_N(\mu)}} =
f^{N,k}_j(\lambda)\frac{N}{w^N_j(\mu)}\big(1+ g^{N,k}_j(\lambda)\big) ,
$$
where
\begin{equation}
\label{da.i1.1}
f^{N,k}_j(\lambda)= \frac{1}{\sqrt[+]{\lambda-\lambda_0^-}}
\frac{4N^2(\sigma^{N,k}_j-\mu)}{w^-_k(\lambda )} 
\prod _{\underset{\ell \not=j, k}{0 < \ell < \infty }}
\frac{\sigma ^{-,k}_\ell - \lambda } {w^-_\ell (\lambda )}
\end{equation}
and $g^{N,k}_j$ is defined by the equation above and satisfies $g^{N,k}_j(\lambda) = O(\frac{1}{L^{5/2}})$.
Actually, $f^{N,k}_j$ and $g_j^{N,k}$ are defined and analytic 
on the rectangle $R_j \subset \mathbb C$  
with top and bottom side given by $\tilde A^-_j \pm i 2\rho$. 
Note that for $\lambda \in R_j$, one has
$4N^2(\sigma^{N,k}_j-\mu) = O(1)$  (Theorem \ref{Theorem 1.1}, 
$\sigma^{N,k}_j \in [\lambda^N_{2j-1}, \lambda^N_{2j}]$), 
$1/\sqrt[+]{\lambda-\lambda_0^-} = O(\frac{1}{j})$, $1/w^-_k(\lambda) =O(\frac{1}{k^2-j^2})$,
and $\prod _{\underset{\ell \not=j, k}{0 < \ell < \infty }}
\frac{\sigma ^{-,k}_\ell - \lambda } {w^-_\ell (\lambda )} = O(1).$
Hence for $\lambda \in R_j$
\begin{equation}\label{fcor2}
f^{N,k}_j(\lambda)=
O\big(\frac{1}{j(k^2-j^2)}\big).
\end{equation}
Moreover, by the same arguments as above, one sees that 
$g^{N,k}_j(\lambda) = O(\frac{1}{L^{5/2}})$
is valid on $R_j$.
To proceed further observe that, by an explicit computation, 
$$
\int_{A_j^N}\frac{d\mu}{w^N_j(\mu)}=0 ,
$$
so that the integral 
$\frac{1}{N}\int _{A_j^N} \frac{\varphi ^N_k(\mu )}{i \sqrt[c]{\chi_N(\mu )}} d\mu\,$
equals
\begin{align*}
\int _{A_j^N}
\frac{f^{N,k}_j(\lambda(\mu))(1+g^{N,k}_j(\lambda(\mu)))-f^{N,k}_j(\lambda(\tau^N_j))(1+g^{N,k}_j(\lambda(\tau_j^N)))}{w^N_j(\mu)}d\mu.
\end{align*}
Now the integrals over the intervals 
$[\lambda^N_{2j-1} - \rho^N, \lambda^N_{2j-1}]$ and
$[\lambda^N_{2j}, \lambda^N_{2j} + \rho^N]$ of $A_j^N$
can be analyzed separately. As the two cases are treated in the same way
we concentrate on the integral $(T)$ over the interval
$[\lambda^N_{2j},\lambda^N_{2j}+\rho^N]$ only. Write $(T) = (T1) - (T2) - (T3)$ where
with $\lambda \equiv \lambda(\mu) = 4N^2(\mu + 2),$ $\tilde \tau^N_j = 4N^2(\tau^N_j + 2)$,
and the identity 
$w_j^N(\mu) = - \frac{1}{2N}\sqrt[+]{\lambda-\nu^N_{2j-1}}\sqrt[+]{\mu-\lambda^N_{2j}} $
on $[\lambda^N_{2j}, \lambda^N_{2j} + \rho^N]$,
\begin{equation}
\label{da.i.3}
(T1) := \int_{\lambda^N_{2j}}^{\lambda^N_{2j}+\rho^N}
\frac{f^{N,k}_j(\lambda)-f^{N,k}_j(\tilde \tau^N_j)) }{w^N_j(\mu)}d\mu, \qquad
\end{equation}
\[
(T2) := \int_{\lambda^N_{2j}}^{\lambda^N_{2j}+\rho^N}f^{N,k}_j(\lambda)
\frac{g^{N,k}_j(\lambda)-g^{N,k}_j(\tilde \tau^N_j)}
{\sqrt[+]{\lambda-\nu^N_{2j-1}}}\frac{2Nd\mu}{\sqrt[+]{\mu-\lambda^N_{2j}}},
\]
\[
(T3) := g^{N,k}_j(\tilde \tau^N_j)\int_{\lambda^N_{2j}}^{\lambda^N_{2j}+\rho^N}
\frac{f^{N,k}_j(\lambda)-f^{N,k}_j(\tilde \tau^N_j)}
{\sqrt[+]{\lambda-\nu^N_{2j-1}}}\frac{2Nd\mu}{\sqrt[+]{\mu-\lambda^N_{2j}}}.
\]
By Lemma \ref{atif} there exists 
$\lambda^N_{2j} <\mu_* < \lambda^N_{2j-1}$ so that with $\lambda_* = \lambda(\mu_*)$
$(T2) = f^{N,k}_j(\lambda_*) \frac{g^{N,k}_j(\lambda_*)-g^{N,k}_j(\tilde \tau^N_j)}
{\sqrt[+]{\lambda_*-\nu^N_{2j-1}}}\, \frac{\sqrt[+]{\rho}}{2}. \,$
As  $f^{N,k}_j(\lambda_*) = O\big(\frac{1}{j(k^2-j^2)}\big)$ (use \eqref{fcor2}) and
by  Cauchy's theorem,
$g^{N,k}_j(\lambda_*)-g^{N,k}_j(\tilde \tau^N_j) 
= O\big(L^{-5/2}(\lambda_* - \tilde \tau^N_j)\big) $
(use $g^{N,k}_j(\lambda) = O\big(L^{-5/2})$)
it then follows from $\nu^N_{2j -1} \le \tilde \tau^N_j$ that
\[
(T2) = O\big(\frac{1}{j(j^2-k^2)}\frac{1}{L^{5/2}} \big).
\]
The term $(T3)$ is studied in a similar way and admits the same bound. 
Similarly, write 
$\int_{A_j^-} \frac{\psi^-_k(\lambda)}{\sqrt{\Delta^2_-(\lambda)-4}}d\lambda
={2\pi k} \int_{A_j^-}\frac{f^{-,k}_j(\lambda)}{w^-_j(\lambda)}d\lambda
$ 
where in view of \eqref{psidelta}, $f^{-,k}_j$ is given by \eqref{da.i1.1} with $4N^2(\sigma^{N,k}_j-\mu)$
replaced by $(\sigma^{-,k}_j-\lambda)$.  
As $\int_{A_j^-} \frac{1}{w^-_j(\lambda)}d\lambda = 0$
\[
\frac{1}{2\pi k} \int_{A_j^-} \frac{\psi^-_k(\lambda)}{\sqrt{\Delta^2_-(\lambda)-4}}d\lambda
=  \int_{A_j^-}\frac{f^{-,k}_j(\lambda)-f^{-k}_j(\tau^-_j)}{w^-_j(\lambda)} d\lambda,
\]
allowing again to consider the integrals over the intervals $[\lambda^-_{2j-1} - \rho, \lambda^-_{2j-1}]$ and $[\lambda^-_{2j}, \lambda^-_{2j} + \rho]$ of $A_j^-$
separately. It then remains to compare $(T1)$ with
\begin{equation}
\label{da.i.4}
(S1) := \int_{\lambda^-_{2j}}^{\lambda^-_{2j}+\rho}\frac{f^{-,k}_j(\lambda)-f^{-,k}_j
  (\tau^-_j)}{w^-_j(\lambda)} d\lambda.
\end{equation}
Make the change of variables
$
[0,\rho] \to [\lambda^N_{2j}, \lambda^N_{2j} +\rho^N], x \mapsto \mu(x) =\lambda^N_{2j}+4N^2x 
$
in $(T1)$ and $[0, \rho] \to [\lambda^-_{2j}, \lambda^-_{2j} + \rho], x \mapsto 
\lambda=\lambda^-_{2j}+x$ in $(S1)$ so that
with $\tilde \gamma_j^N=4N^2\gamma^N_j$, $\tilde\tau^N_j=4N^2(\tau^N_j+2)$, and
the sign of the standard root
\[
(T1) = - \int_0^\rho\frac{f^{N,k}_j(x+\nu^N_{2j})-f^{N,k}_j(\tilde\tau^N_j)}
{\sqrt[+] x\sqrt[+]{x+\tilde\gamma^N_j}}dx\\
=- \int_0^\rho\frac{F^{N,k}_j(x) \cdot (x+\tilde\gamma^N_j/2)}{\sqrt[+]{x+\tilde\gamma^N_j}}
\frac{dx}{\sqrt[+] x}
\]
where 
$$
F^{N,k}_j(x)=\int_0^1\partial_\lambda f^{N,k}_j(\tilde \tau^N_j+t(x+\tilde \gamma^N_j/2))dt .
$$
Similarly, with 
$F^{-,k}_j (x) = \int_0^1\partial_\lambda f^{-,k}_j( \tau^-_j+t(x+ \gamma^-_j/2))dt,$  
one has 
\begin{equation}
\label{S1}
(S1)= -
\int_0^\rho\frac{F^{-,k}_j(x) \cdot (x+\gamma^-_j/2)}{\sqrt[+]{x+\gamma^-_j}} \frac{dy}{\sqrt[+] x}.
\end{equation}
Thus, by Lemma \ref{atif}, there exists $0 < z < \rho$ so that
\[
\big( (S1) - (T1) \big) \frac{2}{\sqrt[+]\rho} \,\,=\,\,
F^{N,k}_j(z) \frac{z+\tilde\gamma^N_j/2}{\sqrt[+]{z+\tilde\gamma^N_j}}-
F^{-,k}_j(z)\frac{z+\gamma^-_j/2}{\sqrt[+]{z+\gamma^-_j}} \qquad
\]
\[
= 
\Big(F^{N,k}_j(z)-F^{-,k}_j(z)\Big)
  \frac{z+\tilde\gamma^N_j/2}{\sqrt[+]{z+\tilde\gamma^N_j}} +
F^{-,k}_j(z) \Big(\frac{z+\tilde\gamma^N_j/2}{\sqrt[+]{z+\tilde\gamma^N_j}}-
  \frac{z+\gamma^-_j/2}{\sqrt[+]{z+\gamma^-_j}} \Big).  
 \]
As  $F^{-,k}_j = O(\frac{1}{j(j^2-k^2)})$  (by Cauchy's estimate and \eqref{fcor2}),
$\tilde\gamma^N_j-\gamma^-_j =  O\big(\frac{M^2}{N})$
(by Theorem \ref{Theorem 1.1}), Lemma \ref{Lemma C.3bis} implies that 
\[
F^{-,k}_j(z) \Big(\frac{z+\tilde\gamma^N_j/2}{\sqrt[+]{z+\tilde\gamma^N_j}}-
  \frac{z+\gamma^-_j/2}{\sqrt[+]{z+\gamma^-_j}} \Big)
=O\big(\frac{1}{j(j^2-k^2)} \frac{M}{N^{1/2}} \big).
\] 
Using that by \eqref{fcor2} and Theorem \ref{theorem2.1}
%\[
%f^{N,k}_j(x+\nu^N_{2j}) -f^{N,k}_j(x+\lambda^-_{2j} ) 
%=  O\big( (\nu^N_{2j} - \lambda^-_{2j}) \frac{1}{j(k^2 - j^2)} \big)
%=O(\frac{M^2}{N} \frac{1}{j(k^2 - j^2)}),
%\]
\[
f^{N,k}_j(x+\lambda^-_{2j}) -f^{-,k}_j(x+\lambda^-_{2j} ) = 
O\big( \frac{\tilde \sigma^{N,k}_{j} - \sigma^{-,k}_{j}}{j(k^2 - j^2)} \big)
=O \big( \frac{1}{L^{5/2}}\frac{1}{j(k^2 - j^2)} \big),
\]
and using  that $f^{N,k}_j$ and $f^{-,k}_j$ are analytic, hence in particular Lipschitz,
one the concludes by Cauchy's estimate
and by the boundedness of  $\frac{x+\tilde\gamma^N_j/2}{\sqrt[+]{x+\tilde\gamma^N_j}}$
 that
\[\Big(F^{N,k}_j(z)-F^{-,k}_j(z)\Big)
  \frac{z+\tilde\gamma^N_j/2}{\sqrt[+]{z+\tilde\gamma^N_j}}
= O\big( \frac{1}{j(k^2 - j^2)} \frac{1}{L^{5/2}} \big).
\] 
By combining the estimates obtained the stated asymptotics follow. 
Going through the arguments of the proof one verifies that the claimed 
uniformity statement holds.\qed

%%%%%%%%%%%%%%%%%%%%%%%%%%%%%%%%%%%%%%%%%%%%%%%%%%%%%%%%%%%%%%%%%%%%%%%%
\begin{corollary}
\label{c.n.2} Uniformly for $1 \le n,k  \le M$, $0 < j < n$ with $j \ne k$, and on bounded
subsets of functions $\alpha , \beta $ in $C_0^2(\mathbb T)$
\begin{equation}
\label{da.i.21}
\int _{A_j^N} \frac{\varphi ^N_k(\mu )}{i \sqrt[c]{\chi_N(\mu )}} d\mu
=O\big(\frac{N}{j(k^2 - j^2)}\big).
\end{equation}      
\end{corollary}
%%%%%%%%%%%%%%%%%%%%%%%%%%%%%%%%%%%%%%%%%%%%%%%%%%%%%%%%%%%%%%%%%%%%%%%%
\proof  
One verifies in a straightforward way that uniformly for any $1 \le n,k \le M$ and
$0 < j < n$ with $j \ne k$, 
\eqref{fcor2} holds and the function $g^{N,k}_j,$ introduced in the proof of Lemma \ref{l.n.2}, 
satisfies $1 + g^{N,k}_j(\lambda(\mu(x))) = O(1)$.
Arguing as in the proofs of Corollary \ref{c.n.1} and  Lemma \ref{l.n.2}
yields the claimed estimate.
\qed

\medskip
Next we analyze the integral $\int _{A_j^N} \frac{\varphi ^N_k(\mu )}{i \sqrt[c]{\chi_N(\mu )}} d\mu$ in the case where $j=k.$
\begin{lemma}
\label{l.n.21} Uniformly for $1 \le n \le L,$ $0 < k < n$,  and on bounded
subsets of functions $\alpha , \beta $ in $C_0^2(\mathbb T)$
\begin{equation}
\label{da.i.22}
\int _{A_k^N} \frac{\varphi ^N_k(\mu )}{i \sqrt[c]{\chi_N(\mu )}} d\mu
=\frac{N}{2\pi k } \int_{A_k^-}
\frac{\psi^-_k(\lambda)}{i\sqrt[c]{\Delta^2_-(\lambda)-4}}d\lambda+
O\Big( \big(\frac{M}{N^{1/2}}+\frac{1}{L^{5/2}}\big) \frac{N} {k} \Big).
\end{equation}      
\end{lemma}
%%%%%%%%%%%%%%%%%%%%%%%%%%%%%%%%%%%%%%%%%%%%%%%%%%%%%%%%%%%%%%%%%%%%%%%%
\proof Proceeding as in the proof of Lemma \ref{l.n.2} one gets 
$$
\frac{\varphi ^N_k(\mu )}{i \sqrt[c]{\chi_N(\mu)}}=f^{N,k}_k(\lambda)
\frac{N}{w^N_k(\mu )}\Big(1+O\big(\frac{1}{L^{5/2}}\big)\Big)
$$
where now
\begin{equation}
\label{da.i1.11}
f^{N,k}_k(\lambda)= \Big( \prod _{\underset{\ell \not=k}{0 < \ell < \infty}}
 \frac{\sigma ^{-,k}_\ell - \lambda } {w^-_\ell (\lambda )}\Big) 
\frac{1}{\sqrt[+]{\lambda-\lambda_0^-}}\end{equation}
which is analytic in the rectangle $R_k$, introduced in the proof of Lemma \ref{l.n.2}.
Since $\lambda\sim  4\pi^2 k^2$ for $\lambda \in R_k$, it satisfies
\begin{equation}\label{fcor3}
f^{N,k}_k(\lambda)=
O\big(\frac{1}{k}\big).
\end{equation}
Arguing as in the proof of Lemma \ref{l.n.2}, the claimed statements follow. \qed

\begin{corollary}
\label{c.n.21} Uniformly for $1 \le n \le M,$ $1 < k < n$,  and on bounded
subsets of functions $\alpha , \beta $ in $C_0^2(\mathbb T)$
\begin{equation}
\label{da.i.23}
\int _{A_k^N} \frac{\varphi ^N_k(\mu )}{i \sqrt[c]{\chi_N(\mu )}} d\mu=O\big(\frac{N} {k}\big).
%{\tt check}
\end{equation}      
\end{corollary}
\proof
Arguing as in the proofs of Corollary \ref{c.n.2} and Lemma \ref{l.n.21},
the claimed statement follows.
\qed

\medskip
Next we analyse the integral $\int _{C_n^N} \frac{\varphi ^N_k(\mu )}{i \sqrt[c]{\chi_N(\mu )}} d\mu$ in the case where $k \ne n$.
%%%%%%%%%%%%%%%%%%%%%%%%%%%%%%%%%%%%%%%%%%%%%%%%%%%%%%%%%%%%%%%%%%%%
\begin{lemma}
\label{l.n.3} 
Uniformly for $1 \le n \le L$, $0 < k < n$, and on bounded
subsets of functions $\alpha , \beta $ in $C_0^2(\mathbb T)$
\begin{equation}
\label{da.i.10}
\int _{C_n^N} \frac{\varphi ^N_k(\mu )}{i \sqrt[c]{\chi_N(\mu )}} d\mu
=\frac{N}{2\pi k } \int_{C_n^-}
\frac{\psi^-_k(\lambda)}{i\sqrt[c]{\Delta^2_-(\lambda)-4 }}d\lambda
+O\big(\frac{N}{n(k^2-n^2)}( \frac{M}{N^{1/2}} + \frac{1}{L^{5/2}})\big).
\end{equation}
\end{lemma}
%%%%%%%%%%%%%%%%%%%%%%%%%%%%%%%%%%%%%%%%%%%%%%%%%%%%%%%%%%%%%%%%%%%%%%
\proof We proceed similarly as in the proof of Lemma \ref{l.n.2}, but have to  
take into account that possibly, $\lambda^N_{2n-1} = \lambda^N_{2n}$. Fortunately, in such a case
$\sigma ^{N,k}_n - \mu = \tau^N_n - \mu = w^N(\mu)$, implying that 
$\frac{\sigma ^{N,k}_n - \mu } {w^N_n (\mu )} \equiv 1$.
 Therefore write $\frac{\varphi ^N_k(\mu )}{i \sqrt[c]{\chi_N(\mu )}}$ as
\[
 \frac{\varphi ^N_k(\mu )}{i \sqrt[c]{\chi_N(\mu )}}=
  \frac{1}{2\sqrt[+]{\mu - \lambda ^N_0}}
\Big(\frac{1}{w^N_k(\mu )}\prod _{\underset{\ell \not= n, k}{1 \leq \ell \leq L}} 
\frac{\sigma ^{N,k}_\ell - \mu } {w^N_\ell (\mu )} \Big) 
 Q^{N,L}_k(\mu) \frac{\sigma ^{N,k}_n - \mu}{w^N_n(\mu )}.
\]
Arguing as in the proof of Lemma \ref{l.n.2}  one sees that for 
$\lambda$ in the rectangle $R_n \subset \mathbb C,$ as introduced there, and 
$\mu = -2 + \frac{\lambda}{4N^2}$
$$
\frac{\varphi ^N_k(\mu )}{i \sqrt[c]{\chi_N(\mu)}} = N
f^{-,k}_n(\lambda) \big(1+ g^{N,k}_n(\lambda)\big) \frac{4N^2(\sigma^{N,k}_n-\mu)}{w^N_n(\mu)} ,
$$
where
\begin{equation}
\label{bvc}
f^{-,k}_n(\lambda)= \frac{1}{\sqrt[+]{\lambda-\lambda_0^-}}
\frac{1}{w^-_k(\lambda )} 
\prod _{\underset{\ell \not=n, k}{0 < \ell < \infty }}
\frac{\sigma ^{-,k}_\ell - \lambda } {w^-_\ell (\lambda )}
= O \big( \frac{1}{n(k^2 - n^2)} \big)
\end{equation}
and $g^{N,k}_n$ is defined by the equation above and satisfies 
$g^{N,k}_n(\lambda) = O(\frac{1}{L^{5/2}})$.
Making the change of variable $\mu(x)=\lambda^N_{2n-1}-\rho^N+\frac{1}{4N^2}x $,
implying that $\lambda \equiv \lambda(\mu (x))=\nu^N_{2n-1}-\rho+x $,
and taking into account the sign of the standard root $w^N_n(\mu),$
the integral $\frac{1}{N}\int_{\lambda_{2n-1}^N-\rho^N}^{\lambda_{2n-1}^N}
\frac{\varphi ^N_k(\mu )}{i \sqrt[c]{\chi_N(\mu)}} d\mu$ equals
\begin{equation}
\label{da.i.30}
\int_0^\rho f^{-,k}_n(\lambda) \big(1+g^{N,k}_n(\lambda)\big)
\frac{\tilde\sigma^{N,k}_n-\nu_{2n-1}^N+\rho-x}{\sqrt[+]{\tilde\gamma^N_n+\rho-x}}
\frac{dx}{\sqrt[+]{\rho-x}}
\end{equation}
where as usual, $\tilde\sigma^{N,k}_n = 4N^2(\sigma^{N,k}_n +2)), $ 
$\nu^N_{2n-1}= 4N^2(\lambda^N_{2n-1} +2)$, and $\tilde \gamma^N_n = 4N^2 \gamma^N_n.$
As $0 \le \tilde\sigma^{N,k}_n-\nu_{2n-1}^N \le \tilde \gamma^N_n$ one has
$\frac{\tilde\sigma^{N,k}_n-\nu_{2n-1}^N+\rho-x}{\sqrt[+]{\tilde\gamma^N_n+\rho-x}} = O(1)$
yielding
\be
\label{fcor4}
\frac{1}{N}\int_{\lambda_{2n-1}^N-\rho^N}^{\lambda_{2n-1}^N}
\frac{\varphi ^N_k(\mu )}{i \sqrt[c]{\chi_N(\mu)}} d\mu 
=(T) + O \big( \frac{1}{L^{5/2}}  \frac{N}{n(k^2 - n^2)} \big)
\ee
\[
(T) := \int_0^\rho f^{-,k}_n(x + \nu^N_{2n-1}-\rho)
\frac{\tilde\sigma^{N,k}_n-\nu_{2n-1}^N+\rho-x}{\sqrt[+]{\tilde\gamma^N_n+\rho-x}}
\frac{dx}{\sqrt[+]{\rho-x}}
\]
Arguing as in the proof of Lemma \ref{l.n.2} it remains to estimate the difference $(T) - (S)$ where
\[
(S) := \int_0^\rho f^{-,k}_n(x + \lambda^-_{2n-1} - \rho)\frac{\sigma^{-,k}_n-\lambda_{2n-1}^-+\rho-x}
{\sqrt[+]{\gamma^-_n +\rho-x}}\frac{dx}{\sqrt[+]{\rho-x}}.
\]
By Lemma \ref{atif}, it suffices to estimate, uniformly for $0 \le x \le \rho,$
$$
f^{-,k}_n(x + \nu^N_{2n-1}-\rho)
\frac{\tilde\sigma^{N,k}_n-\nu_{2n-1}^N+\rho-x}{\sqrt[+]{\tilde\gamma^N_n+\rho-x}}
-
f^{-,k}_n(x + \lambda^-_{2n-1}-\rho)
\frac{\sigma^{-,k}_n-\lambda_{2n-1}^-+\rho-x}{\sqrt[+]{\gamma^-_n +\rho-x}}.$$ 
Note that by Cauchy's theorem, 
$f^{-,k}_n(x + \nu^N_{2n-1}-\rho) - f^{-,k}_n(x + \lambda^-_{2n-1}-\rho)$ is
\[
O\big( (\nu^N_{2n-1} - \lambda^N_{2n-1}) \frac{1}{n(k^2 - n^2)} \big)
= O\big( \frac{M^2}{N}\frac{1}{n(k^2 - n^2)} \big)
\]
whereas 
\[
(\tilde\sigma^{N,k}_n-\nu_{2n-1}^N)
-(\sigma^{-,k}_n-\lambda_{2n-1}^-) = O\big( \frac{1}{L^{5/2}} + \frac{M^2}{N} \big)
= O\big(\frac{1}{L^{5/2}} \big)
\]
 and, by Lemma \ref{Lemma C.3bis}, 
\[
\big|\frac{\rho-x}{\sqrt[+]{\tilde\gamma^N_n+\rho-x}} -
\frac{\rho-x}{\sqrt[+]{\gamma^-_n +\rho-x}} \big|
\le  \sqrt[+]{\big|\tilde\gamma^N_n - \gamma^-_n \big|} = O\big( \frac{M}{N^{1/2}}\big).
\]
Combining these estimates, the stated asymptotics follow.
Going through the arguments of the proof the claimed 
uniformity statement follows.\qed

\begin{corollary}
\label{c.n.3} 
Uniformly for $1 \le n \le M$, $0 < k < n$, and on bounded
subsets of functions $\alpha , \beta $ in $C_0^2(\mathbb T)$
\begin{equation}
\label{da.i.10c}
\int _{C_n^N} \frac{\varphi ^N_k(\mu )}{i \sqrt[c]{\chi_N(\mu )}} d\mu
=O\big(\frac{N}{n(k^2-n^2)}\big).
\end{equation}
\end{corollary}
\proof
One verifies in a straightforward way that uniformly for any $1 \le k,n \le M$ with $k \ne n$, 
\eqref{bvc} holds and the function $g^{N,k}_n,$ introduced in the proof of Lemma \ref{l.n.3}, 
satisfies $1 + g^{N,k}_n(\lambda(\mu(x))) = O(1)$. Following the arguments of the proof
of Lemma \ref{l.n.3}, the claimed statement then follows.
\qed

\medskip
%%%%%%%%%%%%%%%%%%%%%%%%%%%%%%%%%%%%%%%%%%%%%%%%%%%%%%%%%%%%%%%%%%%%%%%%%%
Next we analyze the integral $\int _{C_n^N} \frac{\varphi ^N_k(\mu )}{i \sqrt[c]{\chi_N(\mu )}} d\mu$ in the case where $k=n$.
\begin{lemma}
\label{l.n.4} Uniformly for $1 \le n \le L$ and on bounded
subsets of functions $\alpha , \beta $ in $C_0^2(\mathbb T)$
\begin{equation}
\label{da.i.10d}
\tilde \gamma^N_n\int _{C_n^N} \frac{\varphi ^N_n(\mu )}{i \sqrt[c]{\chi_N(\mu )}} d\mu
=\frac{N\gamma_n^-}{2\pi n }  \int_{C_n^-}
\frac{\psi^-_n(\lambda)}{i\sqrt[c]{\Delta^2_-(\lambda)-4 }}d\lambda
+O\big(\frac{N}{n} ( \frac{M}{N^{1/2}} + \frac{1}{L^{5/2}}) \big).
\end{equation}
\end{lemma}
%%%%%%%%%%%%%%%%%%%%%%%%%%%%%%%%%%%%%%%%%%%%%%%%%%%%%%%%%%%%%%%%%%%%%%%%%%
\proof We proceed similarly as in the proof of Lemma \ref{l.n.3}, but note that in the case
$k=n,$ the factor $\gamma^N_n$ is a substitute for the missing factor $\sigma^{N,k}_n - \mu$.
Using the same terminology as in the proof of Lemma \ref{l.n.3}, we have
\begin{equation}
\label{da.i.101}
\frac{\tilde\gamma^N_n}{N} 
\int_{\lambda_{2n-1}^N-\rho^N}^{\lambda_{2n-1}^N}\frac{\varphi ^N_n(\mu )}
{i \sqrt[c]{\chi_N(\mu )}}=\int_0^\rho
f^{-,n}_n(\lambda) \big(1+g^{N,k}_n(\lambda)\big) 
\frac{\tilde \gamma^N_n}{\sqrt[+]{\tilde\gamma^N_n+\rho-x}}
\frac{dx}{\sqrt[+]{\rho-x}}
\end{equation}
where in the case at hand, $f^{-,n}_n,$ defined on the rectangle $R_n$, satisfies
\begin{equation}
\label{fcor5}
f^{-,n}_n(\lambda)= \frac{1}{\sqrt[+]{\lambda-\lambda_0^-}}
 \prod _{\underset{\ell \not =n}{0 < \ell < \infty}} 
\frac{\sigma ^{-,n}_\ell - \lambda } {w^-_\ell (\lambda )}
= O\big( \frac{1}{n} \big)
\end{equation}
and $g^{N,k}_n(\lambda)$ is again $O(\frac{1}{L^{5/2}})$.
Following the arguments of  the proof of Lemma \ref{l.n.3} yields the claimed statements.
\qed

\begin{corollary}
\label{c.n.4} Uniformly for $1 \le n \le M$ and on bounded
subsets of functions $\alpha , \beta $ in $C_0^2(\mathbb T)$
\begin{equation}
\label{da.c.10}
\tilde \gamma^N_n\int _{C_n^N} \frac{\varphi ^N_n(\mu )}{i \sqrt[c]{\chi_N(\mu )}} d\mu
=O\big(\frac{N}{n }\big) .
\end{equation}
\end{corollary}
\proof
 Arguing as in the proofs of Corollary \ref{c.n.3} and Lemma \ref{l.n.4}  yields the claimed
statement.
\qed

\medskip
Finally we analyze the integrals over the intervals $B^N_j$ with $1 \le j \le n-1$.
%%%%%%%%%%%%%%%%%%%%%%%%%%%%%%%%%%%%%%%%%%%%%%%%%%%%%%%%%%%%%%%%%%%%%%
\begin{lemma}
\label{l.n.5} 
Uniformly for $1 \le n, k \le L,$ $0< j <  n$, and on bounded
subsets of functions $\alpha , \beta $ in $C_0^2(\mathbb T)$
\begin{equation}
\label{da.i.1012}
\int _{B_j^N} \frac{\varphi ^N_k(\mu )}{i \sqrt[c]{\chi_N(\mu )}} d\mu
=\frac{N}{2\pi k } \int_{B_j^-}
\frac{\psi^-_k(\lambda)}{i\sqrt[c]{\Delta^2_-(\lambda)-4 }}d\lambda
+O\big(\frac{N}{j \cdot d_{k,j}}\frac{1}{L^{5/2}}\big)
\end{equation}
where $d_{k,j} = 1 + \mbox{min}(|j^2 -k^2|, |k^2 - (j+1)^2|) $.
\end{lemma}
%%%%%%%%%%%%%%%%%%%%%%%%%%%%%%%%%%%%%%%%%%%%%%%%%%%%%%%%%%%%%%%%%%%%%
\proof As on the intervals $B^N_j,$ the integrand 
$\frac{\varphi ^N_k(\mu )}{i \sqrt[c]{\chi_N(\mu )}}$ is not
singular, the integral is easier to handle. To get the claimed results, one argues 
as in the proof of Lemma \ref{l.n.1},
taking into account that $\frac{1}{w^-_k(\lambda)} = O(\frac{1}{d_{k,j}})$.
\qed

\medskip

Following the arguments of the proofs of Lemma \ref{l.n.5} and Corollary \ref{c.n.1} leads to the following
\begin{corollary}
\label{c.n.5} Uniformly for $1 \le n, k \le M,$ $0< j < n$, and on bounded
subsets of functions $\alpha , \beta $ in $C_0^2(\mathbb T)$
\begin{equation}
\label{da.c.1012}
\int _{B_j^N} \frac{\varphi ^N_k(\mu )}{i \sqrt[c]{\chi_N(\mu )}} d\mu
=O\big(\frac{N}{j \cdot d_{j,k}}\big).
\end{equation}
\end{corollary}

\medskip

%%%%%%%%%%%%%%%%%%%%%%%%%%%%%%%%%%%%%%%%%%%%%%%%%%%%%%%%%%%%%%%%%%%%%%%%%%%
%%%%%%%%%%%%%%%%%%%%%%%%%%%%%%%%%%%%%%%%%%%%%%%%%%%%%%%%%%%%%%%%%%%%%%%%%%%
%%%%%%%%%%%%%%%%%%%%%%%%%%%%%%%%%%%%%%%%%%%%%%%%%%%%%%%%%%%%%%%%%%%%%%%%%%%
\section{Symmetry of the Toda chain}
\label{TodaSymmetry}

In this Appendix we discuss a symmetry of the Toda chain used
to reduce the proof of the claimed asymptotics \eqref{freq2inSec6} 
of Theorem \ref{thm1.4} of the frequencies $\omega^N_n$ 
at the right edge to the one of the asymptotics \eqref{freq1inSec6} of the frequencies 
at the left edge. It is more convenient to discuss the symmetry property for arbitrary Toda chains.
Extend a vector $(b,a) = ((b_n)_{1 \le n \le N}, (a_n)_{1 \le n \le N})$ in 
$\mathbb R^N \times \mathbb R^N_{>0}$
and extend the vectors $b$ and $a$ periodically, $(b_n)_{n \in \mathbb Z},$ respectively
$(a_n)_{n \in \mathbb Z}.$  Define for any $n \in \mathbb Z$
\[
\tilde b_n:=-b_{N-n}   \quad \tilde a_n:=a_{N-1-n}
\]
 and consider the $2N \times 2N$ matrix $Q(\tilde b,\tilde a)$.
In this appendix, the superindex $N$ is used to indicate the number of particles of the Toda chain considered and we use tilde for quantities such as eigenvalues, characteristic functions, discriminants, actions, \dots
when evaluated at $(\tilde b, \tilde a)$ instead of $(b,a).$ So e.g. 
$\tilde \lambda^N_n = \lambda^N_n (\tilde b, \tilde a)$ or 
$\tilde \Delta_N(\mu) = \Delta_N(\mu, (\tilde b, \tilde a)).$
Then the following proposition holds.

\begin{proposition}
\label{symmetry.main}
The following formulae hold:
\[
(i)\,\,\,\, \tilde\lambda^N_n=-\lambda^N_{2N-1-n}\quad \forall \,\, 0 \le n < 2N \quad
\mbox{and hence} \quad\tilde \chi_N(\mu)=\chi_N(-\mu). \quad \qquad
\]
\[
(ii) \,\quad \tilde \Delta_N(\mu)=(-1)^N\Delta_N(-\mu);
\quad \partial_{\mu} \tilde \Delta_N (\mu)=-(-1)^N\dot\Delta_N(-\mu).\qquad \qquad \qquad
\]
\[
(iii) \,\,\quad \tilde I^N_n= I^N_{N-n}\quad \mbox{and} \quad
\tilde \varphi^N_n(\mu)= (-1)^N \varphi^N_{N-n}(-\mu)\quad \forall 0 < n < N.\qquad \qquad \qquad
\]
Combining these formulas one obtains 
\begin{equation}
\label{sym7}
 \tilde \omega^N_n=\omega^N_{N-n} \quad \forall  \,\, 0 < n < N\ .
\end{equation}
\end{proposition}
\proof To prove (i) note that an arbitrary vector $(F_n)_{1 \le n \le 2N} \in \mathbb R^{2N}$ is an eigenvector of $Q(b,a)$ with eigenvalue $\lambda$ iff $((-1)^nF_{2N-n})_{1 \le n \le 2N}$ is an eigenvector of  $Q(\tilde b,\tilde a)$ with eigenvalue $-\lambda$. The identity for the charcteristic polynomial 
$\tilde \chi_N (\mu)$ then follows from its product representation. Similarly, the first identity in (ii) 
follows from the product representation of $\Delta_N (\mu) - 2$, implying the claimed identity for the derivatives. Towards (iii) note that 
by the definition \eqref{cNroot} of the c-root, one has
$\sqrt[c]{\Delta^2 _N(-\mu) - 4}  = (-1)^N \sqrt[c]{\Delta^2 _N(\mu) - 4}.$
The stated identity for the actions then  follows from formula \eqref{Arnold} 
whereas the one for $\tilde \varphi^N_n$ follows from \eqref{boot} and the fact that they are
polynomials with leading term $\mu^{N-2}$. Finally, the claimed formula for the frequencies
then follows from \eqref{B.1}. \qed

\medskip

By applying the identity \eqref{B.5} to $\tilde\omega_{n}^N$, Proposition \ref{symmetry.main} 
leads to the following

\begin{corollary}
\label{cor.sym}
For any $0 < n < N$
\[
iN\omega^N_{N-n}= \sum ^n_{j=1} \int ^{\lambda ^N_{2N-1-(2j-2)}}
                    _{\lambda ^N_{2N-1-(2j-1)}} \frac{(\mu - \frak p_N/N)\dot \Delta _N
                     (\mu )}{\sqrt[c]{\Delta^2_N(\mu ) - 4}} d\mu 
\]
\begin{equation}
\label{sym.8}
                - \sum _{k \in {\cal J_N}} I^N_{N-k} \omega ^N_{N-k} \sum ^n_{j=1}
                    \int ^{\lambda ^N_{2N-1-(2j-2)}}_{\lambda ^N_{2N-1-(2j-1)}}
                     \frac{\varphi ^N_{N-k}(\mu )}{\sqrt[c]{\chi_N(\mu )}} d\mu
\end{equation}
\end{corollary}

%%%%%%%%%%%%%%%%%%%%%%%%%%%%%%%%%%%%%%%%%%%%%%%%%%%%%%%%%%%%%%%%%%%%%%%%%%%%
%%%%%%%%%%%%%%%%%%%%%%%%%%%%%%%%%%%%%%%%%%%%%%%%%%%%%%%%%%%%%%%%%%%%%%%%%%%%
%%%%%%%%%%%%%%%%%%%%%%%%%%%%%%%%%%%%%%%%%%%%%%%%%%%%%%%%%%%%%%%%%%%%%%%%%%%%

\medskip

%%%%%%%%%%%%%%%%%%%%%%%%%%%%%%%%%%%%%%%%%%%%%%%%%%%%%%%%%%%%%%%%%%%%%%%%%%
%%%%%%%%%%%%%%%%%%%%%%%%%%%%%%%%%%%%%%%%%%%%%%%%%%%%%%%%%%%%%%%%%%%%%%%%%%%
%%%%%%%%%%%%%%%%%%%%%%%%%%%%%%%%%%%%%%%%%%%%%%%%%%%%%%%%%%%%%%%%%%%%%%%%%%%

\end{document}